
	%
	%
	%
	\newfam\symboles
	

%
 %

 %

	\font\ninerm=cmr9
\font\eightrm=cmr8
 \font\ninei=cmti9

%


	%
	%

\def\flead{\leaders\hbox to 5pt {\hss.}\hfill}
\def\somr#1|#2|{#1\flead \rlap{\hbox to 25pt{\hfill #2}}\par
\goodbreak}
\def\ssomr#1|#2|{\qquad #1\flead \rlap{\hbox to 25pt {\hfill #2}}\par
\goodbreak}

 \magnification=1200 
\pageno=1
	\baselineskip7mm
	\hsize=130mm \vsize=200mm 
	\hoffset=35mm \voffset=10mm
	\hfuzz=50pt \vfuzz=10pt
	
	\headline={\hfill-- \tenrm\folio\hbox{ }--\hfill}
	\footline={\hfill}

\newtoks\auteurcourant      \auteurcourant={\hfil}
\newtoks\titrecourant       \titrecourant={\hfil}

\newtoks\hautpagetitre      \hautpagetitre={\hfil}
\newtoks\baspagetitre       \baspagetitre={\hfil}

\newtoks\hautpagegauche   \newtoks\hautpagedroite 
  
\hautpagegauche={\tenrm\rlap{\folio}\tenit\hfil\the\auteurcourant\hfil}
\hautpagedroite={\tenit\hfil\the\titrecourant\hfil\tenrm\llap{\folio}}

\newtoks\baspagegauche      \baspagegauche={\hfil} 
\newtoks\baspagedroite      \baspagedroite={\hfil}

\newif\ifpagetitre          \pagetitretrue
\nopagenumbers  

\def\nopagenumbers{\def\folio{\hfil}}  

\headline={\ifpagetitre\the\hautpagetitre
            \else\ifodd\pageno\the\hautpagedroite
             \else\the\hautpagegauche
              \fi\fi}

\footline={\ifpagetitre\the\baspagetitre\else
            \ifodd\pageno\the\baspagedroite
             \else\the\baspagegauche
              \fi\fi
               \global\pagetitrefalse}


\def\raggedbottom{\topskip 10pt plus 36pt\r@ggedbottomtrue}

	%
	%
	%
	%
	%
	%
	
	%
        \def\page#1{\leaders\hbox to 5mm{\hfil\hbox{\ninerm.}\hfil}
        \hfill\rlap{\hbox to 5mm{\hfill#1}}\par}

%
  \def\Th#1{\medskip\goodbreak \nod{\bf  \underbar{Th\'eor\`eme} #1.}
	\quad\nobreak \sl }
%
  

%
\def\Prop#1{\medskip\goodbreak \nod{\bf
	 Proposition #1.} \quad\nobreak \sl }
	
%
	\def\Lemme#1{\medskip\goodbreak \nod{\bf  Lemme #1.}
	\quad\nobreak \sl }
%

%

%

%

%

%
	
%

%
        
%

%
\def\Cor#1{\medskip\goodbreak \nod{\bf
	Corollaire#1.} \quad\nobreak \sl }
%
	\def\dem{\medskip\goodbreak \nod \ninei D\'emonstration. \rm}

  \def\ds{\displaystyle}
	\def\ss{\scriptstyle}
	
	\def\nod{\noindent}
	
	\def\hb{\hbox}

	\def\dbar{d\!\!\hbox to 4.5pt{\hfill\vrule height 5.5pt
	depth -5.3pt width 3.5pt}}

	\def\tas#1_#2^#3{\mathrel{\mathop{\kern 0pt#1}
					\limits_{#2}^{#3}  }}

	\def\up#1{\raise \hbox{\ss \rm #1}}


	\def\hfl#1#2{\smash{\mathop{\hbox to 10mm{\rightarrowfill}}
	\limits^{\ss#1}_{\ss#2}}}
        \def\hfla#1#2{\smash{\mathop{\hbox to 6mm{\rightarrowfill}}
	\limits^{\ss#1}_{\ss#2}}}
	\def\hflo#1#2{\smash{\mathop{\hbox to 10mm{\leftarrowfill}}
	\limits^{\ss#1}_{\ss#2}   }}

	\def\pile#1#2{\smash{\mathop{\hbox to 20mm{}}
	\limits^{\hb{\rm#1}}_{\hb{\rm#2}} }}
	\def\pila#1#2#3{\smash{\mathop{\hbox to 20mm{#2}}
	\limits^{\hb{\rm#1}}_{\hb{\rm#3}}   }}

 \auteurcourant={ M. Kaddar }
 \pageno=1
\baselineskip=5mm

 \hoffset=-4mm
  \voffset=8mm
 \hsize=1cm
\vsize=180mm \hfuzz=-2cm
\input amssym.def
\input amssym.tex%
\magnification = 1200 \hoffset=-4mm
  \voffset=8mm
 \hsize=14cm
\vsize=190mm
\input xy
\xyoption{all}
\par
\topskip = 10 truemm \overfullrule=0pt \rightskip=10mm

\def \cqfd{\unskip\kern 6pt\penalty 500
 \raise -2pt\hbox{\vrule\vbox to10pt{\hrule width 4pt
\vfill\hrule}\vrule}\par}

\def\adots{\mathinner{\mkern2mu\raise1pt\hbox{.}
\mkern3mu\raise4pt\hbox{.}\mkern1mu\raise7pt\hbox{.}}}
\def\pmb#1{\setbox0=\hbox{#1}%
\kern-.025em\copy0\kern-\wd0 \kern.05em\copy0\kern-\wd0
\kern-.025em\raise.0433em\box0}

 \font \tite=cmbxti10 scaled 1400


\def\build#1_#2^#3{\mathrel{\mathop{\kern 0pt#1}\limits_{#2}^{#3}}}

\def\build#1_#2^#3{\mathrel{\mathop{\kern 0pt#1}\limits_{#2}^{#3}}}
\def\hfl#1#2{\smash{\mathop{\hbox to 18mm{\rightarrowfill}}\limits^{\scriptstyle#1}_{
\scriptstyle#2}}}

\def\boxit#1#2{\setbox1=\hbox{\kern#1{#2}\kern#1}%
\dimen1=\ht1 \advance\dimen1 by #1 \dimen2=\dp1 \advance\dimen2 by
#1
\setbox1=\hbox{\vrule height\dimen1 depth\dimen2\box1\vrule}%
\setbox1=\vbox{\hrule\box1\hrule}%
\advance\dimen1 by .4pt \ht1=\dimen1 \advance\dimen2 by .4pt
\dp1=\dimen2 \box1\relax}
\def\page#1{\leaders\hbox to 5mm{\hfil.\hfil}\hfill\rlap{\hbox to
10mm{\hfill#1}}\par} \rightskip=10mm\sl

\titrecourant={Morphismes g\'eom\'etriquement plats  et classes
fondamentales relatives en g\'eom\'etrie analytique complexe.}

\input xy

\headline={\ifnum\pageno=1 {\hfill} \else{\hss \tenrm -- \folio\ --
\hss} \fi} \footline={\hfil}



\input xy
\xyoption{all}
\par

\par\centerline{{\tite{PLATITUDE GEOMETRIQUE ET }}}\bigskip
\centerline{\tite{ CLASSES FONDAMENTALES
RELATIVES}}\bigskip\noindent\centerline{\tite{ PONDEREES
I.}}\bigskip\bigskip
\centerline{MOHAMED KADDAR}\bigskip\bigskip\bigskip\noindent
$${\vbox{\hsize=12cm{{\bf Abstract.~} {\sevenrm{Let $X$ and $S$ be complex
spaces with $X$ countable at infinity and $S$ locally finite reduced
 pure dimensional. Let $\pi:X\rightarrow S$ be an open  morphism with constant
 pure $n$-dimensional fibers (we call such morphism
 universally-$n$-equidimensional). If there is a cycle $\goth{X}$
  of $X\times S$ such that, his support coincide fiberwise set-theorically
  with the fibers  of $\pi$ and endowed this with a good
  multiplicities in  such a  way that $(\pi^{-1}(s))_{s\in S}$ becomes
  a local analytic (resp. continuous) family of cycles in the sense of [B.M],
  we say that $\pi$ is analytically geometrically flat (resp. continuously geometrically flat) according to the weight  $\goth{X}$. The purpose of this article is to give a  sheaf-theoritic  characterization of  analytically  geometrically flat maps. One of many results obtained in this work  say that
an  universally-$n$-equidimensional morphism is analytically geometrically flat if and only if admit a
weighted relative fundamental  class morphism satisfies many nice functorial properties which giving, for a finite Tor-dimensional morphism or in the embedding case, the  relative fundamental class of Angeniol-Elzein  [E.A] or Barlet [B4]. From this, we deduce that universally-$n$-equidimensional morphism are  analytically geometrically flat if and only if there exist a  Kunz-Waldi sheaf of regular meromorphic
  relative forms which is  compatible with arbitrary  base change between reduced complex spaces. We have the  algebraic statement :\par\noindent
Let $X$ and  $S$ be reduced (or more generally with no embedded points)
  locally noetherian schemes of finite Krull dimension  with $S$ excellent
  on field  ${\rm k}$ of  characteristic 0 . Let $\pi:X\rightarrow S$ be
  a  finite type, universally open with $n$- pure dimensional fibers ( generically smooth)  morphism. Then the Kunz- Waldi relative sheaf of meromorphic regular
   forms, $\tilde{\omega}^{n}_{X/S}$,  is  compatible with any base change
   between ${\rm k}$-schemes like  $S$ if and only if $\pi$
   define an  algebraic family  of  $n$-cycles parametrized by
   $S$.}}}}}$$\rm\bigskip\bigskip\noindent
AMS Classification (2000): 32C15, 32C30, 32C35, 32C37, 14C05.\bigskip\noindent

\vfill\eject\noindent \centerline{{\tite Table des
mati\`eres:}}\bigskip\bigskip\noindent
 \S 0. {\bf
Introduction et \'enonc\'e des r\'esultats}\page {3} \noindent\S 1.
{\bf Notion de morphisme contin\^ument et
 analytiquement g\'eom\'etriquement plat}\par
 \leftskip=1em
{\bf 1.0. Terminologie relative \`a l'espace des cycles.}\page
{13}\smallskip{\bf 1.1. Quelques  exemples.}\page {27} {\bf 1.2. Sur quelques
notions fondamentales.}\page {30} {\bf 1.3. Platitude
g\'eom\'etrique.}\page{36} {\bf 1.4. Quelques petits
r\'esultats.}\page {47} {\bf 1.5. Analogie avec la platitude
alg\'ebrique.}\page {50}
\smallskip\noindent\S 2. {\bf Morphisme g\'eom\'etriquement plat et
int\'egration}\par\leftskip=1em {\bf 2.0. Position du probl\`eme et
th\'eor\`eme de Reiffen.}\page {51} {\bf 2.1. Classes de cohomologie
et repr\'esentant $\bar\partial$-ferm\'e sur un espace \indent
singulier}\page {52} {\bf 2.2. Morphisme universellement ouvert
adapt\'e \`a l'int\'egration.}\page {54} {\bf 2.3. Propri\'et\'es
fonctorielles de l'int\'egration sur les fibres d'un mor\par\indent
phisme g\'eom\'etriquement plat.}\page {57}\smallskip\noindent
 \S 3. {\bf
Faisceaux dualisants, formes r\'eguli\`eres et
 m\'eromorphes r\'eguli\`eres: cas absolu et relatif.}\par\leftskip=1em {\bf 3.0. Faisceaux dualisants.}\page {60} {\bf 3.1.
Formes r\'eguli\`eres et m\'eromorphes r\'eguli\`eres: cadre
alg\'ebrique.}\page {61}
 {\bf 3.2. Formes r\'eguli\`eres et
m\'eromorphes r\'eguli\`eres: cadre analytique.}\page
{65}\smallskip\noindent \S 4. {\bf Classe fondamentale: cas absolu
et relatif.}\par\leftskip=1em
 {\bf 4.0. Cas absolu.}\page {77}
{\bf 4.1. Cas relatif.}\page {80}\smallskip\noindent \S 5. { \bf
Quelques petits r\'esultats utiles.}\page {84} \smallskip\noindent
\S 6.  {\bf
Bibliographie}\page {86}\bigskip\bigskip\bigskip\bigskip\eject\vfill
\noindent\centerline{{\bf $\S 0.$ \tite{INTRODUCTION}.}}
\smallskip\bigskip\smallskip\noindent {\bf 0.0.}  Soient $n$ un
entier et $S$ un espace analytique complexe r\'eduit de dimension
pure. Notons ${\cal E}(S,n)$ (resp. ${\cal G}_{a}(S,n)$) l'ensemble
des morphismes $\pi:X\rightarrow S$ ouverts et \`a fibres de
dimension pure $n$ avec $X$ d\'enombrable \`a l'infini (resp. les
morphismes analytiquement g\'eom\'etriquement plats c'est-\`a-dire
les \'el\'ements de ${\cal E}(S,n)$ munis  d' un cycle poids de
$S\times X$  induisant des multiplicit\'es convenables sur les
fibres pour en faire une famille analytique (locale) de cycles au
sens de [B1],[B.M])\footnote{$^{(0)}$}{ Les \'el\'ements de ${\cal
E}(S,n)$ seront appel\'es par la suite {\it universellement
$n$-\'equidimensionnels}. Nous insistons sur le terme {\it
analytiquement} car, dans ce qui suit, se d\'egagera naturellement
la notion de {\it contin\^ument g\'eom\'etriquement plat} dont la
classe sera not\'ee  ${\cal G}_{c}(S,n)$.}. ${\cal G}_{a}(S,n)$
forme une classe de morphismes particuli\`erement int\'eressante
puisqu'elle contient les morphismes plats et plus g\'en\'eralement
ceux de
 Tor-dimension finie, elle est  stable par changement de base entre espaces
 complexes r\'eduits et de nature locale sur $X$ et $S$. De plus dans toute factorisation
 locale de $\pi$ (par rapport \`a
 l'une quelconque de ses fibres) en un morphisme fini et surjectif suivi d'une projection lisse
  sur la base, le morphisme fini est toujours analytiquement
  g\'eom\'etriquement plat. Remarquons au passage que cette propri\'et\'e est
  v\'erifi\'ee par les morphismes de Tor-dimension finie, mais pas  par les
 plats car sinon  les fibres seraient toujours de Cohen-Macaulay !
 Signalons que cette classe
n'est toutefois  pas stable pour la composition des morphismes et
rend, de ce fait, difficile l'\'elaboration d'une th\'eorie de
l'intersection avec param\`etre g\'en\'eral.\smallskip\noindent
L'objectif premier de
 cet article est de donner une caract\'erisation simple  des \'el\'ements de  ${\cal G}_{a}(S,n)$ qui va  s'exprimer aussi bien en terme de morphisme  d'{\it int\'egration}  (op\'eration de nature globale)  ou en terme  de morphismes  {\it trace} (op\'eration de nature purement locale). Tout ceci passe par la mise en \'evidence de deux  faisceaux particuli\`erement importants, l'un pouvant apparaitre comme l'avatar local de l'autre.\par\noindent
Pour cela, on montre, dans un premier temps (cf th\'eor\`eme {\bf
1}), qu'\`a tout \'el\'ement $\pi$ de ${\cal E}(S,n)$ est
canoniquement associ\'e  un faisceau ${\cal O}_{X}$-coh\'erent,
$\omega^{n}_{\pi}$, de
 profondeur au moins deux fibre par fibre sur $S$, munissant $\pi$  d'une fl\`eche
 myst\'erieuse
  ${\rm I}\!{\rm R}^{n}\pi_{!}\omega^{n}_{\pi}\rightarrow {\cal O}_{S}$ de formation compatible
  aux restrictions ouvertes sur $X$ et $S$, stable par changement de
  base plat et coincidant, pour $\pi$ propre, gr\^ace \`a la dualit\'e analytique
  relative de Ramis-Ruget-Verdier ([R.R.V]) et \`a l'isomorphisme de
  Verdier ([Ve]), avec ${\cal
H}^{-n}(\pi^{!}({\cal O}_{S}))$ qui est le $n$-\`eme faisceau
d'homologie du complexe \`a cohomologie coh\'erente $\pi^{!}({\cal
O}_{S})$.\par\noindent On prouve, alors (cf th\'eor\`eme {\bf 2}),
 qu'un \'el\'ement $\pi$ de  ${\cal E}(S,n)$ muni d'un certain cycle
poids $\goth{X}$ est dans ${\cal G}_{a}(S,n)$ si et seulement si
l'une des conditions \'equivalentes
  suivantes est r\'ealis\'ee  :\smallskip\noindent
 {{\bf(1)} il existe un unique morphisme d'int\'egration
 $\int_{\pi,{\goth{X}}}:{\rm I}\!{\rm R}^{n}\pi_{!}\Omega^{n}_{X/S}\rightarrow
 {\cal O}_{S}$, continu
 relativement aux structures {\bf Q.F.S} et {\bf F.S} dont on peut munir ces groupes
  de cohomologies, compatible \`a l'additivit\'e des pond\'erations $\goth{X}$, stable par
 changement de base entre espaces complexes  r\'eduits et donnant, en particulier,
  l'int\'egration usuelle dans
 le cas absolu.\par\noindent
 {\bf (2)} il existe un  unique morphisme de faisceaux coh\'erents
 ${\cal C}_{\pi,\goth{X}}:\Omega^{n}_{X/S}\rightarrow \omega^{n}_{\pi}$ v\'erifiant les propri\'et\'es suivantes:\par
{\bf(i)} il est de formation compatible \`a l'additivit\'e des
pond\'erations,   aux
 changement de base r\'eduit et satisfait  {\it{la propri\'et\'e de
 la
 trace relative}},\par
{\bf(ii)} il induit  un morphisme de complexes  diff\'erentiels
gradu\'es\par
 \centerline{${\cal C}^{\bullet}_{\pi,\goth{X}}:\Omega^{\bullet}_{X/S}\rightarrow \omega^{\bullet}_{\pi}$}
{\bf(iii)} $\bullet$ si $\goth{X}$ est la pond\'eration standard,
${\cal C}_{\pi,\goth{X}}$  donne  la classe fondamentale relative
\par de
 [B4],\par
 $\bullet$ si  $\goth{X}$ est la pond\'eration alg\'ebrique (i.e si $\pi$ est plat ou plus g\'en\'eralement de\par {\it Tor-dimension} finie),  ${\cal C}_{\pi,\goth{X}}$ est le morphisme  classe fondamentale relative de [A.E]\par  \'etendu \`a ce cadre.
 \smallskip\noindent
Une cons\'equence importante de ces  r\'esultats donn\'ee par le
th\'eor\`eme {\bf 3} peut s'\'enoncer sous la forme suivante:
\par\noindent {\it{un morphisme  $\pi$ ouvert et \`a fibres de
dimension pure $n$ est analytiquement g\'eom\'etriquement plat si et
seulement si le faisceau $\omega^{n}_{\pi}$ (d\'efini dans le
th\'eor\`eme 1) est stable par changement de base arbitraire entre
espaces complexes r\'eduits et caract\'eris\'e par la propri\'et\'e
de la trace relative}}.
\par\noindent Ainsi, si $X$ est analytiquement g\'eom\'etriquement
plat sur $S$, $\omega^{n}_{\pi}$ est non seulement stable par
changement de base mais incarne parfaitement le faisceau des formes
m\'eromorphes r\'eguli\`eres
$\tilde{\omega}^{n}_{X/S}$.\smallskip\noindent Pour une discussion
d\'etaill\'ee de ce qui va suivre, nous renvoyons  le lecteur au
{\bf\S.3}. Rappelons simplement  que, dans le cadre alg\'ebrique
Kunz et Waldi ([K.W]) ont construit, pour tout morphisme de
sch\'emas noeth\'eriens sans composantes immerg\'ees(par exemple
r\'eduits) $\pi: X\rightarrow S$ de type fini,
$n$-\'equidimensionnel et g\'en\'eriquement lisse et tout entier
$k\leq n$,  des faisceaux ${\cal O}_{X}$-coh\'erents
$\tilde{\omega}^{k}_{X/S}$ enti\`erement caract\'eris\'es par la
{\it propri\'et\'e de la trace relative} et appel\'es faisceaux des
$k$-formes  {\it m\'eromorphes r\'eguli\`eres relatives}. De plus,
pour  $\pi$ plat en caract\'eristique nulle, on dispose d'un
morphisme canonique  ${\cal C}_{\pi}:\Omega^{n}_{X/S}\rightarrow
\tilde{\omega}^{n}_{X/S}$ qui est le morphisme {\it classe
fondamentale }de Elzein ([E]). Cela montre, que, dans cette
situation particuli\`ere, le faisceau $\omega^{n}_{X/S} :={\cal
H}^{-n}(\pi^{!}({\cal O}_{S}))$ est exactement le faisceau de
Kunz-Waldi des $n$-formes m\'eromorphes r\'eguli\`eres. D\`es lors
deux questions naturelles se posent \`a savoir:\par {\bf(a)} {\it
 quelle est la plus large  classe de morphismes satisfaisant les
hypoth\`eses de Kunz-Waldi et dans laquelle le faisceau
$\omega^{n}_{X/S}$ s'identifie canoniquement \`a
$\tilde{\omega}^{n}_{X/S}$?}\par{\bf(b)} {\it dans cette classe
quels sont les morphismes pour lesquels les faisceaux
$\omega^{n}_{X/S}$ sont stables par changement de base arbitraire?}
\smallskip\noindent La classe d\'ecrite par {\bf (a)} n'est pas vide puisqu'elle
 contient, par exemple, les morphismes lisses ou plats, les morphismes sur une base  $S$ qui
 est un sch\'ema r\'egulier (cf [Y]) ou, un peu plus
g\'en\'eralement, un sch\'ema excellent sans composantes immerg\'ees
et v\'erifiant la condition de prolongement ${\bf S}_{2}$ de Serre
(cf [L.S]).\par\noindent En ce qui concerne le changement de base,
on peut dire que si  $X$ est plat sur $S$, $\omega^{n}_{X/S}$ est
stable par tout changement de base. Dans le cas g\'en\'eral, il ne l'est
jamais \`a moins d'avoir un changement de base plat ou d'imposer des
 hypoth\`eses fortes sur $\pi$ ou sur la base $S$ (on peut renvoyer
le lecteur \`a [S] qui traite du cas  d'un morphisme de
Cohen-Macaulay).\smallskip\noindent Sachant que les espaces analytiques complexes sont excellents, notre construction nous
sugg\`ere (toute pr\'ecautions gard\'ees) l'\'enonc\'e alg\'ebrique
suivant:\par\noindent
 {\it Soient $X$ et $S$ des sch\'emas localement noeth\'eriens de dimension de Krull finie
  et sans composantes immerg\'ees
 sur un corps de
caract\'eristique nulle avec $S$ excellent et  $\pi:X\rightarrow S$
un morphisme de type fini, g\'en\'eriquement lisse et
$n$-\'equidimensionnel. Alors $\tilde{\omega}^{n}_{X/S}$ est stable
par changement de base arbitraire si et seulement si $\pi$ d\'efinit
une famille alg\'ebrique de $n$-cycles param\'etr\'ee par
$S$.}\par\noindent On en d\'eduira, \'evidemment, un morphisme
canonique ${\cal C}_{X/S}:\Omega^{n}_{X/S}\rightarrow
\tilde{\omega}^{n}_{X/S}$ coincidant, si $\pi$ est plat,  avec le
morphisme classe fondamentale de Elzein ([E]).\smallskip\noindent
Dans le cadre de la g\'eom\'etrie analytique complexe locale, la
construction des formes m\'eromorphes r\'eguli\`eres relatives, pour
$X$ r\'eduit de dimension pure et plat sur $S$, est implicitement
contenue dans [K.W]. Par des m\'ethodes compl\`etement diff\'erentes
et avec des hypoth\`eses plus g\'en\'erales, Kersken propose, dans
un travail cons\'equent ([Ke], (Ke1]), une construction d'un
complexe de formes r\'eguli\`eres donnant le complexe des formes
m\'eromorphes r\'eguli\`eres relatives de Kunz-Waldi dans le cas
r\'eduit et de dimension pure. Plus pr\'ecisemment, il montre que,
pour tout morphisme plat $\phi:{\cal P}\rightarrow {\cal A}$
d'alg\`ebres analytiques locales (avec ${\cal A}$ non
n\'ecessairement r\'eduite ni de dimension pure!), muni de
l'alg\`ebre diff\'erentielle $(\Omega^{\bullet}_{{\cal A}/{\cal
P}},d)$ des formes diff\'erentielles relatives dot\'ees  de la
diff\'erentielles relative usuelle, il existe un complexe
diff\'erentiel de $(\Omega^{\bullet}_{{\cal A}/{\cal P}},d)$-modules
${\cal D}_{\Omega}({\cal A}/{\cal P})$  s'exprimant au moyen d'un
complexe de cousin relatif et  dont la cohomologie de degr\'e $0$
donn\'ee par ${\rm Ker}: {\cal D}_{\Omega}({\cal A}/{\cal
P})^{0,\star}\rightarrow {\cal D}_{\Omega}({\cal A}/{\cal
P})^{1,\star}$
 est le  $(\Omega_{{\cal A}/{\cal P}}, d)$- module
des formes (g\'en\'eriquement holomorphes) r\'eguli\`eres
$\omega^{\bullet}_{{\cal A}/{\cal P}}$. Signalons, au passage, que
ce type de formes, dont la construction utilisent aussi les
complexes r\'esiduels dans un cadre alg\'ebrique,  se trouvent dans
les travaux de Elzein [E].\par\noindent
  Si ${\cal A}$ est r\'eduite de dimension pure, ces formes r\'eguli\`eres
  sont des formes m\'eromorphes r\'eguli\`eres au sens de Kunz-Waldi. Ainsi se trouvent d\'efinies, pour tout morphisme plat d'espaces
analytiques complexes $\pi:X\rightarrow S$ et  en tout degr\'e $k$,
des faisceaux $\omega^{k}_{X/S}$ de formes dites {\it
r\'eguli\`eres} qui, pour $X$ r\'eduit de dimension pure, coincident
avec les faisceaux de Kunz-Waldi
$\tilde{\omega}^{k}_{X/S}$.\par\noindent Si $\pi$ est un morphisme
propre d'espaces analytiques complexes d\'enombrables \`a l'infini
de dimension localement finie ayant des fibres de dimension pure
$n$, la dualit\'e analytique relative de Ramis-Ruget-Verdier
([R.R.V]) permet de montrer que le faisceau ${\cal
O}_{X}$-coh\'erent ${\cal H}^{-n}(\pi^{!}({\cal O}_{S}))$ (qui est
de profondeur au moins deux fibre par fibre sur $S$) joue le r\^ole
d'un faisceau de Grothendieck relatif. Il est, alors,  facile de
voir que, pour $\pi$ plat,   $\omega^{n}_{X/S}$ s'identifie
canoniquement \`a ${\cal H}^{-n}(\pi^{!}({\cal O}_{S}))$ et, par la
force des choses, aussi \`a  $\tilde{\omega}^{n}_{X/S}$ si $X$ est r\'eduit de dimension pure.\smallskip\noindent Par essence, les  faisceaux  $\tilde{\omega}^{n}_{X/S}$ et ${\omega}^{n}_{X/S}$ sont, en g\'en\'eral,  de nature radicalement diff\'erente. Le premier, de nature locale puisque sa description s'exprime en terme de germes de param\'etrisations locales, est parfaitement appropri\'e pour une description ``analytique'' des multiplicit\'es au moyen des morphismes r\'esidus et le second, de nature globale, doit son existence \`a une int\'egration globale rev\^etant un certain caract\`ere universel et d'aspect plus topologique.\par\noindent
Se d\'ebarasser des  hypoth\`eses  de platitude impos\'ees par la construction de Kersken   et de propret\'e, inh\'erente  \`a la dualit\'e relative analytique de Ramis-Ruget-Verdier n\'ecessite de comprendre en profondeur la probl\'ematique qui nous occupe. On s'apper\c coit que la g\'en\'eralisation repose  fondamentalement sur une dualit\'e relative partielle et une classe de  morphismes dont la factorisation locale produise des morphismes traces.\par\noindent
 L'un des inter\^et du  pr\'esent travail est de mettre clairement en \'evidence ces constats en  proposant deux m\'ethodes reposant, pour l'une, sur une
 approche globale et, pour l'autre, sur une approche locale.
Nous avons volontairement
 commencer par  d\'evelopper la premi\`ere , dans laquelle le
 morphisme d'int\'egration occupe une place centrale, pour assurer une
 certaine continuit\'e avec [KIII]
 qui traite de la notion de {\it paire dualisante}; ce qui fait l'objet du
 {\it th\'eor\`eme 2}. Le {\it th\'eor\`eme 4}, quant \`a lui, utilise l'approche
 locale partant de  la situation plong\'ee et de la classe fondamentale
  relative  pond\'er\'ee donn\'ee par le {\it th\'eor\`eme 0}
  aboutissant au fait que, pour tout morphisme $\pi:X\rightarrow S$ $n$-analytiquement
   g\'eom\'etriquement plat d'espaces, il existe un (unique)
   faisceau ${\cal O}_{X}$-coh\'erent (unique \`a
isomorphisme canonique pr\`es) ${\Lambda}^{n}_{X/S}$ sur $X$ de
profondeur au moins deux fibre par fibre sur $S$, compatible aux
changements de bases et tel que, en chaque point $x$ de $X$ en
lequel $\pi$ est plat, le germe ${\Lambda}^{n}_{X/S,x}$ coincide
avec le faisceau des formes m\'eromorphes r\'eguli\`eres relatives
caract\'eris\'e par la propri\'et\'e de la trace relative au sens de
Kunz-Waldi-Kersken.\par\noindent  De plus, le complexe
$\Lambda^{\bullet}_{X/S}:={\cal H}om(\Omega^{n-\bullet}_{X/S},
\Lambda^{n}_{X/S})$ est munie d'une diff\'erentielle non triviale
${\rm D}$, faisant de $(\Lambda^{\bullet}_{X/S}, {\rm D})$  un
complexe diff\'erentiel de $(\Omega^{\bullet}_{X/S},
d_{X/S})$-modules. \par\noindent Comme on l'a d\`ej\`a dit plus
haut, la classe des morphismes analytiquement g\'eom\'etriquement
plats contient strictement celle des morphismes plats d'espaces
complexes r\'eduits.
\bigskip\par\noindent
Pour mener \`a bien ce programme, on dispose de deux ingr\'edients
essentiels \`a savoir :\smallskip\noindent {\bf(a)} le th\'eor\`eme
de Reiffen ou une de ses variantes assurant l'annulation des
faisceaux de cohomologie ${\rm I}\!{\rm R}^{k}\pi_{!}{\cal F}$, en
tout degr\'e $k>n$, et dont une cons\'equence est le lemme du
d\'ecoupage (cf [B2], [B.V] ou [K1])\par\noindent {\bf(b)} une
connaissance chirurgicale de l'int\'egration sur les cycles ou
op\'eration de Andr\'eotti-Norguet (cf [A.N], [B.V] ou [K2]) et
quelques techniques usuelles.\smallskip\noindent
 La strat\'egie consiste \`a  localiser puis globaliser en recollant convenablement les donn\'ees locales recueillies. Pour cela, on utilise les  techniques de d\'ecoupage des classes de cohomologie et
installation locale du morphisme relativement \`a des \'ecailles (ou
cartes), pour  nous ramener au cas o\`u $\pi$ est factoris\'e en un
plongement $\sigma$ dans $Z$ lisse sur $S$ de dimension relative
$n+p$, suivi d'une projection lisse sur $S$ (ou bien d'un morphisme
fini sur un espace lisse sur $S$ suivi d'une projection sur $S$).
Dans une telle  configuration locale o\`u $X$ est  plong\'e avec la
codimension $p$ fibre par fibre dans $Z$, cette int\'egration
accouche naturellement  d' un morphisme ${\cal
C}^{\sigma}_{\pi}:\sigma_{*}\Omega^{n}_{X/S}\rightarrow {\cal
E}xt^{p}({\cal O}_{X}, \Omega^{n+p}_{Z/S})$, v\'erifiant les
propri\'et\'es de {\bf(2)}. Par ailleurs, une version
``pond\'er\'ee'' de  [B4] donn\'ee par le {\it th\'eor\`eme 0} du
{\bf\S.4}) nous dit que, dans cette situation locale,
   $\pi$ d\'efinit une famille analytique de cycles si et seulement si
 la famille g\'en\'eriquement holomorphe $(\pi^{-1}(s))_{s\in S}$ est repr\'esent\'ee par une
  classe de cohomologie
 dans ${\rm H}^{p}_{|X|}(Z, \Omega^{p}_{Z/S})$, induisant la classe
 fondamentale de $X_{s}$ pour chaque $s$ fix\'e. Mais il est   facile de voir  qu'un tel objet caract\'erise et est caract\'eris\'e par
 la donn\'ee d'un  morphisme  ${\cal C}^{\sigma}_{X/S}: \sigma_{*}\Omega^{n}_{X/S}\rightarrow {\cal H}^{p}_{|X|}(\Omega^{n+p}_{Z/S})$
 compatible
  aux changements de bases entre espaces complexes  r\'eduits  et poss\'edant la
   fameuse propri\'et\'e de la trace. On v\'erifie, sans peine, que
   ${\cal C}^{\sigma}_{\pi}$ est un rel\`evement naturel de  ${\cal C}^{\sigma}_{X/S}$.
   En proc\'edant par recollement, on produit un morphisme global
$\Omega^{n}_{X/S}\rightarrow \omega^{n}_{\pi}$ ayant toutes les
propri\'et\'es voulues et induisant  un morphisme de complexes
diff\'erentiels ${\cal
C}^{\bullet}_{\pi}:\Omega^{\bullet}_{X/S}\rightarrow
\omega^{\bullet}_{\pi}$ (avec diff\'erentielle non triviale sur
$\omega^{\bullet}_{\pi}$).
 \par\noindent Il est facile de voir que, pour un morphisme plat ou plus g\'en\'eralement de {\it Tor-dimension finie}  d'espaces complexes
  avec base r\'eduite, notre construction g\'en\'eralise celle de  [E] et
  [A.E].
\bigskip\noindent
Ce travail est compos\'e de trois  parties. Dans  la premi\`ere, on \'enonce les principaux r\'esultats et  donnons  les notions fondamentales  \'etay\'ees de nombreux exemples et contre-exemples agr\'ement\'ees de quelques petits r\'esultats plus ou moins connus. La seconde est enti\`erement consacr\'ee \`a la preuve des th\'eor\`emes et leurs cons\'equences  annonc\'es pr\'ec\'edemment. Dans la troisi\`eme, on examine, dans ce contexte,  la notion de paire dualisante introduite par Kleinman ([Kl]). \bigskip\noindent
 {\bf 0.1. Enonc\'es des principaux r\'esultats.}\smallskip\noindent
 \Th{1}{} Soit $\pi$ un \'el\'ement de
  ${\cal E}(S,n)$.  Alors, il existe un unique faisceau de ${\cal O}_{X}$-modules
${\omega}^{n}_{\pi}$, v\'erifiant :\par\noindent {\bf (i)} Si
$\xymatrix{X\ar@/_/[rr]_{\pi}\ar[r]^{\sigma}&Z\ar[r]^{q}&S}$ est une
factorisation locale de $\pi$
 dans laquelle  $\sigma$ est un plongement local dans $Z$ lisse sur $S$ et
 de dimension relative $n+p$, alors
$$\omega^{n}_{\pi}\simeq \sigma^{*}{\cal E}xt^{p}(\sigma_{*}{\cal
O}_{X}, \Omega^{n+p}_{Z/S})$$ {\bf (ii)} Il est ${\cal O}_{X}$
coh\'erent, de profondeur au moins deux fibres par fibres sur
$S$\smallskip\noindent {\bf(iii)} $\omega^{n}_{\pi}$ et
$\Omega^{n}_{X/S}$ coincident canoniquement sur  la partie
r\'eguli\`ere du morphisme $\pi$.\smallskip\noindent {\bf (iv) } Sa
construction est compatible aux inclusions ouvertes sur $X$ dans le
sens suivant:\par\noindent
 si   $\pi_{i}:X_{i}\rightarrow S,\,\,i=1,2$ sont deux morphismes analytiques  universellement $n$-\'equidimensionnels  et $U$ un ouvert de $X_{1}$ muni
  de deux inclusions ouvertes $j_{i}:U\rightarrow X_{i}$
  tels  que le diagramme
  $$
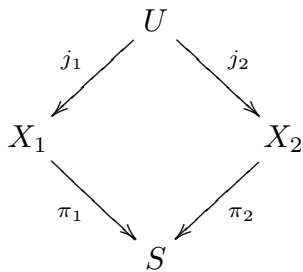
\xymatrix{&U\ar[ld]_{j_{1}}\ar[rd]^{j_{2}}&\\
  X_{1}\ar[rd]_{\pi_{1}}&&X_{2}\ar[ld]^{\pi_{2}}\\
  &S&}$$
  soit commutatif, on a
  $j_{1}^{*}(\omega^{n}_{\pi_{1}}) =
  j_{2}^{*}(\omega^{n}_{\pi_{2}})$
et donc, en particulier,  pour tout ouvert  $U$ de $X$ muni de
l'injection naturelle $j:U\rightarrow X$ et de la restriction de
$\pi$ \`a $U$ que l'on note $\pi_{U}$, on a
$\omega^{n}_{\pi}|_{U}=\omega^{n}_{\pi_{U}}$. De plus, il est stable par tout changement de base plat entre espaces complexes.\smallskip\noindent {\bf
(v)} Il munit $\pi$ d'un morphisme canonique $\int_{\pi}:{\rm I}\!{\rm
R}^{n}\pi_{!}\omega^{n}_{\pi}\rightarrow {\cal O}_{S}$ de formation
compatible aux restrictions ouvertes sur $X$ (resp. $S$) et aux changements de bases plats sur $S$ .\smallskip\noindent
{\bf (vi) } Si $\pi$ est propre, $\omega^{n}_{\pi}= {\cal H}^{-n}(\pi^{!}({\cal O}_{S}))$\smallskip\noindent
{\bf(vii)} Tout  diagramme commutatif d'espaces analytiques complexes
$$\xymatrix{X_{2}\ar[rr]^{\Psi}\ar[rd]_{\pi_{2}}&&X_{1}\ar[ld]^{\pi_{1}}\\
&S&}$$
 avec $\pi_{1}$ (resp. $\pi_{2}$)  universellement \'equidimensionnels et propres  de dimension relative $n_{1}$ (resp. $n_{2}$), $\Psi$ propre de dimension relative born\'ee par l'entier   $d:=n_{2}-n_{1}$,
 donne un morphisme canonique ${\rm I}\!{\rm
 R}^{d}{\Psi}_{*}\omega^{n_{2}}_{\pi_{2}}\rightarrow \omega^{n_{1}}_{\pi_{1}}$ induisant le diagramme commutatif de faisceaux analytiques
$$\xymatrix{{\rm I}\!{\rm
 R}^{n_{2}}{\pi_{2}}_{*}\omega^{n_{2}}_{\pi_{2}}\ar[rr]\ar[rd]_{\int_{\pi_{2}}}&&
 {\rm I}\!{\rm
 R}^{n_{1}}{\pi_{1}}_{*}\omega^{n_{1}}_{\pi_{1}}\ar[ld]^{\int_{\pi_{1}}}\\
&{\cal O}_{S}&}$$
{\bf(viii)} Si $S$ est un point, il
coincide avec  le faisceau dualisant de Grothendieck, Andr\'eotti-
Kas ou Golovin :
$$\omega^{n}_{\pi}\simeq {\cal H}_{n}({\cal O}_{X})\simeq {\cal D}^{n}({\cal O}_{X})=
\omega^{n}_{X}$$\rm\smallskip\noindent
 \Cor{1.1}{}Soient $\pi\in{\cal E}(S,n)$ et $k$ un entier naturel. Alors, il existe un unique faisceau ${\cal O}_{X}$-coh\'erent  $\omega^{k}_{\pi}$  v\'erifiant les propri\'et\'es suivantes :\smallskip\noindent
 {\bf(i)} il est de profondeur au moins $2$ fibre
par fibre sur $S$, coincide avec le faisceau des $k$-formes
holomorphes relatives sur la partie r\'eguli\`ere de $\pi$ et, pour
toute installation locale de $\pi$ (cf Thm 1), on a
$$\omega^{k}_{\pi}:=\sigma^{*}({\cal E}xt^{p}(\sigma_{*}\Omega^{n-k}_{X/S},
  \Omega^{n+p}_{Z/S}))$$
\par\noindent {\bf(ii)} on a des isomorphismes canoniques
  $\omega^{k}_{\pi}\simeq {\cal
H}om_{{\cal O}_{X}}(\Omega^{n-k}_{X/S},
\omega^{n}_{\pi})$\par\noindent
\par\noindent
{\bf(iii)} Si $S$ est un point, ces faisceaux coincident avec les
dualis\'es de Andr\'eotti- Kas- Golovin  des faisceaux
$\Omega^{n-k}_{X}$. \rm\smallskip\smallskip\noindent
\smallskip\noindent \Th{2}{} Soit  $\pi$ un \'el\'ement de ${\cal
E}(S,n)$ muni d'une pond\'eration ${\goth X}$.  Alors, on a les
\'equivalences :\smallskip\noindent {\bf(1)} $\pi$ est
analytiquement g\'eom\'etriquement plat\smallskip\noindent {\bf(2)}
il existe un unique morphisme ${\cal O}_{S}$-lin\'eaire continu
\footnote{$^{(1)}$}{On entend par l\`a, une continuit\'e topologique
au niveau des  groupes  des  sections globales munis de leurs
structures d'espaces vectoriels topologiques {\bf{F.S}} ou
{\bf{Q.F.S}}} $\ds{{\int}_{\pi,{\goth X}}:{\rm I}\!{\rm
R}^{n}\pi_{!}\Omega^{n}_{X/S}\rightarrow {\cal O}_{S}}$  de
formation compatible \`a l'additivit\'e des pond\'erations sur
$\pi$, stable par changement de base entre espace complexes
r\'eduits, donnant, en particulier,  l'int\'egration usuelle dans le
cas o\`u $S$ est un point.\par\noindent {\bf(3)}  il existe  un
unique morphisme canonique  de faisceaux coh\'erents
 ${\cal C}_{\pi,\goth{X}}:\Omega^{n}_{X/S}\rightarrow \omega^{n}_{\pi}$,
prolongeant  le morphisme
naturel $j_{*}j^{*}\Omega^{n}_{X/S}\rightarrow
j_{*}j^{*}\omega^{n}_{\pi}$ et
 v\'erifiant les propri\'et\'es suivantes:\par
{\bf(i)} il est de formation compatible \`a l'additivit\'e des
pond\'erations,   aux
 changement de base r\'eduit et satisfait  {\it{la propri\'et\'e de la
 trace relative}},\par
{\bf(ii)} il induit  un morphisme de complexes  diff\'erentiels
gradu\'es \footnote{$^{(2)}$}{Ce n'est pas un morphisme d'alg\`ebres
car  $\omega^{\bullet}_{\pi}$ n'en est pas une puisque l'on a pas de
produit interne. Cette lacune est d\'ej\`a pr\'esente dans le cas
absolu comme on peut le voir sur  $X=\{(x,y)\in {\Bbb C}^{2}:
x^{2}-y^{3}=0\}$. En effet, la fonction m\'eromorphe ${y\over{x}}$
(resp. la forme ${dx\over{y}}$) d\'efinit une section de
$\omega^{0}_{X}$ (resp. $\omega^{1}_{X}$) mais leur produit
${dx\over{x}}$ ne d\'efinit pas une section de  $\omega^{1}_{X}$ !
Cependant, ce produit est interne si $X$ est normal.}\par
 \centerline{${\cal C}^{\bullet}_{\pi}:\Omega^{\bullet}_{X/S}\rightarrow \omega^{\bullet}_{\pi}$}
 {\bf(iii)} $\bullet$ si $\goth{X}$ est la pond\'eration standard,  ${\cal C}_{\pi,\goth{X}}$  donne  la classe fondamentale relative \par de
 [B4],\par
 $\bullet$ si  $\goth{X}$ est la pond\'eration alg\'ebrique (i.e si $\pi$ est plat ou plus g\'en\'eralement de\par {\it Tor-dimension} finie),   ${\cal C}_{\pi,\goth{X}}$ est le morphisme  classe fondamentale relative de [A.E]\par  \'etendu \`a ce cadre.
\rm\smallskip\medskip\par\noindent La fl\`eche ${\cal C}_{\pi,\goth{X}}$  sera
appel\'ee  {\it{morphisme classe fondamentale relative pond\'er\'e }} de $\pi$.
\smallskip\medskip\noindent
\Cor{2.1}{}Soit $\pi_{2}:X_{2}\rightarrow X_{1}$ (resp.
$\pi_{1}:X_{1}\rightarrow S$) un morphisme analytiquement
g\'eom\'etriquement plat de dimension relative $n_{2}$ (resp.
$n_{1}$). Alors, le morphisme compos\'e $\pi:X_{2}\rightarrow S$ est
analytiquement g\'eom\'etriquement plat si et seulement si le
morphisme naturel ${\rm I}\!{\rm R}^{n_{1}+n_{2}}{\pi}_{!}(
\Omega^{n_{2}}_{X_{2}/X_{1}} \otimes
\pi_{2}^{*}\Omega^{n_{1}}_{X_{1}/S})\rightarrow{\cal O}_{S}$ se
prolonge en unique morphisme ${\rm I}\!{\rm
R}^{n_{1}+n_{2}}{\pi}_{!}\Omega^{n_{1}+n_{2}}_{X_{2}/S}\rightarrow{\cal
O}_{S}$ ayant les propri\'et\'es du th\'eoreme 2 et rendant
commutatif le diagramme
$$\xymatrix{ {\rm I}\!{\rm R}^{n_{1}+n_{2}}{\pi}_{!}(\Omega^{n_{2}}_{X_{2}/X_{1}} \otimes
\pi_{2}^{*}\Omega^{n_{1}}_{X_{1}/S} )\ar[rd]\ar[rr]&&
 {\rm I}\!{\rm R}^{n_{1}+n_{2}}{\pi}_{!}\Omega^{n_{1}+n_{2}}_{X_{2}/S}\ar[ld]^{\Phi}\\
&{\cal O}_{S}&}$$
\Th{3}{} Soit  $\pi$ un \'el\'ement de ${\cal
E}(S,n)$ muni d'une pond\'eration ${\goth X}$. Alors, les assertions
suivantes sont \'equivalentes
\smallskip\noindent
 {\bf(i)} $\pi$ est analytiquement g\'eom\'etriquement plat\smallskip\noindent
{\bf (ii)} le faisceau  $\omega^{n}_{\pi}$ est caract\'eris\'e par
la propri\'et\'e de la trace (relative) et stable par changement de
base entre espaces complexes r\'eduits.\rm\smallskip\par\noindent
Comme nous l'avions mentionn\'e dans l'introduction,  la
construction, de nature purement alg\'ebrique,  du faisceau des
formes  m\'eromorphes r\'eguli\`eres relatives, entreprise par  Kunz
et Waldi([K.W]), se transporte \`a la cat\'egorie des germes
d'espaces  analytique complexes  pour un morphisme plat gr\^ace \`a
l'important et cons\'equent travail de Kersken ([Ke]). En dehors du
cas propre, pour lequel on utilise la dualit\'e analytique relative
de [R.R.V],
 la globalisation d'un tel travail ne semble pas du tout \'evidente.
 Il s'av\`ere que, sans emprunter les m\'eandres d'une dualit\'e
 analytique relative g\'en\'erale  inachev\'ee, l'existence de l'analogue analytique relatif
  du faisceau de Kunz-Waldi existe non seulement dans le cas plat mais plus
  g\'en\'eralement dans le cas analytiquement g\'eom\'etriquement plat comme
  le montre le
\Cor{3.1}{} Soient\par\noindent $\bullet$  $\pi:X\rightarrow S$ un
morphisme $n$-analytiquement g\'eom\'etriquement plat  d'espaces
analytiques complexes r\'eduits de dimension localement
fini,\par\noindent $\bullet$ $X_{0}$ (resp. $S_{0}$) l'ouvert dense
de $X$ (resp. $S$) sur lequel $\pi$ est plat et
 $\tilde{\omega}^{n}_{X_{0}/S_{0}}$ le faisceau des
 formes m\'eromorphes r\'eguli\`eres de Kunz-Waldi-Kersken.\smallskip\noindent
 Alors, $\omega^{n}_{\pi}$ est l'unique prolongement coh\'erent du
 faisceau  $\tilde{\omega}^{n}_{X_{0}/S_{0}}$. De plus, la famille de
faisceaux ${\cal O}_{X}$-coh\'erents $\omega^{\bullet}_{\pi}:={\cal
H}om(\Omega^{n-\bullet}_{X/S}, \omega^{n}_{\pi})$ est munie d'une
diff\'erentielle non triviale ${\rm D}$, faisant de
$(\omega^{\bullet}_{X/S}, {\rm D})$  un complexe diff\'erentiel de
$(\Omega^{\bullet}_{X/S}, d_{X/S})$-modules.
\rm\smallskip\par\noindent
Le r\'esultat suivant, qui  peut-\^etre \'etabli de fa\c con presqu'ind\'ependante de ce qui pr\'ec\`ede, montre que le caract\`ere ``analytiquement g\'eom\'etriquement plat'' d'un morphisme contient suffisamment d'informations pour pouvoir accoucher d'une structure dualisante sans faire explicitement appel \`a une quelconque th\'eorie de la dualit\'e.
\Th{4}{} Soit $\pi:X\rightarrow S$ un morphisme analytiquement
g\'eom\'etriquement plat. Alors, il existe un unique faisceau ${\cal
O}_{X}$ coh\'erent, $\Lambda^{n}_{X/S}$ v\'erifiant:\par\noindent
{\bf(i)} il est  de profondeur au moins deux fibre par fibre sur
$S$, caract\'eris\'e par la propri\'et\'e de la trace relative et
stable par changement de base entre espaces complexes r\'eduits de
dimension finie,\par\noindent {\bf(ii)} il munit $\pi$ d'un
morphisme d'int\'egration $\int_{\pi}:{\rm I}\!{\rm
R}\pi_{!}\Lambda^{n}_{X/S}\rightarrow {\cal O}_{S}$ compatible avec
l'additivit\'e des cycles et stable par changement de
base entre espaces complexes r\'eduits. De plus, tout  diagramme commutatif d'espaces analytiques complexes De plus, tout  diagramme commutatif d'espaces analytiques complexes
$$\xymatrix{X_{2}\ar[rr]^{\Psi}\ar[rd]_{\pi_{2}}&&X_{1}\ar[ld]^{\pi_{1}}\\
&S&}$$
 avec $\pi_{1}$ (resp. $\pi_{2}$) analytiquement g\'eom\'etriquement plats  de dimension relative $n_{1}$ (resp. $n_{2}$)  et $\Psi$ universellement \'equidimensionnel et propre de dimension relative born\'ee par l'entier   $d:=n_{2}-n_{1}$,  donne un morphisme canonique ${\rm I}\!{\rm
 R}^{d}{\Psi}_{*}\Lambda^{d+n_{1}}_{X_{2}/S}\rightarrow \Lambda^{n_{1}}_{X_{1}/S}$ induisant le diagramme commutatif de faisceaux analytiques
$$\xymatrix{{\rm I}\!{\rm
 R}^{n_{2}}{{\pi_{2}}_{!}}\Lambda^{n_{2}}_{X_{2}/S}\ar[rr]\ar[rd]_{\int_{\pi_{2}}}&&
 {\rm I}\!{\rm
 R}^{n_{1}}{{\pi_{1}}_{!}}\Lambda^{n_{1}}_{X_{1}/S}\ar[ld]^{\int_{\pi_{1}}}\\
&{\cal O}_{S}&}$$
\rm
\smallskip\smallskip\bigskip\noindent\eject\vfill\noindent {\tite I. Notion
de  morphisme contin\^ument ou analytiquement g\'eom\'etriquement
plat.}\smallskip\smallskip\smallskip\noindent {\bf 1.0. Terminologie
usuelle relative \`a l'espace des cycles  selon
[B1].}\smallskip\noindent Dans la construction d'une structure
analytique
 sur l'ensemble des cycles compacts effectifs de dimension pure d'un espace complexe
 donn\'e, l'espace quotient normal
${\rm Sym}^{k}({\Bbb C}^{p}):= {({\Bbb C}^{p})^{k}/\sigma_{k}}$
($\sigma_{k}$ \'etant le groupe sym\'etrtique d'ordre $k$) joue un
r\^ole capital et fondamentale. En effet, l'id\'ee directrice de
Barlet dans [B1] \'etait de dire qu'une famille de $n$-cycles
effectifs d'un espace complexe $X$  est analytique si et seulement
si pour tout plongement local "adapt\'e", l'intersection locale de
chaque membre de cette famille avec des plans transverses (et de
codimension compl\'ementaire) donne une famille analytique de
$0$-cycles. Or, l'espace des $0$-cycles d'un espace complexe
quelconque $X$, s'identifie naturellement \`a la somme disjointe
$\ds{\coprod_{r\geq 0}{\rm Sym}^{r}(X)}$. On peut remarquer que
cette construction  est propre \`a la g\'eom\'etrie analytique
complexe (et en grande partie \`a la caract\'eristique nulle du
corps de base !) puisqu'elle utilise des outils de l'analyse
complexe tels que les formules int\'egrales de Cauchy, les relations
non triviales entre {\it Trace} et {\it D\'eterminant} d'un
endomorphisme et le fait que l'on puisse exprimer en termes de
fonctions de Newton les fonctions sym\'etriques des racines d'un
polyn\^ome donn\'e. Ce dernier point extr\^ement important permet,
une fois rep\'er\'es  les points de  ${\rm Sym}^{k}({\Bbb C}^{p})$
par leurs fonctions sym\'etriques,de  produire un plongement naturel
de ${\rm Sym}^{k}({\Bbb C}^{p})$ dans un espace num\'erique
convenable. De ce fait, la structure analytique globale recherch\'ee
est produite par recollement de structures analytiques locales
parfaitement et pr\'ecisement d\'ecrites.
\smallskip\par\noindent
{\bf 1.0.1. Quelques d\'efinitions.}\smallskip\noindent {\bf 1.0.1.1.}   Si $Z$ est
un espace analytique complexe r\'eduit et $n$ un entier naturel, on
appelle  $n$- {\it{cycle analytique effectif}} de $Z$  la donn\'ee
d'une combinaison lin\'eaire formelle localement finie du type
$\ds{W:= \sum_{i\in I} n_{i}W_{i}}$\footnote{$^{(3)}$}{Pour un sous
ensemble analytique $Y$ de composantes irr\'eductibles
$(Y_{i})_{i\in I}$, le cycle associ\'e \`a $Y$ est not\'e
$\ds{[Y]:=\sum_{i\in I}[Y_{i}]}$ o\`u $[Y_{i}]:= 1.Y_{i}$ }  o\`u
les $W_{i}$ sont des sous ensembles analytiques complexes
irr\'eductibles de dimension pure $n$ dans  $Z$, et $n_{i}$ des
entiers strictement positifs appel\'es multiplicit\'es; le sous
ensemble analytique  $|W|:=\ds\bigcup_{i\in I}W_{i}$ est appel\'e
{\it  support}  de $W$.\par\noindent Soit $V$ un ouvert de $Z$ muni
d'un plongement $\sigma:V\rightarrow U\times B\subset {\Bbb C}^{N}$
avec $U$ (resp. $B$) un polydisque ouvert de  ${\Bbb C}^{n}$ (resp.
${\Bbb C}^{N-n}$). La carte   $E:=(\sigma,V,U\times B)$ est dite
{\it{adapt\'ee}} au cycle $W$ si $\sigma$ se prolonge en un
plongement $\sigma':V'\rightarrow U'\times B'$ avec  $V\Subset
V'\subset Z$,  $U\Subset U'\subset {\Bbb C}^{n}$, $B\Subset
B'\subset {\Bbb C}^{N-n}$ et $\sigma(|W|\cap
V)\cap(\overline{U}\times {\partial
B})=\emptyset$\footnote{$^{(4)}$}{Dans cette situation, pour tout
$t$ dans $U$,  le cycle intersection $W\bullet(\{t\}\times B)$ est
bien d\'efini dans $U\times B$ et se projette en un $0$-cycle de
$B$.}\par\noindent  Dans ce cas, pour tout $i$ dans  $I$, la restriction de la
projection
 canonique $p_{1}:U\times B\rightarrow U$ \`a  $W_{i}\cap V$ est un  rev\^etement ramifi\'e d'un certain degr\'e $k_{i}$ sur $U$ et d\'efinit donc une application analytique $F_{E,i}:U\rightarrow {\rm Sym}^{k_{i}}({\Bbb
C}^{p})$, associant  \`a $t$ fix\'e dans $U$  les fonctions
sym\'etriques des branches locales de ce  rev\^etement ramifi\'e.
L'application naturelle $\ds{({\Bbb C}^{p})^{r}\times ({\Bbb
C}^{p})^{l}\rightarrow ({\Bbb C}^{p})^{r+l}}$ qui induit, par passage au
quotient, une application canonique  ${\rm Sym}^{r}({\Bbb
C}^{p})\times {\rm Sym}^{l}({\Bbb C}^{p})\rightarrow {\rm
Sym}^{r+l}({\Bbb C}^{p})$  (qui est un hom\'eomorphisme et m\^eme un
isomorphisme local) permet  de d\'efinir le {\it degr\'e} du
cycle $W$ dans cette carte $E$  comme \'etant l'entier $\ds{{\rm
deg}_{E}(W)=k:=\sum_{i}k_{i}n_{i}}$  qui est le degr\'e total du
rev\^etement ramifi\'e donn\'ee par l'application analytique
$\ds{F_{E}:U\rightarrow {\rm Sym}^{k}({\Bbb C}^{p})}$ d\'efinie par
$\ds{F_{E}(t):=\prod_{i}(F_{E,i}(t))^{n_{i}}}$. En d\'esignant par
 $S_{h}({\Bbb C}^{p})$ la composante homog\`ene de degr\'e $h$ de l'alg\`ebre sym\'etrique de ${\Bbb C}^{p}$ qui s'identifie naturellement \`a l'espace vectoriel des polyn\^omes homog\`enes de degr\'e $h$ sur $({\Bbb C}^{p})^{*}$, il est facile de voir que la somme directe des applications $\ds{s_{h}(X_{1},\cdots, X_{k}):=\sum_{1\leq
i_{1}<\cdots<i_{h}\leq h}X_{i_{1}}\cdots X_{i_{h}}}$, d\'efinies de  $({\Bbb C}^{p})^{k}$ \`a valeurs dans $S_{h}({\Bbb C}^{p})$, induit, par
passage au quotient,  un  plongement naturel $\ds{\sigma:{\rm
Sym}^{k}({\Bbb C}^{p})\rightarrow V:=\bigoplus^{k}_{h=1} S_{h}({\Bbb
C}^{p})}$. Remarquons, au passage,  que le polyn\^ome unitaire
\footnote{$^{(5)}$}{Ce n'est rien d'autre que la version
multi-lin\'eaire du polyn\^ome {\bf{discriminant}} donn\'e par
$\ds{\rm D}(z)=\prod_{i<j}(w_{i}(z)-w_{j}(z))^{2}
=(-1)^{{n(n-1)\over{2}} }\prod^{j=n}_{j=1}{{\partial
{P}\over{\partial w}}(z,w_{j})}$ pour $\ds{P(z,w):=
\sum^{j=n}_{j=0}a_{j}(z)w^{j}}$}
$${\rm D}(X_{1},\cdots,X_{k})[T]:=\prod_{i<j}(T^{2} - (X_{i} - X_{j})^{2})= \sum^{{k(k-1)\over{2}}}_{h=0}(-1)^{h}D_{h}(x)T^{2h}$$
\`a coefficients dans l'alg\`ebre sym\'etrique de ${\Bbb C}^{p}$,
d\'efinit, pour chaque entier $h$, une application $\sigma_{k}$-
invariante de $({\Bbb C}^{p})^{k}$ dans $S_{k(k-1) - 2h}({\Bbb
C}^{p})$ et donc des applications
$${\rm D}_{h}: {\rm Sym}^{k}({\Bbb C}^{p})\rightarrow S_{k(k-1) - 2h}({\Bbb C}^{p})$$
jouant un r\^ole notoire dans la  stratification  alg\'ebrique de
${\rm Sym}^{k}({\Bbb C}^{p})$.\smallskip\par\noindent {\bf 1.0.1.2.
Multiplicit\'e d'un point dans un cycle.}\smallskip\noindent
  Selon [B.M],  la {{\it multiplicit\'e}\rm}
de $x\in {\Bbb C}^{p}$ dans $(x_{1},\cdots,x_{k})\in {\rm
Sym}^{k}({\Bbb C}^{p})$ est le nombre de fois qu'est r\'ep\'et\'e
$x$ dans ce $k$-uplet. Si $U$ est un ouvert de ${\Bbb C}^{n}$,
$F:U\rightarrow {\rm Sym}^{k}({\Bbb C}^{p})$ une application
holomorphe associ\'ee \`a un certain rev\^etement ramifi\'e
$f:X\rightarrow U$ et $z=(t,x)\in U\times {\Bbb C}^{p}$, la {{\it
multiplicit\'e}\rm}  de $z$ dans $X$ est l'entier not\'e
$\mu_{z}(X)$ donnant la multiplicit\'e de $x$ dans le $k$-uplet
$F(t)$.\par
 Si $X$ est un
cycle d'un espace complexe $Z$. Pour toute \'ecaille adapt\'ee $E$,
on note $X_{E}$ le rev\^etement ramifi\'e associ\'e \`a $X$ dans
$E$. Soit $z\in Z$. On appelle {\it{multiplicit\'e}} de $z$ dans le
cycle $X$, l'entier $\mu_{z}(X)$ donn\'ee par\par
\centerline{$\ds{\mu_{z}(X):=min_{E}(\mu_{z}(X_{E})}$} \noindent
quand $E$ d\'ecrit l'ensemble des \'ecailles adapt\'ees \`a
$X$.\par\noindent On peut remarquer que pour tout $z$ dans $Z$, il
existe toujours une \'ecaille adapt\'ee \`a $X$ telle que
$\mu_{z}(X)= deg_{E}(X)$.\par\noindent Soit $Y=\ds{\sum_{i\in
I}n_{i}Y_{i}}$ un cycle et $z$ un point d'un espace complexe $Z$, la
{\it{ multiplicit\'e de $z$ dans le cycle $Y$}} est donn\'ee par\par
\centerline{$\mu_{z}(Y):=\sum_{i\in I}n_{i}\mu_{z}(Y_{i})$}
\smallskip\noindent
On peut signaler au lecteur que la d\'efinition ``alg\'ebrique'' donn\'ee dans \$3.1 de [Fu] coincide avec la d\'efinition pr\'ec\'edente. Par ailleurs, il trouvera expos\'e dans [Si] une d\'efinition d\^ue \`a [Tu] adapt\'ee au cas d'un morphisme ouvert $\pi:X\rightarrow S$ d'espaces complexes
avec $X$ localement de dimension pure  et $S$ localement
irr\'eductible.
\smallskip\par\noindent {\bf 1.0.1.3.
Fonctions et formes de Newton.}\par\noindent {\bf(a)}  Soit $F$ un
${\Bbb C}$- espace vectoriel de dimension finie. Soient $x$ (resp.
$y$) un $k$- ulpet de $({\Bbb C}^{p})^{k}$ (resp. de $F^{k}$). la
$l$-{\it{\`eme fonction de Newton}} de $x$ {\it{pond\'er\'ee par}}
$y$ est l'\'el\'ement de  $S_{l}({\Bbb C}^{p})\otimes F$ d\'efini
par l'expression $\ds{N_{l}(x,y):= \sum^{k}_{j=1}y_{j}\otimes
x^{l}_{j}}$. Alors, pour tout entier $l\leq k$, on a les relations
de Newton vectorielles
$$\sum^{k}_{h=0}(-1)^{h}N_{l-h}(x,y).S_{h}(x) = 0,\,\,\,{S}_{0}:=1$$\par\noindent Le lecteur  trouvera dans le chapitre 0 de [B1],
l'adaptation vectorielle des r\'esultats classiques tournant autour
de cette notion et pourra consulter [Ang] sur les  polyn\^omes et
formes multilin\'eaires de  W\"aring permettant d'exprimer les
fonctions sym\'etriques en termes de fonctions de Newton dans ce
cadre vectoriel.\smallskip\noindent {\bf(b)} Comme nous le verrons,
la notion de {\bf{Formes de Newton}} (cf [B4])  jouent un r\^ole
cruciale dans le cas des familles de  cycles de dimension
strictement positive.   Soit
$$\xymatrix{({\Bbb C}^{p})^{k}\ar[rr]^{q}\ar[rd]_{p_{j}}&&{\rm Sym}^{k}({\Bbb C}^{p})\ar[ld]^{p}\\
&{\Bbb C}^{p}&}$$
le diagramme commutaif naturel dans lequel $q$ est l'application quotient, $p$ la projection
naturelle et  $p_{j}$ la projection canonique de $({\Bbb
C}^{p})^{k}$ sur son $j$-\`eme facteur.\par\noindent Soit $U$ un
ouvert de ${\rm Sym}^{k}({\Bbb C}^{p})$ et  $\ds{p(U)=\lbrace{x\in
{\Bbb C}^{p}: \exists x_{2},\cdots,x_{k}\in {\Bbb C}^{p}\,{\rm
tel}\,{\rm que}\,\,q(x,x_{2},\cdots,x_{k})\in U}\rbrace}$ son image
dans ${\Bbb C}^{p}$.\par\noindent
 Soit
$\omega^{r}:=(q_{*}\Omega^{r})^{\sigma_{k}}$, la partie
$\sigma_{k}$-invariante de l'image directe par $q$ des formes
holomorphes sur $({\Bbb C}^{p})^{k}$. Alors, une forme   $\xi\in
\Gamma(U, \omega^{r})$ est dite {\it{de Newton}} s'il existe une
forme $\eta$ d\'efinie sur l'ouvert $\ds{p(U)}$ telle que $\ds{\xi =
\sum^{k}_{1}p_{j}^{*}\eta}$.\smallskip\noindent $\bullet$ Le
faisceau des $m$- formes de Newton, not\'e  ${\cal N}_{m}$, est
coh\'erent et  engendr\'e par les formes de Newton globales
$\ds{w_{I,J}= \sum^{j=k}_{j=1}p_{j}^{*}(x^{I}dx^{J})}$ o\`u $I$ et $J$ sont des
  ensembles d'indices de longueures  respectives  $a$ et $b$. De plus, pour tout entier
$m\geq 0$, on a  les relations
$$\sum^{h=k}_{h=0}(-1)^{h}S_{h}(x_{1},\cdots,x_{k})w_{m+h,I}
(x_{1},\cdots,x_{k})=0$$
 Desquelles r\'esultent
l'identification $\ds{\omega^{r}=\sum^{r}_{i=1}\omega^{r-i}\wedge
{\cal N}_{i}}$, pour tout entier $r\geq 1$.
\smallskip\par\noindent
{\bf 1.0.2.  Famille analytique de cycles.}\smallskip\noindent {\bf
1.0.2.1. Ensemble des cycles analytiques effectifs d'un espace
complexe.}\par\noindent Si $X$ est un espace analytique complexe
donn\'e, on note ${\cal C}_{*}(X)$(ou ${\cal
C}^{loc}_{*}(X)$)l'ensemble des cycles analytiques (positifs) de
dimension donn\'ee s'\'ecrivant localement comme combinaison
lin\'eaire localement finie \`a coefficients entiers de sous
ensembles analytiques  irr\'eductibles de $X$. Utilisant les
r\'esultats de Federer-Fleming-King, A. Cassa a montr\'e dans [C]
 que l'ensemble des  cycles positifs de dimension pure $n$ d'une
vari\'et\'e analytique complexe est un espace topologique m\'etrique
complet dont la topologie est caract\'eris\'ee par celle de la
convergence en "masse" d\'efinie sur les courants. Il montre, en
particulier, l'existence d'une unique application injective
$\ds{{\cal C}_{n}(X)\rightarrow {\cal D}^{'}_{2n}(X)}$ faisant
correspondre \`a tout cycle de $X$ le courant d'int\'egration qui
lui est naturellement associ\'e gr\^ace au r\'esultat fondamental de
Lelong [L]. Une proc\'edure standard de localisation (par plongement local) et recollement permet d'\'etendre [C] au cas  o\`u $X$ est singulier et d'en d\'eduire  que  ${\cal C}_{*}(X)$ est
un espace topologique de Hausdorff, \`a bases d\'enombrables de
voisinages (i.e "first countable"); ce qui permet de tester la compacit\'e ou fermeture de sous espaces de ${\cal C}_{n}(X)$ en utilisant des suites de cycles.   \par\noindent   Par ailleurs la
construction de Barlet ([B1])  montre que les ouverts
$V_{k,E}:=\{Y\in {\cal C}_{n}(X)/ {\rm E}\,{\rm
adapt\acute{e}}\,{\rm pour }\,Y\,{\rm et}\,\,{\rm Deg}_{E}(Y)=k\}$
forment, \`a $k$ et $E$ variables,
 une base d'ouverts naturels pour ${\cal C}_{*}(X)$. Il se trouve (cf \$2.4, [Fu]) que la topologie de Hausdorff induite par la topologie faible de ${\cal D}^{'}_{2n}(X)$  coincide avec la topologie de Barlet. Pour plus de d\'etails, nous renvoyons le lecteur \`a [Fu],[Si],[Ma] en gardant \`a l'esprit que cette question est d\'
ej\`a trait\'ee presque dans sa totalit\'e  dans [C]. Il va sans dire que la topologie de ${\cal C}_{*}(X)$ est naturellement induite par les projections
${\cal C}_{*}(X)\rightarrow {\cal C}_{n}(X)$ associant \`a chaque cycle de dimension mixte sa partie de dimension pure $n$. \par\noindent
En g\'en\'eral,  l'espace toplogique  ${\cal C}_{*}(X)$ n'est jamais
 localement compact et ne  peut, donc, \^etre dot\'e  d'une
  structure analytique complexe. Cependant, dans le cas des
   cycles compacts, c'est un espaces analytique complexe
   comme l'a montr\'e Barlet dans son cons\'equent travail consign\'e
   dans [B1] et dont nous allons rappeler le th\'eor\`eme fondamental
    de repr\'esentabilit\'e.\smallskip\par\noindent
 {\bf 1.0.2.2. Le th\'eor\`eme de Barlet.}
\par\noindent Dans toute la suite $S$ d\'esignera un espace complexe r\'eduit, localement de dimension finie.\par\noindent
Si $S'$ est une partie de $S$, {\it{une  \'ecaille ou carte
$S'$-adapt\'ee}} \`a une famille $(X_{s})_{s\in S}$ d'un espace
complexe $Z$  correspond \`a la donn\'ee d'une carte $E:=(\sigma,
V,U\times B)$ v\'erifiant la condition
$$\overline{\bigcup^{}_{s\in S'}(\sigma(|X_{s}|\cap \bar{V})}\cap(\bar{U}\times\partial{B})= \emptyset$$\smallskip\noindent
Cela signifie exactement qu'il existe un compact $K$ de $B$ tel que
$$\sigma(|X_{s}|\cap \bar{V})\subset \bar{U}\times K\,\,\,\,\,\,\forall\,s\in S'$$
Sachant que l'espace des $0$-cycles d'un espace complexe $Z$ donn\'e
s'identifie canoniquement \`a la somme disjointe $\ds{\coprod_{r\geq
0}{\rm Sym}^{r}(Z)}$, on aimerait dire
 naturellement qu'une famille continue de
cycles $(X_{s})_{{s}\in S}$ est analytique en un point $s_{0}$ de
$S$ si il existe un certain voisinage ouvert $S_{0}$ de $s_{0}$ tel
que pour toute \'ecaille, $E=(V,U,B,\sigma)$, $S_{0}$-adapt\'ee,
l'application qui \`a $\{s\}$ associe ${p_{1}}_{*}([\{t\}\times
B]\bullet (X_{s}\cap V))$ est holomorphe ($p_{1}$ \'etant la
projection canonique sur le premier facteur de $U\times B$). Cela
traduit le fait que l'intersection de chaque $X_{s}$ par des
$p$-plans transverses donne une famille analytique de
points.\smallskip\noindent
Soit $S$ un espace analytique complexe
r\'eduit et $(X_{s})_{{s}\in S}$  une famille de $n$- cycles
 de $Z$ param\'etr\'ee par $S$. On dit  que cette famille est {\it{analytique}} si pour tout  $s_{0}\in S$ et  {{\bf toute}\rm} \'ecaille $E$ adapt\'ee \`a
$X_{s_{0}}$,
 il existe  un voisinage ouvert $S_{0}$ de $s_{0}$ dans  $S$ tel que \par i)  $E$ est adapt\'ee \`a $X_{s}$ pour tout $s\in S_{0}$ (i.e  $E$ est  $S_{0}$-adapt\'ee).\par
ii) $deg_{E}(X_{s}) = deg_{E}(X_{s_{0}})$ ( en d\'esignant, par
  $deg_{E}(X_{s})$, le degr\'e du rev\^etement ramifi\'e sous
jacent au cycle $X_{s}$ dans l'\'ecaille $E$).\par
 ii)  L'application induite
$F_{E}: S_{0}\times U\rightarrow {\rm Sym}^{k}(B)$ est analytique.
\smallskip\noindent
Dans [B1] chapitre 5, D. Barlet  montre le  r\'esultat suivant:
  \Th{}
Soient $Z$ un  espace  analytique  complexe,  ${\cal C}$  la
cat\'egorie  des  espaces complexes  r\'eduits de dimension finie et
morphismes  analytiques, ${\cal F}_{n}$  le  foncteur contravariant
de ${\cal C}$ dans ${\cal E}$  qui \`a $S$ associe ${{\cal
F}_{n}}(S)$ l'ensemble des familles analytiques de $n$-cycles
compacts de $Z$  param\'etr\'ee par $S$. Alors, pour tout entier
naturel $n$, ${\cal F}_{n}$  est un foncteur
repr\'esentable.\smallskip\noindent \rm En d'autres termes, il
existe un espace analytique complexe r\'eduit, not\'e ${\cal
B}_{n}(Z)$, et une famille analytique de $n$-cycles compacts
$\{X_{Y} / Y\in B_{n}(Z)\}$ de $Z$ param\'etr\'ee par ${\cal
B}_{n}(Z)$ v\'erifiant la propri\'et\'e universelle suivante:\par
Pour toute famille analytique $(X_{s})_{s\in S}$ de $n$ cycles
compacts de $Z$ param\'etr\'ee par un espace complexe  r\'eduit $S$,
il existe un unique morphisme  $\eta :S\longrightarrow {\cal
B}_{n}(Z)$, tel que\par \centerline{ $\eta(s)=
X_{s}$}\bigskip\noindent {\bf 1.0.2.2. Morphisme Douady-
Barlet.}\par\noindent Rappelons que pour tout espace complexe $Z$,
Douady ([D]) a muni l'ensemble des sous
 espaces compacts de $Z$ d'une structure analytique complexe. Il montre qu'il existe un
 couple d'espaces complexes  (${\cal D}_{\#}(Z)$, ${\cal D}(Z)$) et un couple de morphismes propres
   ($\sigma$, $\pi_{{\cal D}}$), $\sigma$ \'etant un plongement et $\pi_{{\cal D}}$ plat, install\'es dans le diagramme commutatif
$\ds{\xymatrix{{\cal D}_{\#}(Z)\ar@/^/[rr]^{\pi_{\cal
D}}\ar[r]_{\sigma}&{\cal D}(Z)\times Z\ar[r]_{p} &{\cal D}(Z)}}$ et
satisfaisant la propri\'et\'e universelle ({\bf D}) :\par\noindent
pour tout  diagramme commutatif d'espaces complexes
$\ds{\xymatrix{X\ar@/^/[rr]^{\pi}\ar[r]_{i}&S\times Z\ar[r]_{p}
&S}}$
 de m\^eme nature que le pr\'ec\'edent, il existe un unique
morphisme $\Theta_{\pi}:S\rightarrow {\cal D}(Z)$ rendant commutatif
le diagramme
 $$\xymatrix{X\ar[r]^{\theta_{\pi}}\ar[d]_{\pi}&{\cal D}_{\#}(Z)\ar[d]^{\pi_{\cal D}}\\
S\ar[r]_{\Theta_{\pi}}&{\cal D}(Z)}$$ dans lequel  $\ds{X\simeq
S\times_{{\cal D}(Z)}{\cal D}_{\#}(Z)}$.\smallskip\noindent En
introduisant la notion de {\it platitude
g\'eom\'etrique}\footnote{$^{(6)}$}{un morphisme propre,  ouvert et
de corang constant $\pi:X\rightarrow S$, avec $S$ r\'eduit et ${\rm
dim}(X)\geq {\rm dim}(S)$, est dit {\it g\'eom\'etriquement plat}
s'il est possible de munir ses fibres ensemblistes (plus
pr\'ecisement leurs composantes irr\'eductibles) de multiplicit\'es
convenables de sorte que la famille $(\pi^{-1}(s))_{s\in S}$
devienne (ou soit induite par ) une famille analytique de cycles
  au sens de [B1]. Cette notion est \`a l'espace de Barlet ce que
   la notion de platitude est \`a l'espace de Douady.}, que nous
   \'etudions en d\'etail dans le \S {\bf (1.3)}, on d\'eduit du th\'eor\`eme de Barlet
  une propri\'et\'e universelle analogue \`a  {\bf (D)} \`a savoir l'existence d'un couple d'espaces
  complexes r\'eduits (et
localement de dimension finie)   (${\cal B}_{\#}(Z)$, ${\cal B}(Z)$)
et d'un couple de morphismes propres ($\sigma$, $\pi_{{\cal B}}$)
avec $\sigma$ plongement  et $\pi_{{\cal B}}$ g\'eom\'etriquement
plat install\'es dans le diagramme commutatif\par\noindent
$\ds{\xymatrix{{\cal B}_{\#}(Z)\ar@/^/[rr]^{\pi_{\cal
B}}\ar[r]_{\sigma}&{\cal B}(Z)\times Z\ar[r]_{p} &{\cal B}(Z)}}$
 et satisfaisant la
propri\'et\'e universelle ({\bf B}) :\par\noindent pour tout
diagramme commutatif d'espaces complexes
$\ds{\xymatrix{X\ar@/^/[rr]^{\pi}\ar[r]_{i}&S\times Z\ar[r]_{p}
&S}}$
 de m\^eme nature que le pr\'ec\'edent et avec $S$ r\'eduit,
il existe un unique morphisme $\Psi_{\pi}:S\rightarrow {\cal B}(Z)$
rendant commutatif le diagramme
$$\xymatrix{X\ar[r]^{\psi_{\pi}}\ar[d]_{\pi}&{\cal B}_{\#}(Z)\ar[d]^{\pi_{\cal B}}\\
S\ar[r]_{\Psi_{\pi}}&{\cal B}(Z)}$$ et dans lequel  $\ds{X\simeq
S\times_{{\cal B}(Z)}{\cal B}_{\#}(Z)}$.\par\noindent
 De plus, si   ${\cal D}_{n}(Z)$ (resp. ${\cal B}_{n}(Z)$)  d\'esigne le
  sous espace de ${\cal D}(Z)$ (resp. de  ${\cal B}(Z)$)  constitu\'e de
  sous espaces de dimension pure $n$ (resp. des cycles analytiques de
  dimension pure $n$), il existe une unique application holomorphe
  appel\'ee {\it{morphisme Douady- Barlet}}
  (cf [B1]) $\rho:  ({\cal D}_{n}(Z))_{red}\rightarrow  {\cal B}_{n}(Z)$
  associant \`a $Y$ sous espace compact de $Z$, la somme
  $\ds{\sum_{i}n_{i}Y_{i}}$, o\`u $Y_{i}$ d\'esignent
   les composantes irr\'eductibles de $Y_{red}$.\smallskip\par\noindent
{\bf 1.0.2.3. Remarques.}\par\noindent $\bullet$ La propri\'et\'e
universelle de l'espace de Douady dit qu 'un morphisme propre
d'espaces complexes $f:X\rightarrow S$ est plat si et seulement si
le graphe de $f$, que l'on notera  ${\cal G}_{f}$, s'installe dans
le diagramme commutatif
 $$\xymatrix{{\cal G}_{f}\ar[r]^{\theta}\ar[d]_{\pi}&{\cal D}_{\#}(X)\ar[d]^{\pi_{\cal D}}\\
S\ar[r]_{\Theta}&{\cal D}(X)}$$ avec  $\ds{{\cal G}_{f}\simeq
S\times_{{\cal D}(X)}{\cal D}_{\#}(X)}$\par\noindent $\bullet$ Dans
le cadre de l'espace des cycles,  cela nous sugg\`ere d'appeler {\it
g\'eom\'etriquement plats}
 les morphismes d'espaces complexes $f:X\rightarrow S$ \`a base r\'eduite dont le graphe ${\cal G}_{f}$
 s'ins\`ere dans un diagramme commutatif analogue au pr\'ec\'edent obtenu en rempla\c cant ${\cal D}$ par  ${\cal B}$
 et v\'erifiant  $\ds{{\cal G}_{f}\simeq S\times_{{\cal B}(X)}{\cal B}_{\#}(X)}$\par\noindent
$\bullet$ Par construction, on a ${\cal B}_{n}(Z):= {\cal
B}_{n}(Z_{red})$ et que ${\cal B}_{0}(Z):=\ds{\coprod_{r\geq 0}{\rm
Sym}^{r}(Z_{red})}$\smallskip\par\noindent {\bf 1.0.3. Morphisme
trace.}\smallskip\noindent {\bf 1.0.3.0. Faisceaux fondamentaux sur
un espace complexe r\'eduit et de dimension pure.}\par\noindent Sans
entrer dans les d\'etails, signalons que sur un espace analytique
r\'eduit de dimension pure quelconque $X$, il existe, en dehors du
faisceau des formes holomorphes (pouvant \^etre d\'efini de
plusieurs mani\`eres d'ailleurs !) deux faisceaux particuli\`erement
interessants \`a savoir $\omega^{k}_{X}$ et ${\cal L}^{k}_{X}$. Le
premier \'etant le dualis\'e de Andr\'eotti-Kaas- Golovin du
faisceau $\Omega^{n-k}_{X}$ enti\`erement caract\'eris\'e par la
propri\'et\'e de la trace et coincidant, en degr\'e maximal, avec le
faisceau dualisant de Grothendieck (cf \S{\bf (3.1.1)},p.56),
l'autre \'etant constitu\'e de formes m\'eromorphes se prolongeant
holomorphiquement sur toute d\'esingularis\'ee de $X$ et coincidant,
en degr\'e maximal, avec le faisceau des formes m\'eromorphes de
carr\'es int\'egrables introduit par Griffiths, Grauert-
Riemmenschneider ([Gri], [G-R])\footnote{$^{(7)}$}{Les dualis\'es de
Andreotti-Kaas-Golovin sont d\'efinis pour $X$ non n\'ecessairement
r\'eduit. Les hypoth\`eses de pure dimensionnalit\'e et de
r\'eduction permettent, respectivement,  une caract\'erisation par
la trace et une interpr\'etation au moyen des formes m\'eromorphes
sur $X$.}. Il se trouve que, dans le cas des
$V$-vari\'et\'es\footnote{$^{(8)}$}{Ce r\'esultat est encore vraie
dans le cas des singularit\'es rationnelles [E].} et donc, en
particulier, sur  ${\rm Sym}^{k}(B)$, ces deux faisceaux coincident.
Mais il est montr\'e dans [K1] que le faisceau  ${\cal L}^{k}_{X}$
est dot\'e de nombreuses propri\'et\'es fonctorielles parmi
lesquelles la stabilit\'e par image r\'eciproque de morphismes
analytiques (qui fait d\'efaut \`a $\omega^{k}_{X}$); les
restrictions aux sous-espaces quelconques auront toujours un sens et
conservent leur nature originelle. De plus, pour tout morphisme
propre et
 g\'en\'eriquement fini d'espaces analytiques complexes $f:X\rightarrow
Y$, l'image directe au sens des courants induit un morphisme trace
$f_{*}\omega^{k}_{X}\rightarrow \omega^{k}_{Y}$ se restreignant
lui-m\^eme en une trace $f_{*}{\cal L} ^{k}_{X}\rightarrow {\cal
L}^{k}_{Y}$.\smallskip\par\noindent{\bf 1.0.3.1. Trace absolue.}
\smallskip\noindent{\bf 1.0.3.1.0.  Trace absolue associ\'ee \`a un
sous ensemble analytique d'une vari\'et\'e de Stein.}\par\noindent
Avant d'entrer dans le vif du sujet, rappelons que, si
$\pi:X\rightarrow Y$ est  un morphisme propre et surjectif d'espaces
analytiques complexes  de dimension pure $n+r$ et $r$ avec $X$
d\'enombrable \`a l'infini et ${\bf D}^{\bullet,\bullet}_{X}$ (resp.
${\bf D}^{\bullet,\bullet}_{Y}$)  d\'esignent le faisceau des
courants de type $(\bullet, \bullet)$ sur $X$ (resp. $S$), on  d'une
{\it image directe} au sens des courants
\par \centerline{${\rm T}^{\bullet, \bullet}_{f}: f_{*}{\bf
D}^{n+\bullet,n+\bullet}_{X}\rightarrow {\bf
D}^{\bullet,\bullet}_{Y}$}\noindent donnant, en particulier, pour
chaque section $\phi$ de ${\cal A}^{n+\bullet, n+\bullet}_{X}$
(faisceau des formes ind\'efiniment diff\'erentiables de type
$(n+\bullet, n+\bullet)$ sur $X$ et d\'efinies par plongement
local),
$$\langle{\rm T}^{\bullet,\bullet}_{f}(\phi), \psi\rangle:=
\int_{X}{\phi\wedge f^{*}\psi},\,\,\,\,\forall\,\psi\in {\cal
A}^{r-\bullet, r-\bullet}_{c}(Y)$$
Forsmis la compatibilit\'e aux inclusions ouvertes, cette correspondance ne jouit d'aucune propri\'et\'e fonctorielle naturelle en les arguments $X$ et $Y$.\smallskip\noindent
Si $X$ est un sous ensemble analytique de dimension pure $n$  d'un certain domaine d'un espace num\'erique ${\Bbb C}^{N}$, le th\'eor\`eme de
param\'etrisation locale ou de pr\'eparation de Weierstrass permet de le  r\'ealiser localement comme un rev\^etement ramifi\'e au dessus d'un certain polydisque
ouvert de ${\Bbb C}^{n}$. Plus pr\'ecisemment,
pour tout point $x$ de $X$, il existe un voisinage ouvert $V$ de
$x$, un plongement $\sigma:V\rightarrow W\subset {\Bbb C}^{n+p}$
o\`u $W:=U\times B\subset {\Bbb C}^{n}\times  {\Bbb C}^{p}$ ( $U$ et
$B$ \'etant des polydisques ouverts relativement compacts de ${\Bbb
C}^{n}$ et ${\Bbb C}^{p}$ respectivement) tel que  $\sigma(V)\cap(\bar{U}\times \partial{B})=\empty$ de sorte que la projection
canonique sur $U$ induise un morphisme fini et surjectif
$f:V\rightarrow U$ qui est un rev\^etement ramifi\'e au sens usuel dont on notera  $k$ son degr\'e et  $(f_{j})_{1\leq j\leq k}$
 ses branches locales (ou feuillets).  Par abus de notation, on notera encore  $X$ l'ouvert $V$ et d\'esignerons par
 $j:{\rm Reg}(X)\rightarrow X$ l'inclusion
naturelle de la partie lisse de $X$ dans $X$. Dans cette situation,
$f$ induit  une image directe au sens des courants qui s'exprime, en
dehors de la ramification, sous la forme\par \centerline{$\ds{{\rm
T}^{\bullet,\bullet}_{f}(\phi):= \sum^{j=k}_{j=1}f^{*}_{j}(\phi)}$}
\noindent En particulier, pour toute $q$- forme holomorphe $\phi$
sur $X$, $\ds{{\rm T}^{q,0}_{f}(\phi)}$ (trace que nous noterons
dor\'enavant $\ds{{\rm T}^{q}_{f}}$  d\'efinit un courant
$\bar\partial$- ferm\'e de type $(q,0)$. Mais $U$ \'etant lisse (et
de Stein), le lemme de Dolbeault- Grothendieck assure, alors, que
c'est une $(q,0)$-forme holomorphe sur $U$.\smallskip\noindent {\bf
D\'efinition: Propri\'et\'e de la trace absolue.}\par\noindent Soit
$x$ un point de $X$ et $\xi$ un germe en $x$ du faisceau
$j_{*}j^{*}\Omega^{\bullet}_{X}$. On dit que $\xi$ v\'erifie la {\it
propri\'et\'e de la trace absolue} si  pour chaque germe de
param\'etrisation locale de $X$ en $x$  et chaque germe $\alpha$  du
faisceau des  formes  holomorphes $\Omega^{n-\bullet}_{X}$ en $x$,
la forme $\ds{{\rm T}^{n}_{f}(\xi\wedge
\alpha):=\sum^{j=k}_{j=1}f_{j}^{*}(\xi\wedge \alpha)}$, d\'efinie
(et holomorphe) en dehors de la ramification, se prolonge
analytiquement \`a $U$ tout entier.\par\noindent On dit qu'une
section  $\xi$ du faisceau $j_{*}j^{*}\Omega^{\bullet}_{X}$, sur un
ouvert $V$ de $X$, v\'erifie la propri\'et\'e de la trace absolue si
son germe en chaque point de $V$ la v\'erifie.  \smallskip\noindent
L'ensemble des sections de $j_{*}j^{*}\Omega^{\bullet}_{X}$
v\'erifiant cette propri\'et\'e constitue un sous faisceau
coh\'erent du faisceau des formes m\'eromorphes sur $X$ et souvent
not\'e $\omega^{\bullet}_{X}$(cf {\bf(3.5)}).\smallskip\par\noindent
{\bf 1.0.3.1.1.  Trace associ\'ee \`a un cycle et trace
universelle.}
\par\noindent Conservons les notations pr\'ec\'edentes et consid\'erons un
$n$- cycle $X=\sum_{i}n_{i}X_{i}$  de $U\times
B$ dont le support est un rev\^etement ramifi\'e de degr\'e total
$k$ sur $U$ et de branches locales $(f_{j})_{1\leq j\leq k}$. Comme
il a \'et\'e dit dans {\bf(1.0.1.1)},  il lui est naturellement
associ\'e une application analytique classifiante
$F_{X}:U\rightarrow {\rm Sym}^{k}(B)$ que l'on peut ins\'erer dans
le diagramme commutatif
$$\xymatrix{&U\times B\ar[r]^{F_{X}\times Id_{B}}\ar[dd]&{\rm Sym}^{k}(B)
 \times B\ar[dd]&\\
X\ar[ur]\ar[rd]_{f}\ar[rrr]&&& {\rm Sym}^{k}(B) \#
B\ar[ld]^{f_{\#}}\ar[ul]
\\&U\ar[r]_{F_{X}}&{\rm Sym}^{k}(B)&}$$
dans lequel, on a ensemblistement $\ds{(F_{X}\times
Id_{B})^{-1}({\rm Sym}^{k}(B)\# B) = X}$ et $f_{\#}$ un morphisme
fini et  g\'eom\'etriquement plat sur ${\rm Sym}^{k}(B)$ qui est
normal.\par\noindent{\bf(i) Trace de fonctions
holomorphes.}
Si $h$ est une fonction holomorphe sur $U\times
B$ que l'on peut supposer $\sigma_{k}$-invariante (sinon on la
sym\'etrise !), on d\'efinit une fonction holomorphe  $\tilde{h}$
sur $U\times B^{k}$ en posant $
\tilde{h}(t,x_{1},\cdots,x_{k}):=\sum_{j}h(t,x_{j})$. Comme elle est
$\sigma_{k}$-invariante, elle correspond \`a une unique application
holomorphe $\bar{h}$ sur $U\times {\rm Sym}^{k}(B)$. On convient,
alors, de la d\'efinition
$${\cal T}^{0}_{f}(h):= {\bar{h}}(t,F_{X}(t))= F^{*}_{X}({\goth T}^{0}_{f_{\#}}(\bar{h}))$$
donn\'ee, en dehors de la ramification, par l'expression
$\ds{\sum_{i}n_{i}\sum_{j}h(t,f_{j}(t))}$.
\par\noindent {\bf(ii) Trace de formes holomorphes.}
Pour d\'efinir notre morphisme trace sur les formes holomorphes, on
 utilise la trace universelle $\ds{{\goth T}^{\bullet}_{f_{\#}}:
{f_{\#}}_{*}{\omega^{\bullet}_{\#}} \rightarrow
\omega^{\bullet}_{{\rm Sym}^{k}({\Bbb C}^{p})}}$ et les formes de
Newton. Comme un tel morphisme est n\'ecessairement ${\cal
O}_{U}$-lin\'eaire, on s'interesse essentiellement aux formes
provenant de $B$. Soient  $w_{i,j}$ une $j$-forme de Newton sur $B$
et $\ds{\tilde{w}_{i,j}:=\sum_{r}p_{r}^{*}w_{i,j}}$ o\`u $p_{r}$
d\'esigne simplement la projection canonique de $B^{k}$ sur le
$r$-\`eme facteur. Comme $\tilde{w}_{i,j}$ est
$\sigma_{k}$-invariantes, elle provient forme m\'eromorphe
$\bar{w}_{i,j}$  section  du faisceau $\omega^{\bullet}_{{\rm
Sym}^{k}(B)\times B}$ (qui est d'ailleurs engendr\'e par ces formes
cf {\bf(1.0.3.1)}. On pose, alors
$${\cal T}^{j}_{f}(w_{i,j}):= F^{*}_{X}({\goth
T}^{j}_{f_{\#}}(\bar{w}_{i,j})$$ Alors, le proc\'ed\'e de
sym\'etrisation et la d\'ecomposition nucl\'eaire des formes sur
$U\times B$, permettent d'\'etendre cette relation aux formes
holomorphes quelconques sur $U\times B$ et de d\'efinir un morphisme
trace ${\cal
T}^{\bullet}_{f}:f_{*}\Omega^{\bullet}_{X}\rightarrow\Omega^{\bullet}_{U}$
qui est ${\cal O}_{U}$-lin\'eaire, commutant \`a la diff\'erentielle
ext\'erieure usuelle, de nature local sur $U$ et $B$ et compatible
\`a l'additivit\'e
 des cycles en faisant la somme point par point des cycles correspondants
  gr\^ace \`a l'application ``somme'' donn\'ee dans {\bf(1.0.1.1)}.
  \smallskip\par\noindent
 {\bf 1.0.3.2. Trace relative.}\par\noindent Pour mieux positionner
  le probl\`eme, commen\c cons par \par\noindent {\bf 1.0.3.2.0.
Quelques remarques d'ordre g\'en\'eral.}\par\noindent Il est bien
connu que tout morphisme ouvert \`a fibres de dimension pure $n$,
$\pi:X\rightarrow S$, entre espaces complexes, admet une
factorisation locale, par rapport \`a l'une quelconque de ses
fibres, du type  $(\clubsuit)$
$$\xymatrix{X\ar[d]^{\pi} \ar [dr]_{f}\ar [r]^{\!\!\!\!\!\!\!\!\!\!\!\!\!\sigma}
 & Z:= S\times U\times B\ar[d]^{p}\\
 S&\ar[l]^{q} S\times U}$$
 o\`u  $U$ et $B$ sont des polydisques relativement compacts de ${\Bbb
C}^{n}$ et ${\Bbb C}^{p}$ respectivement, $f$ est ouvert, fini et surjectif,
$\sigma$ un plongement local, $p$ et $q$ les projections
canoniques.\par\noindent
L'enjeu est de voir, dans quelle mesure, il est possible de produire (par analogie avec le cas absolu!) un morphisme $\ds{{\cal
T}^{\bullet}_{f}:f_{*}\Omega^{\bullet}_{X/S}\rightarrow\Omega^{\bullet}_{S\times
U/S}}$  v\'erifiant naturellement  les
propri\'et\'es suivantes:\par {\bf(i)} il est de nature locale sur
$U$ et $S$\par {\bf(ii)} il est   $\Omega^{\bullet}_{S\times
U/S}$-lin\'eaire et commute aux diff\'erentielles $S$-relatives  usuelles.\par
{\bf(iii)} il  est   additive en les pond\'erations sur  $\pi$ dans
le sens suivant:
\smallskip\noindent
soient ${\goth X}_{1}$ et ${\goth X}_{2}$ deux pond\'erations de
$\pi$ et $(U,B,\sigma)$ une \'ecaille $S$-adapt\'ee, alors
$${\cal
T}^{\bullet}_{f,{\goth X}_{1}\oplus{\goth X}_{2} }= {\cal
T}^{\bullet}_{f,{\goth X}_{1}}+ {\cal T}^{\bullet}_{f,{\goth
X}_{2}}$$
On peut tenter d'utiliser l'image directe des courants mais alors on se retrouve devant de s\'erieues difficult\'es. En effet, {\bf (1.0.3.1)} donne,
en particulier, un morphisme {\it{trace}} ${\cal
T}^{\bullet,0}_{f}:f_{*}i^{*}(\Omega^{\bullet}_{Z/S})\rightarrow
{\cal D}^{\bullet,0}_{S\times U}$, \`a valeurs dans le faisceau des
courants $\bar\partial$- ferm\'es de type $(\bullet,0)$. Cette
trace est bien \'evidemment $\Omega^{\bullet}_{S\times U}$-
lin\'eaire, g\'en\'eriquement holomorphe sur $S\times U$ et
holomorphe fibre par fibre en vertu du lemme de Dolbeault-
Grothendieck. Mais des exemples simples montrent qu'en   pr\'esence de
singularit\'es arbitraires sur $S$ et sans aucune condition sur le
mode de variance des fibres du morphisme, il n'y'a pratiquement
aucune chance de voir ce courant d\'efinir une forme holomorphe
globale $S$-relative sur $S\times U$. En fait, le
passage de la situation absolue \`a la situation relative n'est pas
du tout une simple formalit\'e puisqu'il fait apparaitre  plusieurs
 obstacles de tailles \`a savoir :\par
 {\bf(1)} un morphisme fini et surjectif sur une  base quelconque n'a aucune raison de d\'efinir un rev\^etement ramifi\'e au sens usuel puisque le degr\'e peut ne pas \^etre localement constant.\par
{\bf(2)} Combien m\^eme ce serait le cas, et disposant ainsi d'une
application analytique classifiante  $F:S\times
U\rightarrow\rm{S}ym^{k}(B)$, la variation de la ramification
\'etant incontr\^olable,  on ne peut emp\^echer qu'il y'ait
certaines valeurs critiques   du param\`etre  $s$ pour lesquelles
l'image de $U$ par la restriction de $F$ \`a $\{s\}\times U$  soit
enti\`erement contenue dans le lieu singulier de $\rm{S}ym^{k}(B)$.
Or, par construction, la trace relative (si elle existe !)  doit
\^etre un objet qui, pour chaque $s$ fix\'e, coincide avec la trace
absolue pr\'ec\'edemment d\'ecrite. Mais, alors,  pour les valeurs
critiques du param\`etre quel sens donn\'e \`a la trace absolue
?\par {\bf (3)} On ne dispose pas d'un analogue du lemme de
Dolbeault- Grothendieck dans cette situation relative g\'en\'erale.
\smallskip\noindent Pour surmonter {\bf(1)}, on impose \`a la
famille d'\^etre au moins {\it continue} ce qui contraint le
morphisme d'\^etre au moins  ouvert et de corang
constant\footnote{$^{(9)}$}{Mais cela ne suffit pas car m\^eme si le
morphisme est ouvert entre espaces complexes irr\'eductibles et donc
\`a fibres de dimension pure constante, les fibres ne d\'efinissent
pas toujours  une famille continue de cycles; pour s'en convaincre,
il suffit, par exemple, de consid\'erer la  normalisation faible
d'une courbe !}. La condition suppl\'ementaire que doit v\'erifier
un tel
 morphisme est que la fonction
 $s\rightarrow {\rm Deg}_{E}(\pi^{-1}(s))$
 \footnote{$^{(10)}$}{ Plus pr\'ecisement, \'etant donn\'e $s_{0}$ dans $S$ et
  une  projection locale quelconque  induite par la donn\'ee d'une
  \'ecaille adapt\'ee \`a $\pi^{-1}(s_{0})$, il existe un voisinage
  ouvert $S_{0}$ de $s_{0}$ tel que $E$ soit adapt\'ee \`a $\pi^{-1}(s)$
  pour tout $s\in S_{0}$ et
   l'application  $s\rightarrow {\rm Deg}_{E}(\pi^{-1}(s))$ est
    localement constante.} soit localement constante; ce qui est toujours le cas si $S$ est localement irr\'eductible. On  contourne {\bf(2)} en
utilisant une stratification convenable de $\rm{S}ym^{k}(B)$, qui,
bri\`evement dit,  assurera l'existence d'une valeur enti\`ere
minimale $l$ telle que $\{s\}\times U$ ne soit pas enti\`erement
incluse dans le lieu singulier de $\rm{S}ym^{l}(B)$. Malheureusement
(ou heureusement !) {\bf(3)} est insurmontable. Sans conditions
suppl\'ementaires sur la variation de la famille, les courants
$\bar\partial$-ferm\'es sur $S\times U$ obtenus par image directe ne
produisent, g\'en\'eralement jamais, des formes holomorphes
$S$-relatives!\smallskip\noindent Pour les d\'etails sur la
stratification de $\rm{S}ym^{k}(B)$, on  renvoie le lecteur \`a
[B2], p18 dans lequel il est montr\'e que les sous ensembles
alg\'ebriques  $M_{\mu} := \left\{x\in \rm{S}ym^{k}({\Bbb C}^{p}) /
mult(x)\ge{\mu}\right\}$  ont pour lieu singulier $M_{\mu+1}$ et
$M_{\mu}\setminus M_{\mu+1}$ est un sous ensemble analytique
localement ferm\'e sans singularit\'es, dont les composantes
connexes sont en bijection avec les suites d'entiers positifs
$(n_{1},\cdots,n_{i})$ v\'erifiant $\ds{\sum_{j=1}^{i}n_{j}=k}$ et
$\ds{\sum_{j=1}^{i}n_{j}{{(n_{j}-1)}\over 2}=\mu}$.\par\noindent Si
$S$ est   un espace complexe r\'eduit irr\'eductible, $F$ une
application analytique de $S\times U$ dans  $\rm{S}ym^{k}({\Bbb
C}^{p})$ et $l$ le plus grand entier pour lequel $F(S\times
U)\subset M_{l}$. Alors, {\bf {la composante g\'en\'erique}} de $F$
est l'unique composante connexe de $M_{l}\setminus M_{l+1}$
contenant $F(S\times U\setminus F^{-1}(
M_{l+1}))$.\smallskip\noindent
 Comme on va le voir les traces des formes holomorphes
$S$ -relatives sur $X$ sont enti\`erement d\'etermin\'ees par les
traces des formes de Newton induisant  elles m\^emes des sections
d'un  faisceau particulier sur ${\rm Sym}^{k}(B)$.
\smallskip\par\noindent
{\bf 1.0.3.2.1. Trace associ\'ee \`a une famille analytique locale
de cycles et trace  universelle.}\par\noindent Notons $\#:={\rm
Sym}^{k}(B)\# B$ le sous espace d'incidence (ou diagonale
g\'en\'eralis\'ee) qui est muni naturellement d'une projection finie
et surjective ${\tilde f}:\#\rightarrow {\rm Sym}^{k}(B)$. Comme
${\rm Sym}^{k}(B)$ est un espace complexe normal donc localement
irr\'eductible, ce morphisme est un rev\^etement ramifi\'e au sens
usuel. Toujours en raison de la normalit\'e, ${\tilde f}$ est
analytiquement g\'eom\'etriquement plat. Comme dans le cas absolu,
on dispose d'un morphisme trace fondamental ou universel
$${\goth T}^{\bullet}_{\#}: {\tilde{f}}_{*}{\omega^{\bullet}_{\#}}
\rightarrow \omega^{\bullet}_{{\rm Sym}^{k}({\Bbb C}^{p})}$$
qui est un morphisme d'alg\`ebres diff\'erentielles
gradu\'ees.\par\noindent La situation g\'en\'erale
locale se traite via le cas universel si une certaine application analytique $F:S\times U\rightarrow {\rm Sym}^{k}(B)$ et donn\'ee. Dans ce cas, on se connecte \`a la situation universelle par le biais du  diagramme commutatif
$$\xymatrix{&S\times U\times B\ar[r]^{F\times Id_{B}}\ar[dd]&{\rm Sym}^{k}(B)
 \times B\ar[dd]&\\
X\ar[ur]\ar[rd]_{f}\ar[rrr]&&& {\rm Sym}^{k}(B) \# B\ar[ld]\ar[ul]
\\&S\times U\ar[r]_{F}&{\rm Sym}^{k}(B)&}$$
dans lequel, on a ensemblistement $\ds{(F\times Id_{B})^{-1}({\rm
Sym}^{k}(B)\# B) = X}$. Notons $\Gamma$ la  composante g\'en\'erique
de $F$ et $S_{\Gamma}\times U$ son image r\'eciproque dans le
changement de base $\Gamma\rightarrow {\rm Sym}^{k}(B)$. Elle  est
\'evidemment  d'int\'erieure vide et connexe puisque $S$ est
suppos\'e irr\'eductible. Alors, les changements de bases
$S_{\Gamma}\times U\rightarrow S\times U$ et $\Gamma\rightarrow {\rm
Sym}^{k}(B)$, donnent un diagramme commutatif dont les carr\'es
extr\^emes sont cart\'esiens
$$\xymatrix{&X\ar[r]^{\tilde{F}}\ar[d]^{f}&\#\ar[d]^{\tilde{f}}&\\
X_{\Gamma}\ar[ur]\ar[d]_{f_{\Gamma}}&S\times U\ar[r]_{F}&{\rm
Sym}^{k}(B)&
\# _{\Gamma}\ar[d]^{\tilde{f}_{\Gamma}}\ar[ul]\\
S_{\Gamma}\times U\ar[ru]\ar[rrr]_{F_{\Gamma}}&&&\Gamma\ar[lu]}$$
On suppose implicitement donn\'e un plongement de $X$ dans $S\times U\times B$.  Soit $\sigma$ une forme de Newton sur $B$ dont on note  $\tilde{\sigma}$
 son image r\'eciproque sur $B^{k}$ (comme dans {\bf(1.0.3.1.1)}). Cette
derni\`ere provient d'une certaine forme m\'eromorphe $\bar{\sigma}$ section du faisceau ${\cal L}^{\bullet}_{{\rm
Sym}^{k}(B)\times B}$. On d\'efinit, alors,  la trace relative \`a $f$,
$${\cal T}^{\bullet,0}_{f}(\sigma):= {\cal
T}^{\bullet,0}_{f_{\Gamma}}(\sigma|_{X_{\Gamma}})=
F_{\Gamma}^{*}({\goth
T}^{\bullet,0}_{\tilde{f}_{\Gamma}}(\tilde\sigma|_{\#_{\Gamma}}))$$
qui a un sens puisque $\sigma|_{X_{\Gamma}}$ est une forme de Newton
provenant de la restriction $\tilde{\sigma}$ \`a $\#_{\Gamma}$
d\'efinissant une section du faisceau  ${\cal
L}^{\bullet}_{\#_{\Gamma}}$.
\par\noindent Pour des raisons d\'ej\`a invoqu\'ees, l'action de ce morphisme sur les formes de Newton suffit \`a d\'eterminer son action sur les formes holomorphes relatives g\'en\'erales pour donner un morphisme trace relative
$\ds{{\cal
T}^{\bullet}_{f}:f_{*}\Omega^{\bullet}_{X/S}\rightarrow\Omega^{\bullet}_{S\times
U/S}}$  v\'erifiant les propri\'et\'es {\bf(i)}, {\bf(ii)} et {\bf(iii)} pr\'ecis\'ees dans {\bf(1.0.3.2.0)}.
\smallskip\par\noindent
 {\bf 1.0.3.3. Degr\'e
g\'eom\'etrique et alg\'ebrique d'un morphisme fini
(cf[Va]).}\par\noindent
 Soient $X$ et $S$ deux espaces analytiques complexes avec $S$ connexe et  $\pi:X\rightarrow S$  un
morphisme fini  ouvert et surjectif. Alors,
\par\noindent $\bullet$ si $\pi$ est  contin\^ument (resp.
analytiquement) g\'eom\'etriquement plat sur $S$ r\'eduit, il existe
un entier $k$, appel\'e {\it le degr\'e  g\'eom\'etrique} de $\pi$,
et une application continue (resp. analytique)
$\ds{\Psi^{\pi}:S\rightarrow \rm{S}ym^{k}(X_{red})}$ donn\'ee par
$\ds{\Psi^{\pi}(s)=\sum_{x\in \pi^{-1}(s)}\{x\}}$ et tenant compte
des multiplicit\'es. De plus, ces conditions g\'en\`ere un morphisme
{\it trace continue} (resp. {\it holomorphe})   ${\cal
T}^{(c)}_{\pi}: \pi_{*}{\cal C}^{0}_{X}\rightarrow {\cal C}^{0}_{S}$
(resp. ${\cal T}^{h}_{\pi}: \pi_{*}{\cal O}_{X}\rightarrow {\cal
O}_{S}$).\smallskip\noindent $\bullet$ si $\pi$ est plat sur $S$
arbitraire, le {\it degr\'e alg\'ebrique} de $\pi$ est d\'efini
comme \'etant le plus grand entier  $r$ tel que le faisceau
coh\'erent $\pi_{*}({\cal O}_{X})$ soit un ${\cal O}_{S}$-module
localement libre de rang $r$.\par\noindent On constate que le
degr\'e alg\'ebrique $r$  est stable par changement de base et qu'il
coincide avec le degr\'e g\'eom\'etrique $k$ si $S$ est
r\'eduit.\par\noindent
 La situation alg\'ebriquement plate  produit aussi une {\it trace holomorphe} naturelle
${\cal T}^{(h)}_{\pi}: {\pi_{*}{\cal O}_{X}}\rightarrow {\cal
O}_{S}$ qui, gr\^ace au changement de base  $S_{red}\rightarrow S$
et le diagramme commutatif qui s'en d\'eduit
$$\xymatrix{X':= S_{red}\times_{S} X\ar[r]\ar[d]_{\pi'}&X\ar[d]^{\pi}\\
S':=S_{red}\ar[r]&S}$$ est compatible avec les traces
pr\'ec\'edentes en un sens \'evident. On pose, d'ailleurs,  ${\cal
T}^{(c)}_{\pi}:= {\cal T}^{(c)}_{\pi'}$ pr\'ec\'edentes en un sens
\'evident. On pose, d'ailleurs,  ${\cal T}^{(c)}_{\pi}:= {\cal
T}^{(c)}_{\pi'}$\smallskip\par\noindent {\bf 1.1. Quelques
exemples.}\smallskip\noindent {\bf(i)} ([Ang]). ${\rm Sym}^{2}({\Bbb
C}^{2})$ est le quotient de $({\Bbb C}^{2})^{2}$ par
$\sigma_{2}=\{Id, -Id\}$ le groupe sym\'etrique d'ordre $2$.
 Soient $(u,v)$ et $(u',v')$ deux couples
de $({\Bbb C}^{2})^{2}$. Alors les applications $N_{1}=u+u'$,
$N_{2}=v+v'$, $N_{11}=u^{2}+u'^{2}$, $N_{12}=uv+u'v'$ et
$N_{22}=v^{2}+v'^{2}$ fournissent un plongement de ${\rm
Sym}^{2}({\Bbb C}^{2})$ dans ${\Bbb C}^{5}$.\smallskip\noindent Soit
$S:=\lbrace{(x,y,z)\in {\Bbb C}^{3}/ xy - z^{2}= 0}\rbrace$ le
c\^one de ${\Bbb C}^{3}$ et  $\pi:X:={\Bbb C}^{2}\rightarrow S$ la
param\'etrisation donn\'ee par $(u,v)\rightarrow (u^{2}, v^{2}, uv)$
qui est un morphisme ouvert, fini, surjectif ( rev\^etement
ramifi\'e de degr\'e $2$), Cohen-Macaulay mais non plat. On a un
morphisme trace
$${\cal T}^{0}_{\pi}:\pi_{*}{\cal O}_{X}\rightarrow {\cal O}_{S}$$
donn\'e par
$${\cal T}^{0}_{\pi}(f)(s):=f(x_{1}(s)) + f(x_{2}(s)),\,\,\,\,\,\forall\,f\in \Gamma(U,{\cal
O}_{X})$$ o\`u les points $x_{i}(s)$ constituent  la fibre
$\pi^{-1}(s)$.\smallskip\noindent Comme les fonctions holomorphes
sur $S$ s'expriment en fonctions de $u$, $v$, $u^{2}$, $v^{2}$ et
$uv$, on voit clairement comment interviennent naturellement les
fonctions de Newton. Dans ce cas, on a
$${\cal T}^{0}_{\pi}(u)= {\cal T}^{0}_{}(v)=0,\,\,{\cal T}^{0}_{\pi}(u^{2})=2x
,\,\,{\cal T}^{0}_{\pi}(v^{2})=2y,\,\,{\cal T}^{0}_{\pi}(uv)=2z$$ Et
le morphisme  $S\rightarrow {\rm Sym}^{2}({\Bbb C}^{2})$
d\'efinissant la famille correspond exactement au morphisme
d'anneaux donn\'e par
$$N_{1}\rightarrow{\cal T}^{0}_{\pi}(u),\,N_{2}\rightarrow {\cal T}^{0}_{\pi}(v),\,\,
N_{11}\rightarrow {\cal T}^{0}_{\pi}(u^{2}),\,\,N_{12}\rightarrow
{\cal T}^{0}_{\pi}(uv),\,\,N_{22}\rightarrow {\cal
T}^{0}_{\pi}(v^{2})$$\par\noindent
 {\bf(ii)} (cf [K1]).  Soit
$S=\{(a,b,c)\in{\Bbb C}^{3} :a^2=cb^2\}$, la surface faiblement
normale (donc non localement irr\'eductible), commun\'ement
appel\'ee ``le parapluie de Whitney'', dont le  lieu singulier est
la droite $\Sigma=\{(a,b,c)\in {\Bbb C}^{3} : a=b=0\}$. Soit $D$ un
disque centr\'e en l'origine  et relativement compact dans ${\Bbb
C}$.\par \noindent Notons $\rm{S}ym^2_{0}({\Bbb C}^{2})$, la partie
homog\`ene de degr\'e 2 de ${\rm Sym}^{2}({\Bbb C}^{2})$
canoniquement  identifi\'ee \`a ${\Bbb C}^{2}/ \{id,-id\}$ qui est
isomorphe  au c\^one de ${\Bbb C}^3$ donn\'e dans l'exemple
{\bf(i)}.\par\noindent Soit ${\cal I}:=( u^2 - ct^2,v^2 - b^2, uv -
at, a^{2}- cb^{2})$ l'id\'eal de ${{\Bbb C}^{6}}$ d\'efinissant
 un certain sous espace que l'on  notera abusivement $X=\{(u,v,a,b,c,t)\in {{\Bbb
C}^{2}}\times S\times D : u^2=ct^2, v^2=b^2, uv=at\}$. On
d\'esignera par $\pi: X\rightarrow S$ (resp. $f:X\rightarrow S\times
D$) les morphismes induits par  les projections canoniques de
$S\times D\times{\Bbb C}^2\rightarrow S\times D\rightarrow S$ et par
 $F:S\times D\rightarrow \rm{S}ym^2_{0}({\Bbb C}^{2})$,
 l'application analytique
donn\'ee par $(a,b,c,t)\rightarrow(ct^2,b^2,at).$\smallskip
\noindent Il apparait clairement que  $X$, qui est  le produit
fibr\'e de $S\times D$ par ${\Bbb C}^2$ au dessus de
$\rm{S}ym^2_{0}({\Bbb C}^{2})$, s'identifie, apr\`es quotient par
l'involution $(a,b,c,t,u,v)\rightarrow(a,b,c,t,-u,-v)$, au graphe de
l'application $F$. \par\noindent Pour $s$ g\'en\'erique, la fibre
$X_{s}:=\pi^{-1}(s)$  est constitu\'ee d'un couple de droites en
position g\'en\'erale dont les \'equations sont
 $$\left.\matrix{
 \left\{\matrix{u_1&=&\ds{a\over b}t\cr
             v_1&=&b\cr}\right.& {\rm et }& \left\{\matrix{u_2&=&-u_1\cr
v_2&=&-v_1\cr}\right.\cr}\right.$$ Remarquons au passage que $\pi$
n'est pas plat puisque de multiplicit\'e $3$ en l'origine alors
qu'elle vaut  $2$ en les points  g\'en\'eriques. \smallskip\noindent
  Comme  les fonctions de Newton fondamentales $u^{2}$, $v^{2}$ et $uv$
ont pour traces $S$-relatives respectives ${\cal
T}^{0}_{f}(u^{2})=2ct^{2}$, ${\cal T}^{0}_{f}(v^{2})=2b^{2}$ et
${\cal T}^{0}_{f}(uv)=2at$, qui sont manifestement holomorphes sur
$S$, on a une famille analytique de rev\^etements ramifi\'es de
degr\'e 2 et de codimension 2 de $D$ dans $D\times {\Bbb C}^{2}$,
``rep\'er\'e'' par le morphisme $ F:S\times D\rightarrow
\rm{S}ym^2_{0}({\Bbb C}^{2})$, envoyant chaque couple $(s,t)$ sur
les fonctions sym\'etriques des branches locales
$f_{1}=(u_{1},v_{1})$ et $f_{2}=-f_{1}$. Il s'en suit que $\pi$
d\'efinit une famille continue de $1$-cycles de ${\Bbb C}^{3}$ (on
dira plus bri\`evement que $\pi$ est contin\^ument
g\'eom\'etriquement plat).
\smallskip\noindent
L'analyticit\'e de la famille se teste en regardant la trace
$S$-relative des formes de Newton. Dans notre cas,  le faisceau des
$1$- formes de Newton est engendr\'e par les formes $udu$, $vdv$,
$udv$ et $vdu$ dont les traces $S$-relatives sont ${\cal
T}^{1}_{f}(udu)= 2ctdt$, ${\cal T}^{1}_{f}(vdv)= 0$ et ${\cal
T}^{1}_{f}(udv-vdu)=- 2adt$ qui sont clairement $S$- holomorphes.
Ainsi $\pi$ d\'efinit une famille analytique locale de 1-cycles de
${\Bbb C}^3 $ ( on dira que $\pi$ est analytiquement
g\'eom\'etriquement plat ).\smallskip\noindent Remarquons que
l'identification $\rm{S}ym^2_{0}({\Bbb
C}^{2})\simeq\{(x,y,z)\in{\Bbb C}^{3}/ xy=z^{2}\}$, permet de voir
que la forme $udv-vdu$ correspond exactement (par image directe) \`a
la forme m\'eromorphe (non holomorphe) $\ds{z({dx\over{x}} -
{dy\over{y}}})$ sur le c\^one; elle d\'efinit naturellement une
section du faisceau ${\cal L}^{1}_{\rm{S}ym^2_{0}({\Bbb
C}^{2})}$.\bigskip\noindent
 {\bf (iii). Exemple
d'une famille continue mais non analytique
([B1]).}\smallskip\noindent La construction de l' exemple se fait de
la fa\c con suivante:\par\noindent Consid\'erons $(t,x_{1}, x_{2},
x_{3}, y_{1},y_{2},y_{3}, z_{1}, z_{2}, z_{3})$ un syst\`eme de
coordonn\'ees sur ${\Bbb C}^{10}$ et
posons\smallskip\bigskip\noindent \centerline{$X(t)=x_{1}t^{2} +
x_{2}t + x_{3}$}\par\noindent \centerline{$Y(t)=y_{1}t^{2} + y_{2}t
+ y_{3}$}\par\noindent \centerline{$Z(t)=z_{1}t^{2} + z_{2}t +
z_{3}$}\smallskip\noindent On impose la condition
$X(t)Y(t)=Z^{2}(t)$ de sorte \`a pouvoir \'ecrire\smallskip\noindent
\centerline{$ u^{2} = (at+b)^{2}= X(t)$}\par\noindent
\centerline{$v^{2} =(ct + d)^{2} = Y(t)$}\par\noindent
\centerline{$uv= (at +b)(ct +d)= Z(t)$}\smallskip\noindent On trouve
alors les relations suivantes :\smallskip\noindent
$x^{2}_{2}=4x_{1}x_{3},\,\,\, z^{2}_{1}=x_{1}y_{1},\,\,
y^{2}_{2}=4y_{1}y_{3},\,\,
z^{2}_{3}=x_{3}y_{3},\,\,4z_{1}z_{3}=x_{2}y_{2},\,\,x_{1}y_{2} +
x_{2}y_{1} =2z_{1}z_{2}$ $2z_{2}z_{3}=x_{2}y_{3} +
x_{3}y_{2},\,\,x_{2}z_{2}=2(x_{1}z_{3}+x_{3}z_{1}),\,\,y_{2}z_{2}=2(y_{1}z_{3}+y_{3}z_{1}),\,\,
x_{1}y_{3} + x_{2}y_{2} + x_{3}y_{1} =2z_{1}z_{3} + z^{2}_{2}$ et
\par\noindent
$x_{1}= a^{2},\,\,x_{2}= 2ab,\,\,x_{3}= b^{2},\,\,y_{1}=
c^{2},\,\,y_{2}= 2cd,\,\,y_{3}= d^{2},\,\, z_{1}= ac,\,\,z_{2}= ad +
bc,\,\,z_{3}= bd,$\par\noindent que nous symboliserons  par
($\spadesuit$) et ($\clubsuit$) respectivement.\smallskip\noindent
Ces relations d\'efinissent une application holomorphe  $\phi:{\Bbb
C}^{4}\rightarrow {\Bbb C}^{9}$ dont l'image s'identifie au c\^one
de dimension $4$ de ${\Bbb C}^{9}$. Notons $S$ le sous espace
r\'eduit (et irr\'eductible puisque c'est l'image de ${\Bbb C}^{4}$)
associ\'ee \`a l'id\'eal d\'efini par les relations $(\spadesuit)$.
Alors, il est facile de voir que  la fonction m\'eromorphe et
localement born\'ee  $\ds{z_{1}x_{2}\over{x_{1}}}$ se prolonge
contin\^ument mais non holomorphiquement  \`a $S$ tout
entier.\par\noindent Soit $D$ le disque unit\'e relativement compact
de ${\Bbb C}$,$\Gamma$ le graphe de la famille de droites
pr\'ec\'edemment d\'ecrite et $\pi:\Gamma\rightarrow S\times D$ le
morphisme fini induit par la projection canonique. Si la famille
\'etait analytique, $\pi$ d\'efinirait une famille analytique de
points param\'etr\'ee par $S$ et par cons\'equent la trace $S$-
relative  de toute forme holomorphe serait $S$- holomorphe. Mais il
n'en est rien puisque l'on a
$${\cal T}r_{\pi}(udv)=2(ct +d)d(at +b) = 2cat.dt + 2da.dt $$
\noindent Or, en vertu de $(\spadesuit)$ et $(\clubsuit)$, on a  sur
$S$,
$$ac = z_{1},\,\,\,\,\,  ad=z_{2} - {z_{1}x_{2}\over{2x_{1}}}$$
\noindent Et comme la fonction m\'eromorphe
$\ds{z_{1}x_{2}\over{x_{1}}}$ ne se prolonge pas holomorphiquement
\`a $S$, cette famille continue ne peut \^etre
analytique.\smallskip\noindent En choisissant convenablement les
param\`etres, il est tout \`a fait possible de transformer $S$ en
espace  r\'eduit, non r\'eduit, Cohen Macaulay, non Cohen Macaulay
et m\^eme intersection compl\`ete. \par\noindent Soit $\xi\in {\Bbb
C}$. Alors, en posant \par\noindent
$$a=\xi^{2}+\xi^{3},\,\,b=\xi^{5},\,\,c=\xi^{6}-\xi^{7},\,\,d=\xi^{9}+\xi^{11}$$
on obtient le syst\`eme
$$\matrix{ x_{1}=\xi^{4}(1+\xi)^{2}&\!\! y_{1}=\xi^{12}(1-\xi)^{2}& z_{1}=\xi^{8}(1-\xi^{2})\cr
x_{2}=2\xi^{7}(1+\xi)&\,\,\,\,\,\,\,\,\,\,\,\,\,\,\,\,\,\,\,
y_{2}=2\xi^{15}(1-\xi)(1+\xi^{2})&\,\,\,\,\,\,\,\,\,\,\,\,\,
z_{2}=\xi^{11}(2+\xi^{2}+\xi^{3}\cr
\!\!\!\!\!\!\!\!\!\!\!\!\!\!\!\!\!\!\!\!
x_{3}=\xi^{10}&y_{3}=\xi^{18}(1+\xi^{2})^{2}&\,\,\,z_{3}=\xi^{14}(1+\xi^{2})}$$
d\'ecrivant  une courbe monomiale  ${\cal C}$  inscrite  sur le
c\^one pr\'ec\'edent et sur laquelle la fonction m\'eromorphe
$\ds{\sigma:= {z_{1}x_{2}\over{x_{1}}}} =
\ds{z_{3}y_{2}\over{y_{3}}}$ v\'erifie une \'equation de
d\'ependance int\'egrale puisque $\ds{\sigma^{2}=
{{x_{3}y_{1}\over{2}}}}$ mais ne se prolonge pas analytiquement sur
${\cal C}$.  Pour s'en convaincre, il suffit de consid\'erer son
image r\'eciproque par cette param\'etrisation (qui est la
normalisation) et de constater que, pour des raisons de degr\'es et
de valuations, la fonction $\ds{{\tilde{\sigma}} =
2\xi^{11}(1-\xi)}$ ne peut \^etre holomorphe sur cette
courbe.\smallskip\par\noindent {\bf 1.2. Sur quelques notions
fondamentales.}\par\noindent Les notions principales  dont il s'agit
seront celles de dimension d'un espace analytique complexe,
\'equidimensionnalit\'e et ouverture d'un morphisme d'espaces
complexes. Nos principales sources auxquelles nous renvoyons le
lecteur pour de plus amples d\'etails, sont [A.S], [G.R1], [Gr2],
[Ka], [Fi] ou [Lo], pour l'aspect analytique complexe et,
\'evidemment, \`a [E.G.A. 4] (et principalement au paragraphe  14.4)
pour l'aspect alg\'ebrique. Il pourra d'ailleurs constater la
validit\'e  d'un bon nombre de r\'esultats alg\'ebriques dans le
cadre de la g\'eom\'etrie analytique complexe. \smallskip\noindent
{\bf 1.2.0.} Il nous arrivera d'utiliser, sans mentionner, les
caract\'erisations suivantes  de l'irr\'eductibilit\'e locale, la
faible normalit\'e ou normalit\'e d'un espace complexe r\'eduit. Si
$X$ est un tel espace, de faisceau structural ${\cal O}_{X}$, on
d\'esigne par ${\cal O}^{c}_{X}$ (resp. ${\cal O}^{0}_{X}$ ) le
faisceau des germes de fonctions g\'en\'eriquement holomorphes et
continues sur $X$ (resp. celui des germes de fonctions
g\'en\'eriquement holomorphes et localement born\'ees). Alors, on a
(cf [Mo] par exemple)\smallskip {\bf (i)} $X$ localement
irr\'eductible en $x$ si et seulement si ${\cal O}^{c}_{X,x}={\cal
O}^{0}_{X,x}$\smallskip {\bf (ii)} $X$ faiblement normal en $x$ si
et seulement si ${\cal O}_{X,x}= {\cal O}^{c}_{X,x}$\smallskip {\bf
(iii)}  $X$ normal en $x$ si et seulement si ${\cal O}_{X,x}={\cal
O}^{0}_{X,x}$\smallskip\noindent {\bf 1.2.1. Morphisme
\'equidimensionnel- morphisme ouvert d'espaces analytiques
complexes.}\smallskip\noindent {\bf 1.2.1.0. Formule des dimensions
dans le cadre des alg\`ebres analytiques
locales.}\smallskip\noindent Rappelons que la cat\'egorie des germes
d'espaces analytiques complexes est
 duale de la cat\'egorie des alg\`ebres analytiques (i.e des anneaux quotients
  non nuls d'anneau de s\'eries convergentes); la dualit\'e \'etant concr\'etis\'ee par la
   correspondance
$\ds{(X,x)\rightarrow {\cal O}_{X,x}}$. Il est bien connu que ceux
sont des anneaux locaux
  noeth\'eriens \`a corps r\'esiduel ${\Bbb C}$ et que tout morphisme d'alg\`ebres analytiques est
  local. On dispose du lemme de normalisation de Noether et du th\'eor\`eme
  de pr\'eparation de Weierstrass disant, respectivement, que\par
  $\bullet$  pour
  toute alg\`ebre analytique ${\cal A}$ de dimension $n$, il existe un morphisme fini et injectif
  $\ds{{\Bbb C}\{X_{1},\cdots, X_{n}\}\rightarrow {\cal A}}$ \par
  $\bullet$ un morphisme d'alg\`ebres analytiques est fini si et seulement si il est quasi-fini.\par\noindent
   La dualit\'e \'etant concr\'etis\'ee par la correspondance
$$(X,x)\rightarrow {\cal O}_{X,x}$$
\par\noindent
Nous renvoyons le lecteur aux excellents expos\'es de Cartan,
Grothendieck et Houzel du seminaire Henri Cartan sur les techniques
de constructions en g\'eom\'etrie analytique ([C]). \par\noindent
Dans les constructions qui nous occupent, on s'apper\c coit que les
fondations de cette architecture repose fondamentalement sur {\it{la
notion de morphisme trace}} et de la possibilit\'e d'associer \`a un
morphisme d'alg\`ebres analytiques un tel objet. Une condition
n\'ecessaire est que  le dit morphisme satisfasse au moins {\it{la
formule des dimensions}} c'est-\`a-dire \'etant donn\'e un morphisme
d'alg\`ebres analytiques $f:{\cal A}\rightarrow {\cal B}$, ${\goth
m}$ le radical de ${\cal A}$, alors
$${\rm Dim}{\cal B} = {\rm Dim}{\cal A} + {\rm Dim}{{\cal B}/{\goth m}{\cal B}}$$
Dans ce cas, il existe une suite d'\'el\'ements
$b_{1}$,$\cdots$,$b_{n}$ de ${\cal B}$  dont les classes
r\'esiduelles forment un syst\`eme de param\`etres sur  ${{\cal
B}/{\goth m}{\cal B}}$ (de dimension de Krull $n$ !) et donnant un
morphisme de ${\cal A}$-alg\`ebres $\phi:{\cal A}\{x_{1},\cdots,
x_{n}\}\rightarrow {\cal B}$ envoyant naturellement $x_{j}$ sur
$b_{j}$ faisant de ${\cal B}$ une ${\cal A}$-alg\`ebre finie sur
${\cal A}\{x_{1},\cdots, x_{n}\}$. L'\'enonc\'e g\'eom\'etrique se
formule de la mani\`ere suivante ([Ki] p3, Hilfssatz
2.3):\par\noindent Soient  $f:X\rightarrow S$ est un morphisme
d'espaces analytiques complexes, $x$ un point de $X$ d'image
$s:=f(x)$. Supposons  que $f$ satisfasse {\it{la formule des
dimensions}}
$${\rm Dim}(X,x) = {\rm Dim}(S,s) + {\rm Dim}(f^{-1}(s),x)$$
Alors, il existe un voisinage ouvert $U$  de $x$ dans $X$, un
voisinage ouvert $V$ de $s$ contenant $f(U)$, un voisinage ouvert
$W$ de l'origine de ${\Bbb C}^{n}$ ($n$ \'etant la dimension de la
fibre $f^{-1}(s)$) et un morphisme fini $g:U\rightarrow W\times V$
tel que $f$ se factorise sous la forme
$\xymatrix{X\ar@/_/[rr]_{f}\ar[r]^{g}&W\times U\ar[r]^{q}&S}$, $q$
\'etant  la projection canonique. Rappelons que toutes les
alg\`ebres analytiques locales peuvent \^etre munies d'une topologie
canonique (cf [Ju]).\par\noindent \smallskip\par\noindent {\bf
1.2.1.1.
 Equidimensionnalit\'e et ouverture d'un morphisme.}\par\noindent
Les d\'efinitions suivantes sont unanimement adopt\'ees dans la litt\'erature (cf [E.G.A].IV, [Fi] ou [A.S]). Un morphisme d'espaces complexes
$\pi:X\rightarrow S$  est dit {\it{\'equidimensionnel}}  si
\par\noindent
 {\bf (i)} la fonction semi-continue sup\'erieurement $x\rightarrow {\rm dim}_{x}(\pi^{-1}(\pi(x)))$ est
localement  constante et \par\noindent {\bf (ii)} chaque composante
irr\'eductible de $X$ domine une composante irr\'eductible de
$S$.\smallskip\par\noindent Signalons au passage quel le
th\'eor\`eme de Chevalley sur la
 semi-continuit\'e de la dimension des fibres est encore valable dans le
  cadre de la g\'eom\'etrie analytique.\smallskip\noindent
  Il est important de signaler que le d\'efaut majeur de cette notion est le manque de stabilit\'e par changement de base; ce qui r\'eduit sensiblement ses applications dans le cadre relatif. On dispose, toutefois,  d'une notion plus adapt\'ee au cadre relatif  qui est celle de morphisme {\it{
ouvert}} dont la propri\'et\'e caract\'eristique est de transformer
tout ouvert de $X$ en ouvert de $S$. Cette notion qui correspond \`a
celle de {\it{ universellement ouvert}} de la  g\'eom\'etrie
alg\'ebrique, est, quant \`a elle, stable par changement de
base. Signalons, \`a titre indicatif, que Grothendieck donne
 un exemple simple de morphisme ouvert non universellement ouvert
 dans [E.G.A] IV.{\it 3,Remarque 14.3.9.i}.\smallskip\par\noindent {\bf 1.2.1.2.
Quelques remarques d'ordre g\'en\'eral.}\smallskip\noindent {\bf 1.2.1.2.1. Sur l '\'equidimensionnalit\'e.}\smallskip\noindent {\bf (i)} C'est une
condition ouverte sur la source du morphisme. L'ensemble des points
en lesquels un morphisme est \'equidimensionnel est toujours un
ouvert non vide.\par\noindent {\bf (ii)} Elle n'est pas stable par
changement de base comme le montre l'exemple suivant donn\'e par
\par\noindent
  $\ds{X_{1}:=\lbrace{(x,y,z)\in {\Bbb C}^{3} :
x=z=0}\rbrace}$, $\ds{X_{2}:=\lbrace{(x,y,z)\in {\Bbb C}^{3} : y=0;
z=1}\rbrace}$,\par\noindent $\ds{S_{1}:= \lbrace{(x,y)\in {\Bbb
C}^{2} : x=0}\rbrace}$ et  $\ds{S_{2}:= \lbrace{(x,y)\in {\Bbb
C}^{2} : y=0}\rbrace}$. En posant  $X:=X_{1}\cup X_{2}$ (resp.
$S:=S_{1}\cup S_{2}$) et $\pi:X\rightarrow S$ le morphisme fini et
surjectif induit par la projection canonique $(x,y,z)\rightarrow
(x,y)$, on voit que   $X$ et $S$ sont de dimension pure $1$ et $\pi$
\'equidimensionnel. Mais le changement de
 base $S_{2}\rightarrow S$  donne le sous espace $\pi^{-1}(S_{2})\simeq \{0\}\bigcup {\Bbb
C}$ qui n'est manifestement pas de dimension pure et  de plus la
restriction de
 $\pi$ \`a $\ds{S_{2}\times_{S} X}$ n'est pas \'equidimensionnelle
 sur
 $S_{2}$.\par\noindent
{\bf (iii)}   Un morphisme dont les fibres sont de dimension pure
constante (m\^eme sur une base localement irr\'eductible!), n'est pas n\'ecessairement \'equidimensionnel. Il suffit
pour cela de consid\'erer l'exemple classique de la r\'eunion d'un
plan et d'une droite donn\'e par\par\noindent
 $X:=X_{1}\bigcup X_{2}:=\lbrace{(x,y,z)\in {\Bbb C}^{3} :
z=0}\rbrace\bigcup\rm\lbrace{(x,y,z)\in {\Bbb C}^{3} : y=0;
x=z}\rbrace$ s'envoyant sur  $S:={\Bbb C}^{2}$ gr\^ace au morphisme
fini et surjectif induit par la projection canonique $p(x,y,z)=
(x,y)$.\par\noindent Si la base d'un morphisme \'equidimensionnel
est de dimension pure alors la source aussi.\smallskip\noindent
 {{\bf 1.2.1.2.2. Sur  l'ouverture.}}\smallskip\noindent
 {\bf(iv)} ce n'est pas une condition ouverte
 (tout comme l'irr\'eductibilit\'e locale). Dans l'exemple
 ({\bf{b.iii}}) le morphisme de $X$ sur $S$ est bien ouvert en $0$
 mais sur aucun voisinage ouvert de ce point! Elle est de nature locale sur la base et non
sur la source. l\`a encore, nous voyons que,  m\^eme sur une base localement irr\'eductible, un morphisme \`a fibres de dimension pure constante n'est pas automatiquement ouvert. Il en sera ainsi si la source est de dimension pure.  L'ensemble des points en lesquels un morphisme
d'espaces complexes est ouvert n'est, en g\'en\'eral, ni ouvert ni
ferm\'e mais son compl\'ementaire est  n\'eanmoins constructible (cf
[Par]) \par\noindent {\bf(v)} Elle est stable par changement de base
(et correspond \`a la notion de "universellement ouvert" de la
g\'eom\'etrie alg\'ebrique).\par\noindent
\par\noindent{\bf(vi)} La restriction d'un
morphisme ouvert sur  les composantes irr\'eductibles de la source
n'est pas n\'ecessairement ouverte comme le montre l'exemple suivant
([Par]) : Soient $X:=\lbrace{(x,y,z,t)\in {\Bbb C}^{4} :
x=y=0}\rbrace\bigcup \lbrace{(x,y,z,t)\in {\Bbb C}^{4} :
z=t=0}\rbrace = X_{1}\bigcup X_{2}$ et $\pi:X\rightarrow {\Bbb
C}^{2}$ donn\'e par\par\noindent \centerline{$\pi|_{X_{1}}(x,y,z,t)
= ((z+t)z, (z+t)t)$}\par\noindent \centerline{$\pi|_{X_{2}}(x,y,z,t)
= ((x-y)x, (x-y)y)$} Alors, il est facile de voir que
$\pi|_{X_{1}}(z,t)$ et  $\pi|_{X_{2}}(x,y)$ ne sont pas  ouvertes
\`a l'origine alors que  $\pi$ l'est.\par\noindent
 {\bf(vi)} L'image r\'eciproque par un morphisme ouvert  d'un ouvert dense est toujours dense. L'image par un morphisme ouvert  de toute composante irr\'eductible est un ouvert de la base; en fait cette image est toujours contenue dans au moins une composante irr\'eductible de la base.\par\noindent
 {\bf(vii)} Soit
$\pi:X\rightarrow S$ est un morphisme d'espaces complexes.
Alors,\par $\bullet$ (Remmert [Re]) Si $X$ et $S$
sont de dimension pure avec $S$ localement irr\'eductible alors on a
les \'equivalences \par {\bf (a)} $\pi$ ouvert\par {\bf (b)} ${\rm
Dim}_{x}\pi^{-1}\pi(x)= {\rm Dim}_{x}X - {\rm
Dim}_{\pi(x)}S=m-r,\,\,\,\forall\,x\in X$\par {\bf (c)}${\rm
Dim}\pi^{-1}(s)=m-r,\,\,\forall\,s\in \pi(X)$\par
 $\bullet$ $X$, $S$ de dimension pure et $\pi$ ouvert
impliquent  que les fibres sont de dimension constante pure. Les morphismes de
normalisation faible ou forte montrent que la r\'eciproque est
\'evidemment fausse.\par $\bullet$ si les fibres et $S$ sont de
dimension pures et $\pi$ ouvert surjectif alors $X$  est de
dimension pure. L'exemple {\bf(b.ii)} montre aussi que la
r\'eciproque est fausse puisque $S_{2}$ et les fibres de $\pi$ sont
de dimension pure mais $\pi$ n'est pas ouvert car sinon
$\pi^{-1}(S_{2})$ serait de dimension pure aussi! On peut aussi
avancer le fait que le sous ensemble ouvert $X_{1}$ s'envoie sur
$S_{1}$ qui n'est pas ouvert !\par
\par $\bullet$ si $\pi$ est ouvert et \`a fibres de dimension pure constante alors $S$ de dimension pure est \'equivalent \`a $X$ de dimension pure.\par
$\bullet$ si $\pi$ est \'equidimensionnel et $S$ localement irr\'eductible alors $\pi$ est ouvert.
\smallskip\par\noindent
{\bf 1.2.1.3. Morphisme universellement
\'equidimensionnel.}\smallskip\noindent {\bf  D\'efinition}:  Soient
$n$ un entier naturel et $\pi:X\rightarrow S$ un morphisme d'espaces
complexes (avec $S$ \'eventuellement r\'eduit). Nous dirons que
$\pi$ est {$n$-\it{universellement \'equidimensionnel}} s'il est
ouvert  et \`a fibres de dimension pure
constante.\rm\smallskip\noindent Comme les propri\'et\'es d'\^etre
ouvert et \`a fibres de dimension constante sont stables par
changement de base arbitraire, il en va de m\^eme de la notion de
{\it universellement \'equidimensionnel}.\par\smallskip\noindent
{\bf 1.2.1.4. Remarques.}\par\noindent
 A ce stade, attirons l'attention du lecteur sur la d\'efinition non
standard  de l'\'equidimensionnalit\'e propos\'ee dans [B.M]. En
effet, cette d\'efinition n'exige pas le point {\bf (ii)} de la
d\'efinition d'usage (cf {\bf (1.2.1.1)}) et tol\`ere implicitement les
fibres vides. Mais alors, il ne sera pas du tout vrai que la
puret\'e de la dimension de $S$ entraine celle de $X$ comme nous le montre
 le premier exemple de {\bf (1.2.1.2.1)} ou tout autre
exemple avec $S={\Bbb C}^{m}$ et $X$ un espace complexe quelconque
de dimension non pure. La d\'efinition usuelle imposerait \`a $X$
d'\^etre de dimension pure si $S$ l'est ! (si $X$ et $S$ sont de dimension pure et $\pi$ \`a fibres de
dimension pure constante, $\pi$ n'est pas n\'ecessairement ouvert
comme le montre le morphisme de normalisation; on peut aussi avoir la puret\'e des
dimensions de $X$ et $S$ et un morphisme \`a fibres de dimension
pure mais non constante comme dans le cas des \'eclatements).\par\noindent
 Par ailleurs, elle nous force \`a introduire les cycles relatifs de dimension mixte (si toutefois ils sont bien d\'efinis!); ce que l'on se garde d'envisager ici.\smallskip\noindent
Cons\'equence de cela, il nous parait pr\'ef\'erable d'appeler {\it
pseudo-\'equidimensionnel} ce qui est appel\'e {\it
\'equidimensionnel} dans [B.M].\smallskip\par\noindent {\bf 1.2.2.
Cycle relatif. }\smallskip\noindent {\bf 1.2.2.1. Cycles associ\'es
aux fibres d'un morphisme.}
\par\noindent Si  $X$ est localement de dimension pure et
$\pi:X\rightarrow S$ est {\it{g\'en\'eriquement ouvert}}, les
proc\'edures d\'ecrites dans [Tu] ou [B.M], permettent de munir
(g\'en\'eriquement sur $S$ suppos\'e r\'eduit) de multiplicit\'es  les fibres de $\pi$
et donc de  d\'efinir  une application naturelle\par\noindent
$\ds{\Psi^{\pi}_{gen}: S_{gen}\rightarrow {\cal
C}_{*}(X)\setminus\{0\}}$ donn\'ee par
 $\ds{s\rightarrow [\pi^{-1}(s)]:=\sum_{i}\mu_{i,s}.[X_{i,s}]}$ o\`u $X_{i,s}$ d\'esignent
  les composantes irr\'eductibles et les entiers $\mu_{i,s}$ les multiplicit\'es associ\'ees.
   Le {\it lemme 3.1}  de [Si] nous dit, d'ailleurs, que c'est  une
   injection continue.\smallskip\par\noindent
{\bf 1.2.2.2. D\'efinition}:  Si $X$ et $S$ sont deux espaces
analytiques complexes r\'eduits, on appellera  {\it{cycle $S$-
relatif}} toute combinaison lin\'eaire localement finie $\ds{{\cal
X}:=\sum_{j}n_{j}X_{j}}$ o\`u les $n_{j}$ sont des entiers relatifs
et les $X_{j}$  des sous espaces analytiques irr\'eductibles
r\'eduits  de $S\times X$ tel que le morphisme $\ds{F:|{\cal
X}|:=\bigcup_{j}X_{j}\rightarrow S}$, induit par la projection
$S\times X\rightarrow S$, soit universellement \'equidimensionnelle.
Il est dit {\it{effectif}} si $n_{j}\geq 0$ pour tout indice $j$ et
de {\it dimension pure}  $m$  si les $X_{j}$ sont tous  de dimension
pure $m$.\rm \smallskip\par\noindent {\bf 1.3. Platitude
g\'eom\'etrique.}\smallskip\noindent {\bf 1.3.1.} Comme nous l'avons
d\'ej\`a dit  ${\cal C}_{*}(X)$ est un espace  topologique de
Hausdorff \`a topologie d\'enombrable.  Mais n'\'etant jamais
localement compact, en g\'en\'eral, il ne peut \^etre muni d'une
structure analytique complexe. Cette lacune ne nous permet pas de
parler d'application analytique globale \`a valeurs dans  ${\cal
C}_{*}(X)$ et   nous force \`a travailler avec  des familles
analytiques locales de cycles (cf [B1]). \par\noindent Une famille
de cycles $(X_{s})_{s\in S}$ d'un espace complexe $X$ param\'etr\'ee
par un espace complexe r\'eduit $S$  est dite {\it continue}  s'il
lui est associ\'e une application continue $\Psi:S\rightarrow {\cal
C}_{*}(X)$ d\'efinie comme  dans {\bf(1.2.2.1.)}.  La famille est
dite {\it analytique locale } si elle est continue et donc
associ\'ee \`a une application continue $\Psi$ qui est de plus
"localement analytique" dans le sens suivant:
\par\noindent pour tout $s_{0}$ de $S$, il existe un voisinage
ouvert $S_{0}$ de $s_{0}$ v\'erifiant que, pour toute \'ecaille
$E:=(V,\sigma, U,B)$ $S_{0}$-adapt\'ee (muni de la projection
canonique $p_{1}:U\times B\rightarrow B$), l'application
``classifiante $F:S_{0}\times U\rightarrow {\rm Sym}^{k}(B)$,
associant \`a $(s,t)$ le $0$-cycle $ {p_{1}}_{*}([\{t\}\times B]\cap
\sigma_{*}(\Psi(s)))$, est analytique.\par\noindent Cette donn\'ee
correspond bien \`a
 une fonction analytique locale sur $S$ \`a valeurs dans un espace
banachique. De fa\c con pr\'ecise, si ${\goth{E}}:=((U_{\alpha},
B_{\alpha}, \sigma_{\alpha}))_{\alpha \in A}$ est une famille
d'\'ecaille adapt\'ee au cycle $X_{s_{0}}$, on obtient un voisinage
ouvert $S_{\goth{E}}$ de $s_{0}$ et une application analytique
$$S_{\goth{E}}\rightarrow \prod_{\alpha}{\rm H}(\overline{U}_{\alpha}, {\rm Sym}^{k_{\alpha}}(B_{\alpha}))$$
Dans le cas des cycles non compacts, le comportement parfois
anarchique des cycles ne permet pas de se
 suffire d'une \'etude locale pour en d\'eduire l'aspect global. Nous
  renvoyons \`a ce propos \`a [Ma] pour se rendre compte des difficult\'es
   inh\'erentes \`a cette situation.\smallskip\par\noindent
{\bf 1.3.2. La notion de platitude g\'eom\'etrique continue et
analytique.}\smallskip\noindent Rappelons que cette notion r\'epond
\`a la question de savoir si les fibres d'un morphisme d'espaces
complexes $\pi:X\rightarrow S$  peuvent \^etre munies de
multiplicit\'es convenables de sorte \`a ce qu'elles soient induites
par une famille analytique locale ou seulement continue de cycles.
Le point {\bf(ii)} de la remarque {\bf (1.3.3)} qui va suivre montre
que, par \'egard au cas propre se raccrochant \`a la situation
universelle, les morphismes candidats se doivent d'\^etre
n\'ecessairement universellement ouverts. Cette condition n'est pas
du tout suffisante comme le montre des exemples simples de
normalisation faible (cf {\bf(1.3.6)}). En r\`egle g\'en\'eral, un
morphisme avec \'eclatement (m\^eme "cach\'es") est aux antipodes de
ce que l'on attend.\par\noindent En s'inspirant de [B.M], nous
allons commencer par donner la :\par\noindent
 {\bf 1.3.2.1.  D\'efinition } Soient
$n$ un entier naturel et  $\pi: X\rightarrow S$ un  morphisme
d'espaces analytiques complexes r\'eduits universellement
 $n$-\'equidimensionnel. Soit
${\goth X}$ un cycle $S$-relatif effectif  de $S\times X$.\par\noindent
 Nous dirons que ${\goth X}$ est une {\bf pond\'eration } pour $\pi$ si  \par \centerline{$
(\star)\,\,\,|{\goth X}|\cap{(\{s\}\times X)}= \{s\}\times
\pi^{-1}(s),\,\,\,\forall\,s\in S$} \noindent Elle sera dite
\par
$\bullet$ {\it  contin\^ument g\'eom\'etriquement plate} si ${\goth
X}$ induit une application continue $\Psi_{\pi}:S\rightarrow {\cal
C}_{n}(X)\setminus \{0\}$\par
 $\bullet$  {\it analytiquement  g\'eom\'etriquement plate} si
  ${\goth X}$ correspond au graphe
d'une famille analytique locale
 de $n$- cycles de $X$ ou que l'application pr\'ec\'edente
 $\Psi_{\pi}$ est localement analytique dans le sens de {\bf(1.3.1)}.
 \smallskip\noindent On dira souvent, par abus de langage, que $\pi$
est contin\^ument (resp. analytiquement) g\'eom\'etriquement plat
s'il est muni d'une pond\'eration ayant cette propri\'et\'e.
\smallskip\noindent {\bf 1.3.2.2. D\'efinition} Soient
$\pi:X\rightarrow S$ un morphisme universellement
$n$-\'equidimensionnel d'espaces complexes r\'eduits et $\Sigma$ un
sous espace  de $S$ tel que  le morphisme d\'eduit de $\pi$ dans le
changement de base $S_{gen}
 :=S\setminus \Sigma\rightarrow S$\footnote{$^{(11)}$}{M\^eme
si $S\setminus \Sigma$ n'est pas n\'ecessairement ouvert, en
g\'en\'eral, $\Psi^{\pi}_{gen}$ est n\'eanmoins bien d\'efini
puisque les multiplicit\'es se d\'eterminent localement sur $X$.}
soit muni d'une pond\'eration ${\goth X}_{gen}$ contin\^ument
(resp.analytiquement) g\'eom\'etriquement plate. Alors $\pi$  est
dit  {\it fortement contin\^ument (resp.analytiquement)
g\'eom\'etriquement plat} s'il existe une pond\'eration $\goth{X}$
contin\^ument (resp.analytiquement) g\'eom\'etriquement plate  de
$\pi$  telle que $\goth{X}|_{S_{gen}\times X}= {\goth
X}_{gen}$.\smallskip\noindent Cela revient \`a dire que
l'application  $\Psi^{\pi}_{gen}: S_{gen}\rightarrow {\cal
C}_{*}(X)$ donn\'ee par
 $s\rightarrow [\pi^{-1}(s)]:=\sum_{i}\mu_{i,s}.[X_{i,s}]$ se prolonge en une application continue (resp. analytique locale)  $\Psi^{\pi}: S\rightarrow {\cal
C}_{*}(X)\setminus\{0\}$.\smallskip\noindent {\bf 1.3.3. Remarques.}
\par\noindent {\bf(i)} La notion de pond\'eration effective est de nature locale sur $X$ et $S$.\par\noindent
{\bf (ii)} Si la base est faiblement normal,  fortement  contin\^ument g\'eom\'etriquement plat
correspond \`a g\'eom\'etriquement plat de Siebert[Si] et
contin\^ument g\'eom\'etriquement plat \`a faiblement
g\'eom\'etriquement plat; fortement analytiquement
g\'eom\'etriquement plat correspond \`a g\'eom\'etriquement plat de
Barlet [B2]. \par\noindent {\bf(iii)} Un morphisme $\pi:X\rightarrow
S$ universellement \'equidimensionnel  d'espaces analytiques
complexes r\'eduits admet toujours  une pond\'eration naturelle
donn\'ee par le graphe de $\pi$. En effet, si $G_{\pi}$ est ce
graphe (qui est un sous espace analytique complexe de $S\times X$)
et $\ds{\bigcup_{j}G^{j}_{\pi}}$ sa d\'ecomposition en composantes
irr\'eductibles, le cycle associ\'e
$\ds{[G_{\pi}]=\sum_{j}[G^{j}_{\pi}]}$ est une pond\'eration au sens
de la d\'efinition pr\'ec\'edente. Si $S$ est normal, cette
pond\'eration naturelle est en fait analytiquement
g\'eom\'etriquement plate. En g\'en\'eral, elle n'est m\^eme pas
contin\^ument g\'eom\'etriquement plate. D'ailleurs  des exemples
simples montrent que la condition  $(\star)$ ne garantit pas du tout
l'existence d'une application continue $\Psi_{\pi}$.
\par\noindent {\bf(iv)} Dans le cas d'un d'un morphisme propre, la
platitude g\'eom\'etrique s'exprime naturellement en terme d'espace
de Barlet ${\cal B}(X)$. En effet, au vu de ce qui a \'et\'e dit
pr\'ec\'edemment dans {\bf(1.0.2.2)}, cela revient exactement
 \`a montrer qu'il existe une application  analytique $\Psi^{\pi} :S\rightarrow {\rm B}(X)$ sur $S$
  tout entier avec $|\Psi^{\pi}(s)|=\pi^{-1}(s)$ pour
tout $s\in S$.\par\noindent
 Notons ${\cal  B}(X)\# X$ le sous espace d'incidence,  $\pi_{\#}:{\cal
B}(X)\# X\rightarrow {\cal B}(X)$ le morphisme universel et ${\cal
B}_{\#}$ le cycle relatif universel associ\'e \`a la famille
universelle. Alors, $\pi_{\#}$ est, par d\'efinition, analytiquement
g\'eom\'etriquement plat et par suite ouvert et \`a fibres de
dimension pure constante. Par r\'ef\'erence \`a cette situation
universelle, $\pi$ sera analytiquement g\'eom\'etriquement plat si
et seulement si il  s'ins\`ere dans le diagramme commutatif
$$\xymatrix{X\ar[r]^{\psi_{\pi}}\ar[d]_{\pi}&{\cal B}_{\#}(X)\ar[d]^{\pi_{\cal B}}\\
S\ar[r]_{\Psi_{\pi}}&{\cal B}(X)}$$ et auquel cas, il sera  alors
pond\'er\'e par le cycle relatif ${\cal X}_{\pi}:=\Theta^{*}({\cal
B}_{\#})$.\par\noindent Les propri\'et\'es de constance des fibres
 et d'ouverture d'un morphisme \'etant stables par changement de
base, tout morphisme $\pi$ rendant commutatif le diagramme
pr\'ec\'edent  est n\'ecessairement ouvert et \`a fibres de
dimension pure constante c'est-\`a-dire universellement
$n$-\'equidimensionnel.
\smallskip\noindent Cela nous montre  qu'une condition n\'ecessaire
pour qu'un morphisme d'espaces complexes soit analytiquement
g\'eom\'etriquement plat est qu'il soit  universellement
\'equidimensionnel et surjectif; ce qui \'elimine la plupart  des
morphismes avec \'eclatements ! Encore une fois, des exemples
simples montreront
  au lecteur que cette condition n'est pas suffisante et dans certain cas
  ne
   garantit m\^eme pas la continuit\'e (en tant que cycles) de la famille
   d\'efinie par les fibres.\smallskip\noindent
{\bf 1.3.4. Notations.}\par\noindent
Nous d\'esignerons par :\par\noindent
$\bullet$ ${\cal E}(S,n)$ l'ensemble de tous
    les morphismes d'espaces analytiques complexes $\pi:X\rightarrow
    S$ universellement $n$-\'equidimensionnels,
\par\noindent
    $\bullet$ ${\cal E}_{pond}(S,n)$ l'ensemble des morphismes
    universellement $n$-\'equidimensionnels munis d'une certaine
    pond\'eration $\goth{X}$,\par\noindent
    $\bullet$ ${\cal G}_{c}(S,n)$ (resp. ${\cal G}_{a}(S,n)$) le
    sous ensemble de ${\cal E}_{pond}(S,n)$ constitu\'e
    d'\'el\'ements dont la pond\'eration est contin\^ument
    g\'eom\'etriquement plate (resp analytiquement
    g\'eom\'etriquement plate). Dans la suite, nous donnerons des  exemples simples  montrant que  les inclusions \par
    \centerline{${\cal G}_{a}(S,n)\subset {\cal G}_{c}(S,n) \subset
    {\cal E}_{pond}(S,n)\subset {\cal E}(S,n)$}
\noindent sont  g\'en\'eralement strictes.\smallskip\noindent
  {\bf 1.3.5. Les familles continues (resp. analytiques locales)
   de cycles g\'en\'eriquement r\'eduits.}\par\noindent
 A l'origine, cette notion de platitude
g\'eom\'etrique \'etait r\'eserv\'ee aux familles g\'en\'eriquement
r\'eduites (cf [B2]) et \'etaient d\'efinies, pour un  morphisme
propre d'espaces complexes irr\'eductibles $\pi:X\rightarrow S$, de
fibre g\'en\'erique de dimension pure $n$, muni d'un ferm\'e
analytique d'int\'erieur vide dans $S$ que l'on note
$\Sigma$\footnote{$^{(12)}$}{On peut prendre $\Sigma$ comme \'etant
la r\'eunion des points non normaux de $S$ et des sous ensembles
$\lbrace{s\in S : {\rm dim}(\pi^{-1}(s))> n} \rbrace\bigcup
\lbrace{s\in S : \pi^{-1}(s)\,\, {\rm est}\, {\rm multiple}\rbrace}$
} et d'une application holomorphe ${\Psi^{\pi}_{gen}}:S_{gen}:=
{S\setminus\Sigma}\rightarrow {\cal B}_{n}(X)$ donn\'ee par
$\ds{s\rightarrow [\pi^{-1}(s)]:=\sum_{i}[X_{i,s}]}$, $X_{i,s}$
\'etant  les composantes irr\'eductibles de $\pi^{-1}(s)$ (en nombre
fini puisque $\pi$ est propre), en disant que $\pi$ est  {\it
g\'eom\'etriquement plat} (c'est-\`a-dire {\it fortement
analytiquement g\'eom\'etriquement plat} dans notre terminologie) si
et seulement si  ${\Psi^{\pi}_{gen}}$ se prolonge analytiquement sur
$S$ tout entier.\smallskip\noindent On peut sensiblement
 g\'en\'eraliser cette d\'efinition (cf [Va]) au cas d'un
 un morphisme $\pi:X\rightarrow S$
propre et $n$- \'equidimensionnel d'espaces complexes r\'eduits.
Ainsi, $\pi:X\rightarrow S$ \'etait  dit g\'eom\'etriquement plat
(ou  {\it{fortement analytiquement g\'eom\'etriquement plat}} selon
notre d\'efinition) si il existe une famille analytique
$(X_{s})_{s\in S}$ de $n$-cycles de $X$ param\'etr\'ee par $S$ telle
que \par {\bf(i)} $X_{s}= [\pi^{-1}(s)]$,  si $s$ est
r\'egulier\footnote{$^{(13)}$}{Cela signifie que $s$ est lisse dans
$S$ et qu' aucune composante irr\'eductible de $\pi^{-1}(s)$ n'est
contenue dans le lieu singulier de $X$ ou dans le lieu constitu\'e
des points de $X$ en lesquels $\pi$ n'est pas de rang maximum.}\par
{\bf(ii)}  $|X_{s}|= \pi^{-1}(s)$ pour tout $s\in
S$.\smallskip\noindent En dehors des cas particuliers o\`u $S$ est
localement irr\'eductible, faiblement normal ou normal, on ne peut
absolument rien dire sur l'existence d'un tel prolongement m\^eme
par continuit\'e!\par\noindent Cette approche pr\'esente deux
inconv\'enients majeurs \`a savoir ses d\'efauts de fonctorialit\'e
flagrants et le fait qu'elle ne s'interesse qu'aux familles
g\'en\'eriquement r\'eduites puisqu'en en \'evitant des ferm\'es
d'int\'erieur vide aussi ``gros'' on \'elimine des situations
g\'eom\'etriques simples telle que famille de coniques
d\'eg\'en\'erant en une droite double ! \footnote{$^{(14)}$}{ Il est
 facile de se convaincre qu'une condition n\'ecessaire assurant l'analyticit\'e d'une  famille
g\'en\'eriquement r\'eduite $(X_{s})_{s\in S}$ est que, pour tout
$s$, les  composantes irr\'eductibles  $X_{i,s}$ de $\pi^{-1}(s)$ ne
soient pas enti\`erement incluse dans $\Sigma_{3}:=\lbrace{(s,x)\in
S\times X : X\, {\rm non}\, S-{\rm lisse}\,\,{\rm  en}\,
(s,x)}\rbrace$.}
\smallskip\noindent
En supposant  $S$  {\bf  faiblement normale}
\footnote{$^{(15)}$}{Cette condition assure l'existence d'un
prolongement analytique de toute fonction continue et
g\'en\'eriquement holomorphe sur $S$ et permet, ainsi, de se suffire
d'une \'etude locale et un traitement ``presque topologique `` de la
question.} et l'image de $\pi$ dense dans $S$ ( donc pas
n\'ecessairement propre ou surjectif et  en autorisant les
composantes multiples!), Siebert ([Si]) proposait les d\'efinitions
suivantes:\par\noindent soient $\pi:X\rightarrow S$ un morphisme
{\it{g\'en\'eriquement ouvert}}, $E$ le lieu de d\'eg\'enerescence
de $\pi$ (cf [Fi], [Si]), $N$ le lieu non normal de $S$ et
$S_{gen}:= \pi(S\setminus \pi^{-1}(N\cup \pi(X)))$. Alors $\pi$ est
{\it{g\'eom\'etriquement plat}}\footnote{$^{(16)}$}{Si $\pi$ est
propre, on retrouve la notion de platitude forte puisque sur un
espace  faiblement normal toute fonction continue et
g\'en\'eriquement holomorphe est holomorphe.}(c'est-\`a-dire {\it
fortement analytiquement  g\'eom\'etriquement plat} dans notre terminologie) si
l'application associ\'ee $\ds{\Psi^{\pi}_{gen}: S_{gen}\rightarrow
{\cal C}_{*}(X)\setminus\{0\}}$ se prolonge contin\^ument sur $S$ et
il est dit {\it{faiblement g\'eom\'etriquement plat}}
(c'est-\`a-dire {\it contin\^ument g\'eom\'etriquement plat} dans
notre terminologie) si il existe une certaine application continue
$\ds{\Psi_{\pi}: S\rightarrow {\cal
C}_{*}(X)\setminus\{0\}}$.\smallskip\noindent Dans [Si] ( {\it proposition  5.3}, p.258), on retrouve le fait que la platitude g\'eom\'etrique forte
 est associ\'ee aux familles analytiques de cycles g\'en\'eriquement
r\'eduits.
\par\noindent La platitude g\'eom\'etrique dont il sera question
dans cet article correspond \`a la platitude g\'eom\'etrique faible
de Siebert sur une base r\'eduite quelconque \`a laquelle on rajoute
une condition d'analyticit\'e locale. Comme nous allons le voir
cette notion b\'en\'eficie, quant \`a elle, de bonnes propri\'et\'es
fonctorielles en les arguments.\smallskip\noindent {\bf 1.3.6.
Quelques remarques g\'en\'erales sur  la platitude g\'eom\'etrique
.}\par\noindent On a d\'ej\`a signal\'e au lecteur que, pour
esp\'erer greffer une structure de cycles sur les fibres d'un
morphisme d'un espace complexe, l'ouverture est une condition
absolument  n\'ecessaire mais loin d'\^etre suffisante. En
cons\'equence, les morphismes avec \'eclatements sont \`a \'eviter
en g\'en\'eral. Les exemples simples suivants (donn\'es en
g\'en\'eral par les normalisations fortes ou faibles)  vont
permettre d'illustrer les principales obstructions rencontr\'ees.
Rappelons que le morphisme de normalisation faible est fini,
surjectif et ouvert. \par\noindent {\bf(i) Morphisme
\'equidimensionnel surjectif non g\'eom\'etriquement plat.}\par\noindent $\bullet$
Consid\'erons l'exemple de la page 76 de [G.P.R],  de la surface faiblement normale
de ${\Bbb C}^{3}$ commun\'ement appel\'ee  ``parapluie de Cartan''
et donn\'ee par  \par\noindent $S:= \lbrace{(x,y,z)\in {\Bbb C}^{3}
: x^{3}-z(x^{2}-y^{2})=0}\rbrace$. Elle est  irr\'eductible, non
localement irr\'eductible (car sinon elle serait normale!) et admet
pour lieu singulier  la r\'eunion des droites $\Sigma_{1}=\{x=y=0\}$
et $\Sigma_{2}=\{x=z=0\}$ dans ${\Bbb C}^{3}$. Soit $X={\Bbb C}^{2}$
et $\pi:X\rightarrow S$ le morphisme fini et surjectif  qui \`a
$(u,v)$ associe $(v(u^{2}-v^{2}), v(u^{2}-v^{2}), v^{3})$. Alors, il
est facile de constater  que   $\pi^{-1}(s)$ est constitu\'e d'un
point pour $s=0$, de six points pour $s\in \Sigma_{1}\setminus\{0\}$
et de trois points si $s\in \Sigma_{2}\setminus\{0\}$ et que, par
cons\'equent, le degr\'e du rev\^etement ramifi\'e g\'en\'eralis\'e
n'est pas localement constant. Cette pathologie se rencontre fr\'equemment quand le morphisme n'est pas ouvert. \par\noindent $\bullet$ Soit
$X:={\Bbb C}$ et $S:=\lbrace{(x,y)\in {\Bbb C}^{2} : x^{3} - x^{2} +
y^{2} = 0}\rbrace$ et \par \centerline{$\pi:X\rightarrow S$}\par
\centerline{ $t\rightarrow (1-t^{2}, t^{2}-t^{3})$}
\smallskip\noindent qui est un morphisme \'equidimensionnel mais non
ouvert (ce n'est pas la normalisation faible)  induisant une
application g\'en\'eriquement holomorphe\par
 \centerline{$\Psi^{\pi}_{gen}:S-\{0,0\}\rightarrow {\cal B}_{0}(X)$}\par
\centerline{$(x,y)\rightarrow \{{y\over{x}}\}$}\smallskip\noindent
La fonction m\'eromorphe $y\over{x}$ est localement born\'ee
puisqu'elle  v\'erifie une \'equation de d\'ependance int\'egrale
mais ne se prolonge pas contin\^ument sur $S$ qui est
irr\'eductible, faiblement normal mais  non localement
irr\'eductible car sinon il serait normal! Dans ce cas, $\pi$ n'est
m\^eme pas contin\^ument g\'eom\'etriquement plat.\par\noindent
$\bullet$  Dans la m\^eme veine, consid\'erons\par\noindent $X:=
\lbrace{(x,y,z)\in {\Bbb C}^{3} : x= z^{2} - 1;\,y= z^{3}
-z}\rbrace$\par\noindent $S:=\lbrace{(x,y)\in {\Bbb C}^{2} : y^{2}=
x^{2}(x+1)}\rbrace$ et $\pi:X\rightarrow S$ induit par la projection
canonique $p:{\Bbb C}^{3}\rightarrow {\Bbb C}^{2}$ envoyant
$(x,y,z)$ sur $(x,y)$. Alors $X$ et $S$ sont de dimension pure $1$,
$S$ est irr\'eductible (puisque c'est l'image de ${\Bbb C}$ par
l'application $t\rightarrow (t^{2}-1, t^{3}-t)$), $\pi$ est
\'equidimensionnel mais non ouvert puisqu' un voisinage ouvert lisse
de $(0,0,1)$ dans $X$ ne peut \^etre appliqu\'e sur  un voisinage
ouvert de $(0,0)$ dans $S$ dont le germe en ce point est la
r\'eunion de deux germes lisses. L\`a encore, il est facile de voir
que la fonction m\'eromorphe et localement born\'ee
$\ds{{y\over{x}}}$ ne se prolonge pas de fa\c con continue sur $S$
qui est de toute \'evidence non localement irr\'eductible. il en
r\'esulte que l'application g\'en\'eriquement holomorphe\par
 \centerline{$\Psi^{\pi}_{gen}:S-\{0,0\}\rightarrow {\cal B}_{0}(X)$}\par
\centerline{$(x,y)\rightarrow \{x,y,{y\over{x}}\}$} \noindent ne se
prolonge pas par continuit\'e sur $S$.\par\noindent {\bf(ii)
Morphisme \'equidimensionnel surjectif  ouvert et non g\'eom\'etriquement plat.}\par\noindent
En g\'en\'eral, les morphismes de normalisation faible r\'epondent \`a ces crit\`eres. Pour l'illustrer, consid\'erons l'exemple simple donn\'e par: \par\noindent
 $X:={\Bbb C}$, $S:=\lbrace{(x,y)\in {\Bbb C}^{2} : x^{3} -
y^{2} = 0}\rbrace$ et $\pi:X\rightarrow S$ le morphisme de
normalisation faible (ou forte!) donn\'e par $t\rightarrow (t^{2},
t^{3})$. \par\noindent $\pi$ est universellement $0$-ouvert et
induit une application (g\'en\'eriquement holomorphe)
$\Psi^{\pi}_{gen}:S-\{0,0\}\rightarrow {\cal
 B}_{0}(X)$ donn\'ee par $(x,y)\rightarrow
\{{y\over{x}}\}$\smallskip\noindent
 Comme $S$ est  localement
irr\'eductible et que la  fonction m\'eromorphe $y\over{x}$
v\'erifie  une \'equation de d\'ependance int\'egrale, elle  se
prolonge contin\^ument sur $S$. Ainsi $\pi$ est fortement
contin\^ument g\'eom\'etriquement plat puisque  l'application
$\Psi^{\pi}_{gen}$ se prolonge contin\^ument sur $S$. Mais ce
dernier \'etant non faiblement normal, ce  prolongement n'est pas
holomorphe et donc  $\pi$ n'est donc pas fortement analytiquement
g\'eom\'etriquement plat.\par\noindent {\bf (iii) La platitude
g\'eom\'etrique forte n'est pas stable par changement de
base.}\par\noindent C'est le d\'efaut majeur de cette notion que l'on peut comprendre ais\'ement puisque, par changement de base, on ne peut emp\^echer, en g\'en\'eral, l'apparition de branches multiples.  Pour s'en convaincre, il suffit de consid\'erer  l'exemple de
Douady d'une r\'eunion de deux plans de ${\Bbb C}^{4}$ dont la
projection canonique sur ${\Bbb C}^{2}$ est ouverte mais non plate.
A un changement lin\'eaire de coordonn\'ees pr\`es, on retrouve cet
exemple dans [V] exprim\'e sous la forme: \par\noindent Soient
$X:=\lbrace{(z_{1}, z_{2}, z_{3}, z_{4})\in {\Bbb C}^{4} :
z_{1}z_{2}=z_{1}z_{4}=z_{2}z_{3}=z_{3}z_{4}=0}\rbrace$, $S:={\Bbb
C}^{2}$,\par\noindent
 $S':= \lbrace{(x_{1}, x_{2})\in {\Bbb C}^{2} : x_{1}x_{2}
=0}\rbrace$.  Consid\'erons les morphismes  $\pi:X\rightarrow S$ et
$\nu:S'\rightarrow S$  donn\'es, respectivement,  par $\pi(z_{1},
z_{2}, z_{3}, z_{4})=(z_{1} + z_{2}, z_{3} + z_{4})$ et $\nu(x_{1},
x_{2})=(x_{1}, 0)$. Alors dans le diagramme cart\'esien
$$\xymatrix{X\times_{S} S'\ar[r]\ar[d]_{\pi_{1}}&X\ar[d]^{\pi}\\
S'\ar[r]_{\nu}&S}$$ $\pi$ est fortement  g\'eom\'etriquement plat
puisqu'il est \'equidimensionnel sur une base normale mais $\pi_{1}$
ne l'est pas puisque ce produit fibr\'e est constitu\'e d'une droite
triple et de deux droites simples au dessus de chacune des deux
branches de $S'$.\smallskip\noindent Remarquons qu'en pond\'erant la
droite triple  par un coefficient $2$ et chacune des droites simples
d'un coefficient $3$, $\pi_{1}$  devient analytiquement
g\'eom\'etriquement plat.\smallskip\noindent
{\bf 1.3.7.  Exemples standards de morphismes analytiquement g\'eom\'etriquement plats.}\smallskip\noindent
Les cas  fr\'equemment
rencontr\'es sont ceux donn\'es par :\smallskip\noindent
 $\bullet$ {\it La pond\'eration standard}:  Soient  $Z$ et $S$ deux espaces analytiques complexes avec $S$ r\'eduit, $(X_{s})_{s\in S}$  une famille analytique de $n$-cycles de $Z$
 param\'etr\'ee par $S$ dont le support du  graphe $X=\{(s,z) \in S\times Z : z\in |X_{s}|\}$ (muni de sa structure r\'eduite) est muni de la projection  $\pi :
 X\rightarrow S$ induite par la projection canonique de $S\times Z$ sur $S$. la projection canonique. Alors le graphe de
 la famille $(\{s\}\times X_{s})_{s\in S}$, qui est un cycle
 $\goth{X}$ de $S\times Z$, est une pond\'eration analytiquement  g\'eom\'etriquement
 plate de $\pi$.\smallskip\noindent
$\bullet$ {\it La pond\'eration alg\'ebrique}:  Tout morphisme plat d'espaces complexes et  $n$-
\'equidimensionnel,  $\pi : X\rightarrow  S$, est naturellement
pond\'er\'e g\'eom\'etriquement. En effet, si $\pi$ est propre,  le
morphisme $ Douady$ -$ Barlet$ (cf {\bf (1.0.2.2)} ou  [B1], chap.
V, thm 8)  permet de donner un sens \`a cela en prenant  pour
pond\'eration
  g\'eom\'etriquement plate le cycle donn\'e par le graphe de la
  famille analytique $ ([\pi^{-1}(s)])_{s\in S}$ (o\`u $ [\pi^{-1}(s)]$ est
  le cycle associ\'e au sous espace  $\pi^{-1}(s)$). Dans le cas g\'en\'eral, la proc\'edure est la m\^eme puisque assigner des multiplicit\'es est une op\'eration locale sur $X$. A la lumi\`ere du th\'eor\`eme 2, ce qui vient d'\^etre dit est valable pour les morphismes d'espaces complexes  de {\it Tor-dimension} finie avec base r\'eduite.\par\noindent
$\bullet$ {\it La pond\'eration normale}: tout morphisme $\pi:X\rightarrow S$ \'equidimensionnel sur
une base normale $S$ est analytiquement  g\'eom\'etriquement
plat; en fait fortement analytiquement g\'eom\'etriquement plat.\smallskip\noindent $\bullet$ Les deux premiers  exemples  du
paragraphes {\bf (1.1)}  donnent des morphismes analytiquement
g\'eom\'etriquement plats sur une base normale et faiblement normale
respectivement  alors que le dernier est seulement contin\^ument
g\'eom\'etriquement plat sur une base non faiblement
normale.\par\noindent
Pour terminer, citons les  exemples  suivants tir\'es de  [B.M] et [Si] pr\'esentant deux situations diff\'erentes  dans lesquelles il est possible ou pas ``d'ajuster'' convenablement les multiplicit\'es du  cycle relatif donn\'e naturellement par le graphe du morphisme, de sorte \`a rendre cette pond\'eration  analytiquement ou seulement contin\^ument  g\'eom\'etriquement plate.\smallskip\noindent
$\bullet$ Soient $S=\lbrace{(x,y)\in {\Bbb C}^{2}/ xy=0}\rbrace$,
$X_{1}=\lbrace{(x,y,z)\in{\Bbb C}^{3} /
z^{2}=x,\,\,y=0}\rbrace$,\par\noindent
$X_{2}=\lbrace{(x,y,z)\in{\Bbb C}^{3} / z^{3}=y,\,\,x=0}\rbrace$  et
$X=X_{1}\bigcup X_{2}$.\smallskip\noindent Soit $\pi:X\rightarrow S$
induit par la projection canonique $p:{\Bbb C}^{3}\rightarrow {\Bbb
C}^{2}$ donn\'ee par  $(x,y,z)\rightarrow (x,y)$. Alors, $\pi$ est
un morphisme fini surjectif et ouvert dont le graphe admet les deux
composantes irr\'eductibles
$\tilde{X}_{1}=\lbrace{(x,y,a,b,c)\in{\Bbb C}^{3} /
c^{2}=a=x,\,\,y=b=0}\rbrace$ et\par\noindent
$\tilde{X}_{2}=\lbrace{(x,y,a,b,c)\in{\Bbb C}^{3} /
c^{3}=b=y,\,\,x=a=0}\rbrace$.\smallskip\noindent Alors, le  cycle
relatif $Y_{1}:= [\tilde{X}_{1}] + [\tilde{X}_{2}]$ est une
pond\'eration naturelle du morphisme $\pi$. Mais vu les
multiplicit\'es, il ne peut faire de $\pi$ un morphisme
contin\^ument g\'eom\'eriquement plat. Par contre, le cycle relatif
$Y_{2}:= 3[\tilde{X}_{1}] + 2 [\tilde{X}_{2}]$ est une pond\'eration
contin\^ument g\'eom\'eriquement plate et donc analytiquement
g\'eom\'etriquement plate puisque la base est faiblement
normale.\smallskip\noindent $\bullet$ Soient
$S_{1}:=\lbrace{(x,y,t)\in {\Bbb C}^{3} : t=0}\rbrace$,
$S_{2}:=\lbrace{(x,y,t)\in {\Bbb C}^{3} : y
=0}\rbrace$,\par\noindent
 $X_{1}:=\lbrace{(x,y,z,t)\in {\Bbb
C}^{4} : x^{2}+y^{2}-z^{2}= t=0}\rbrace$ et
$X_{2}:=\lbrace{(x,y,z,t)\in {\Bbb C}^{4} : x-z = y=0}\rbrace$.
Posons $X:= X_{1}\bigcup X_{2}$, $S:=S_{1}\bigcup S_{2}$ et $\pi:
X\rightarrow S$ le morphisme ouvert, fini et surjectif  induit par
la projection lin\'eaire $p:{\Bbb C}^{4}\rightarrow{\Bbb C}^{3}$
envoyant $(x,y,z,t)$ sur $(x,y,t)$. Alors, pour $i=1, 2$, on a
$\pi(X_{i})=S_{i}$ et  les multiplicit\'es  des restrictions de
$\pi$ \`a $X_{1}$ et $X_{2}$ sur $X_{1}\cap X_{2}$ au dessus de
$(0,0)$ coincident et valent  $1$ (en dehors de ce point, elles
valent $2$ et $1$ respectivement). La platitude g\'eom\'etrique
imposerait l'\'egalit\'e des poids associ\'es aux diff\'erentes
composantes irr\'eductibles ce qui ne peut avoir lieu ici puiqu'il
est  impossible de choisir un couple d'entiers non nuls  $(a,b)$
pond\'erant chacune des composantes irr\'eductibles de sorte \`a
r\'ecup\'erer, par continuit\'e, la valeur $1$
!\smallskip\par\noindent {\bf 1.3.8. Propri\'et\'es fonctorielles de
la notion de platitude g\'eom\'etrique analytique ou
continue.}\smallskip\noindent {\bf 1.3.8.1. Fonctorialit\'e
basique.}\par\noindent
Tout comme la platitude, cette notion est de
nature locale sur la source; ce qui ram\`ene bon nombre de
probl\`emes sur les morphismes g\'eom\'etriquement plats \`a des
probl\`emes d'alg\`ebres analytiques locales.\par\noindent Dans la
suite, nous \'ecrirons ``g\'eom\'etriquement plat'' quand
l'\'enonc\'e propos\'e est valable aussi bien pour {\it
contin\^ument g\'eom\'etriquement plat} que pour {\it analytiquement
g\'eom\'etriquement plat}.\smallskip\noindent
 {\bf (1) Compatibilit\'e avec les
changements de bases sur $S$:} :\par\noindent Soient
$\pi:X\rightarrow S$ un morphisme $n$-g\'eom\'etriquement plat et $\eta:T\rightarrow S$ un morphisme
d'espaces complexes r\'eduits de dimension localement pure. Le
morphisme $\hat{\pi}:{T\times_{S} X}\rightarrow T$ d\'eduit du
changement de base $\eta$ est naturellement $n$- universellement
\'equidimensionnel (rappelons \`a ce sujet, que de fa\c con
g\'en\'erale, un morphisme est ouvert si et seulement si sa
restriction au r\'eduit est ouverte). La platitude g\'eom\'etrique
se lit donc sur $\hat{\pi}_{red}:(T\times_{S} X)_{red}\rightarrow
S$. Soit ${\cal X}$ le cycle relatif d\'efinit par $\pi$. Comme la
notion de famille analytique (resp. continue) de cycles est stable par changement de
base, $\hat{\pi}$ est aussi g\'eom\'etriquement plat de
pond\'eration $\Theta^{*}({\cal X}):= \hat{X}$ associ\'ee \`a la
famille $(X_{\eta(t)})_{t\in T}$.\smallskip\noindent
    {\bf (2)  Compatibilit\'e avec les inclusions ouvertes sur
$X$:}\par\noindent Si $U$ est un ouvert de $X$ dont l'image par
$\pi$ est l'ouvert  $S_{U}$ et  $j:U\rightarrow X$ l'injection
naturelle. Alors, de fa\c con naturelle, le morphisme induit
$\pi:U\rightarrow S_{U}$ est aussi g\'eom\'etriquement
plat.\smallskip\noindent
 {\bf (3)  Compatibilit\'e avec les images directes de cycles}\par\noindent
Soit $\xymatrix{ X\ar[rr]^{f}\ar[rd]_{\pi}&&Z\ar[ld]^{\tilde{\pi}}\\
&S&}$  un diagramme commutatif de $S$-morphismes d'espaces complexes
r\'eduits  avec $f$ propre et surjectif.  Soient  $F:=f\times Id:
S\times X\rightarrow S\times Z$ le morphisme propre d\'eduit de $f$
et ${\cal X}$ le cycle relatif pond\'erant $\pi$. Il est facile de
voir que $F_{*}({\cal X})$ est une pond\'eration g\'eom\'etriquement
plate naturelle de $\tilde{\pi}$.\smallskip\par\noindent {\bf
1.3.8.2.Image r\'eciproque de cycles par un morphisme
g\'eom\'etriquement plat et composition.}} \smallskip\noindent
Il est bien connu que si $\pi:X\rightarrow S$ est un morphisme {\it plat} d'espaces complexes, l'image r\'eciproque d 'un cycle $Y$ de $S$ est bien d\'efinie (cf [Fu]) puisque pour tout sous espace $Z$ de $S$ on a $[\pi^{*}(Z)]=\pi^{*}[Z]$. Il se trouve qu'elle garde encore un sens si $\pi$ est analytiquement g\'eom\'etriquement plat. En effet, sa  pond\'eration permet de ``distribuer'' convenablement les multiplicit\'es sur les composantes irr\'eductibles de l'image r\'eciproque ensembliste. D'ailleurs, le  {\it th\'eor\`eme 2} (cf \'enonc\'e p.9) permet de d\'ecrire et d'expliciter, par le biais du {\it th\'eor\`eme 0} (cf p.77), le calcul des multiplicit\'es. Comme ces questions sont plus ou moins li\'ees \`a la th\'eorie de l'intersection que nous n'abordons pas, nous renvoyons le lecteur, pour plus d'informations en la mati\`ere,  au \S{\bf 8}, p.130 de  [Fu]  ou  \`a la page 840 de  [B.M] proposant les r\'esultats connus de Fulton mais exprim\'es dans le langage de l'espace des cycles et en terme de la notion de famille analytique de cycles param\'etr\'ee par un espace complexe {\bf lisse}. \smallskip\noindent
  {\bf
1.3.8.2.0.  Image r\'eciproque d'un cycle.} \par\noindent Soit
$\pi:X\rightarrow S$ un morphisme analytiquement g\'eom\'etriquement
plat de poids $\goth{X}$. Soit
 $Y$ un sous ensemble irr\'eductible de $S$ et
$\ds{\pi^{-1}(Y)=\bigcup_{m}{\cal Y}_{m}}$ la d\'ecomposition en
composantes irr\'eductibles de son image inverse. Soit alors
$\mu_{m}$ le plus petit degr\'e de ramification donn\'ee par les
 factorisations locales de $\pi$ le long de ${\cal Y}_{m}$ (cf {\bf (1.0.1.2)}). On pose, alors,
$$\pi^{*}(Y)=\sum_{m}\mu_{m}{\cal Y}_{m}$$
et on l'appelle {it image r\'eciproque analytiquement
g\'eom\'etriquement plate} de $Y$. Il est absolument clair que cette
 op\'eration s'\'etend par lin\'earit\'e \`a un cycle quelconque
$\ds{Y=\sum_{i\in I}k_{i}Y_{i}}$ (la somme \'etant toujours
suppos\'ee localement finie et $Y_{i}$ irr\'eductible) en
d\'efinissant une application $\pi^{*}:{\cal C}(S)\rightarrow {\cal
C}(X)$.\par\noindent On peut remarquer que cette image r\'eciproque
 s'exprime naturellement en fonction de la pond\'eration. En effet, pour
 chaque entier $i$, le changement de base $Y_{i}\hookrightarrow S$,
 donne un cycle bien d\'efini par la relation
 ${\goth Y}_{i}:=\goth{X}\bullet(Y_{i}\times X)$ correspondant \`a
 la pond\'eration analytiquement g\'eom\'etriquement plate de ${\rm
 R}ed(Y_{i}\times_{S} X)\rightarrow Y_{i}$ et
 repr\'esentant le graphe de la famille analytique (resp. continue) de cycles
 d\'eduite de la famille analytique
 (resp. continue) $([\pi^{-1}(s)])_{s\in S}$. Cela nous permet de
 distribuer convenablement les multiplicit\'es sur chacune des
 composantes irr\'eductibles de $\pi^{-1}(Y_{i})$ et de poser
 $$\pi^{*}(Y):=\sum_{i\in I}k_{i}{\goth Y}_{i}$$
{\bf 1.3.8.2.1.  Image r\'eciproque d'une famille
cycles.}\par\noindent L'image r\'eciproque d'une famille analytique
(resp. continue) de cycles par un morphisme g\'eom\'etriquement plat
se d\'efinit de fa\c con similaire (cf [B.K] ou [Si]). \par\noindent
En effet, si $\pi:X\rightarrow Z$ est un morphisme
g\'eom\'etriquement plat de poids $\goth{X}$ et $({\cal C}_{s})_{s
\in S}$ une famille analytique ou continue de cycles de $Z$
param\'etr\'ee par l'espace complexe r\'eduit $S$. Alors, pour
chaque $s \in S$, on d\'efinit, comme pr\'ec\'edemment, le cycle
${X}_{{\cal C}_{s}}:=\pi^*({\cal C}_t)$  comme le graphe de la
famille analytique ou continue d\'eduite du changement de base
${\cal C}_{s}\hookrightarrow Z$.
\smallskip\noindent
{\bf 1.3.8.2.3. Composition de morphismes g\'eom\'etriquement
plats.}\par\noindent Disons tout de suite que la classe des
morphismes analytiquement g\'eom\'etriquement plats n'est, en
g\'en\'eral, pas stable pour la composition comme le montre le
contre exemple de [B1] cit\'e dans \S{\bf(1.1)}, {\bf(iii)}, p.28.
Par contre, celle des contin\^ument g\'eom\'etriquement plats
l'est.\par\noindent
 Commen\c cons par voir, dans le cas propre, o\`u se nichent les
 \'eventuelles obstructions. Soit $\pi_{1}:X\rightarrow S$ (resp. $\pi_{2}:S\rightarrow T$) un morphisme analytiquement g\'eom\'etriquement plat de dimension relative $n_{1}$ (resp. $n_{2}$).  Si ces morphismes sont propres,  on a des plongements analytiques
$\eta_{1}: S\rightarrow {\rm B}_{n_{1}}(X)$ (resp.  $\eta_{2}:
T\rightarrow {\rm B}_{n_{2}}(S)$) donnant une application analytique
$\eta:T\rightarrow{\rm B}_{n_{2}}(S)\rightarrow  {\rm
B}_{n_{2}}({\rm B}_{n_{1}}(X))$. Alors, montrer que la compos\'ee
$\pi:=\pi_{2}o\pi_{1}$ est analytiquement g\'eom\'etriquement plate
revient exactement \`a montrer que l'application g\'en\'eriquement
holomorphe $\psi:{\rm B}_{n_{2}}({\rm B}_{n_{1}}(X))\rightarrow {\rm
B}_{n_{1}+n_{2}}(X)$ se prolonge analytiquement. Il en r\'esultera
que l'application continue et g\'en\'eriquement holomorphe
$T\rightarrow {\rm B}_{n_{1}+n_{2}}(X)$ est en fait analytique; ce
qui traduit la platitude g\'eom\'etrique analytique de
$\pi$.\smallskip\noindent Si la base est {\bf faiblement normale} et
sans condition de propret\'e,  on a le  r\'esultat positif (cf [Si],
{\it proposition 5.12}, p.262) pour les morphismes fortement
g\'eom\'etriquement plats selon notre terminologie. Il n'est pas
difficile de le donner sous la forme g\'en\'erale \Prop{1}{} Soient
$\pi_{1}:X\rightarrow S$ et $\pi_{2}:S\rightarrow T$ des morphismes
analytiquement  g\'eom\'etriquement plats  d'espaces analytiques
complexes r\'eduits  (d\'enombrables \`a l'infini) avec $T$
faiblement normal. Alors la compos\'ee $\pi:X\rightarrow S$ est un
morphisme analytiquement  g\'eom\'etriquement
plat.\rm\smallskip\noindent \dem Comme le probl\`eme est de nature
locale, consid\'erons un point $x_{0}$ de $X$ dans $S\times U\times
B$ et  une \'ecaille adapt\'ee \`a $S_{t_{0}}$ pr\`es de
$\pi_{1}(x_{0})$ donn\'ee par le $T$-plongement $S \hookrightarrow
T\times U' \times B'$
 et une \'ecaille adapt\'ee \`a $X_{\pi_{1}(x_{0})}$ (autour de $x_{0}$ ) relative au plongement
$X\hookrightarrow S \times U \times B \hookrightarrow T\times U'
\times U \times B' \times B$. \noindent associ\'ees aux applications
analytiques $f_{X/S}: S\times U\rightarrow{\rm Sym}^{k}(B)$ et
$g_{S/T}:T\times U'\rightarrow{\rm Sym}^{k'}(B')$ elles m\^emes
associ\'ees aux morphismes {\it{trace}}: ${\cal T}^{0}_{f}:
{\pi_{1}}_{*}({\cal O}_{X})\rightarrow{\cal O}_{S\times U}$ et
${\cal T}^{0}_{g}:{\pi_{2}}_{*}({\cal O}_{S})\rightarrow{\cal
O}_{T\times U'}$.\par\noindent
 $f_{X/S}$ et $g_{S/T}$ permettent de construire un morphisme
 $h:T\times U\times U'\rightarrow {\rm Sym}^{kk'}(B\times B')$
 obtenu en composant les applications holomorphes
suivantes:\smallskip\noindent \centerline{$g_{S/T}\times Id_{U}: T\times
U'\times U\rightarrow {\rm Sym}^{k'}(S\times U)$}\smallskip\noindent
\centerline{${\rm Sym}^{k'}(f_{X/S}):{\rm Sym}^{k'}(S\times
U)\rightarrow {\rm Sym}^{kk'}(S\times B)$}\smallskip\noindent
\centerline{${\rm Sym}^{kk'}(p\times Id_{B}):{\rm Sym}^{kk'}(S\times
B)\rightarrow {\rm Sym}^{kk'}(B'\times B)$}\smallskip\noindent o\`u
$p:S\rightarrow B'$ est la compos\'ee de la projection sur $B'$ avec
le plongement de $S$ dans $T\times U'\times B'$. $h$ est manifestement  analytique  puisqu'il est associ\'e au morphisme trace
(d\'eduit de ${\cal T}^{0}_{f}$ et ${\cal T}^{0}_{g}\otimes
Id_{U}$):
$${\cal T}_{h}:(\pi_{2}o\pi_{1})_{*}({\cal O}_{X})\rightarrow {\cal O}_{T\times U'\times U} $$
La faible normalit\'e de l'espace des param\`etres assure
l'analyticit\'e du changement de projection ([B1], {\it th\'eor\`eme
2}, p.42) et donc le r\'esultat.
$\,\,\,\blacksquare$\smallskip\noindent {\bf 1.3.8.2.3.
Remarques.}\par\noindent Dans ce qui pr\'ec\`ede, nous avons
implicitement utiliser les propri\'et\'es fondamentales des
morphismes traces :\par\noindent $\bullet$ Si $f:X\rightarrow S$ et
$g:S\rightarrow T$ sont deux morphismes finis et g\'eom\'etriquement
plats, la famille $(f(g^{-1}(U)))_{U\in T}$ index\'ees par les
ouverts de $T$ est une base d'ouverts de $S$\smallskip\noindent
$\bullet$ la compos\'ee des traces est la trace des
compos\'ees.\smallskip\noindent $\bullet$ Par ailleurs la
pond\'eration naturellement induite par la compos\'ee des morphismes
$\pi_{1}$ et $\pi_{2}$ est d\'efinie par une intersection de cycles
d\'ecrites de la fa\c con suivante :\smallskip\noindent Soit ${\cal
X}_{1}$ (resp. ${\cal X}_{2}$) la pond\'eration de $\pi_{1}$ (resp.
$\pi_{2}$) qui est un cycle de $X\times S$ (resp. $S\times T$) et
consid\'erons le diagramme
$$\xymatrix{&&X\times S\times T\ar[ld]_{p_{2}}\ar[d]^{p}\ar[rd]^{p_{1}}&&\\
&S\times T\ar[d]&X\times T\ar[l]_{\pi_{1}\times id}\ar[d]&X\times S\ar[l]_{id\times \pi_{2}}\ar[d]&\\
&S\ar[r]_{\pi_{2}}&T&S\ar[l]^{\pi_{2}}}$$ Alors la pond\'eration de
$\pi$ est donn\'ee par le cycle ${\cal X}$ de $X\times S$
v\'erifiant
$$p^{*}({\cal X}) = p^{*}_{1}({\cal X}_{1})\bullet p^{*}_{2}({\cal X}_{2})$$
ou bien (\`a \'equivalence analytique pr\`es)
$${\cal X} =  p^{*}_{1}({\cal X}_{1})\bullet p^{*}_{2}({\cal X}_{2})\cap [S\times  S]$$
La condition d'incidence, qui est  naturellement  satisfaite, permet
de donner un sens \`a ce cap produit.\smallskip\noindent On peut
sensiblement am\'eliorer la proposition pr\'ec\'edente en donnant la
\Prop{2}{} Avec les hypoth\`eses de la {\it  proposition 1} et en
supposant $S$ faiblement normal et $T$ seulement r\'eduit. La
compos\'ee est analytiquement g\'eom\'etriquement plate. \dem En
vertu de la {\it proposition 7} de {\bf(1.3.9)}, $T$ est
n\'ecessairement faiblement normale et dans ce cas le r\'esultat
d\'ecoule de la {\it proposition 1}.$\,\,\,\blacksquare$
\smallskip\noindent
{\bf 1.3.9. Quelques petits r\'esultats.}\smallskip\noindent Les
r\'esultats les plus significatifs font l'objet de [KII].
 \Prop{3}{}\par\noindent
{\bf(i)} ([B.M],Prop.2,p.814).  Soit $(X_{s})_{s\in S}$ une famille
analytique de cycles d'un espace complexe $Z$ param\'etr\'ee par un
espace r\'eduit $S$. Alors la fonction
$$\nu: S\times Z\rightarrow {\Bbb N};\,\,(s,z)\rightarrow {\rm mult}_{z}(X_{s})$$
est semi-continue sup\'erieurement pour la topologie de
Zariski.\par\noindent {\bf(ii)}([B.M],Prop.A1,p.839). Soit
$\pi:X\rightarrow S$ un morphisme \'equidimensionnel et ouvert
d'espaces complexes r\'eduits. Alors, l'ensemble des points de $X$
en lesquels $\pi$ n'est pas g\'eom\'etriquement plat est un ouvert
dense (et m\^eme de Zariski !).\smallskip\noindent
\Prop{4}{}\par\noindent {\bf (i)}  ([Si],Prop.5.8) , [B1].  Soit
$\pi:X\rightarrow S$ un morphisme ouvert et surjectif d'espaces
complexes dont les fibres g\'en\'eriques sont de dimension constante
sur $S$ faiblement normal.  Alors $\pi$ est (fortement)
g\'eom\'etriquement plat si et seulement si pour tout point $x\in X$
et une certaine projection (donc pour toute) $\eta:W\rightarrow U$,
$x\in W$ ouvert de $X$ et $V=\pi(W)$, le rev\^etement ramifi\'e
$(\pi,\eta):W\rightarrow V\times U$ est de degr\'e
constant.\par\noindent {\bf(ii)} [Sib],  Soit $\pi:X\rightarrow S$
un morphisme ouvert et surjectif d'un espace complexe $X$ localement
de dimension pure sur un espace complexe  $S$ localement
irr\'eductible. Alors $\pi$ est contin\^ument g\'eom\'etriquement
plat.\smallskip\noindent\rm On remarque que sur une base faiblement
normal, les notions de `` contin\^ument g\'eom\'etriquement plat''
et  `` analytiquement  g\'eom\'etriquement plat'' coincident et se
v\'erifient sur {\bf une } installation locale choisie
!\smallskip\noindent Comme corollaires du {\it th\'eor\`eme 2} de
[KII], on obtient \Prop{5}{}  Soit $\pi:X\rightarrow S$ un morphisme
ouvert et de corang constant, avec $S$ r\'eduit de dimension pure.
 Alors \par\noindent
 {\bf(1)} $\pi$ est contin\^ument g\'eom\'etriquement plat si et seulement si
pour tout point $s_{0}$ de $S$, il existe un voisinage ouvert
$S_{0}$ de $s_{0}$ dans $S$  tel que pour {\bf toute}  installation
locale $S_{0}$-adapt\'ee,
$\xymatrix{X\ar[r]^{f}\ar@/_/[rr]_{\pi}&Y\ar[r]^{p}&S_{0}}$ avec $f$
ouvert, fini et surjectif sur $Y$ lisse sur $S_{0}$, il existe une
application continue $\Psi^{\pi}:Y\rightarrow {\ds{\coprod_{r\geq
0}}}{\rm Sym}^{r}(X_{red})$\par\noindent {\bf(2)} $\pi$ est
analytiquement  g\'eom\'etriquement plat (ou plus simplement
g'eom\'etriquement plat) si et seulement si pour tout point $s_{0}$
de $S$, il existe un voisinage ouvert $S_{0}$ de $s_{0}$ dans $S$
tel que pour {\bf toute}  installation locale $S_{0}$-adapt\'ee
(comme ci-dessus), il existe un entier $k$ (appel\'e {\it degr\'e
g\'eom\'etrique}) et une application analytique
$\Psi^{\pi}:Y\rightarrow {\ds{\coprod_{r\geq 0}}}{\rm
Sym}^{r}(X_{red})$ se factorisant dans le diagramme commutatif
$$\xymatrix{Y\ar[rr]^{\Psi_{\pi}}\ar[rd]&&\ds{\coprod_{r\geq 0}}{\rm Sym}^{r}(X_{red})\\
&{\rm Sym}^{k}(X_{red})\ar[ru]&}$$ \Prop{6}{} ([KII],{\it corollaire
4}, p.35) Soient $\pi:X\rightarrow S$ un morphisme universellement
$n$-\'equidimensionnel d'espaces complexes r\'eduits et
$\ds{\Psi^{\pi}_{gen}: S_{gen}\rightarrow {\cal
C}_{*}(X)\setminus\{0\}}$ l'application naturelle associ\'ee
donn\'ee par
 $\ds{s\rightarrow [\pi^{-1}(s)]:=\sum_{i}\mu_{i,s}.[X_{i,s}]}$ o\`u $X_{i,s}$ d\'esignent
  les composantes irr\'eductibles et les entiers $\mu_{i,s}$ les multiplicit\'es associ\'ees. On a \par\noindent
{\bf(i)} si $S$ est localement irr\'eductible alors $\pi$ est
contin\^ument g\'eom\'etriquement plat.\par\noindent {\bf(ii)} si
$S$ est faiblement normal et $\pi$ contin\^ument g\'eom\'etriquement
plat alors il est analytiquement g\'eom\'etriquement
plat.\par\noindent {\bf(iii)} si $S$ est normal alors $\pi$ est
analytiquement g\'eom\'etriquement plat.\rm\smallskip\noindent
\Prop{7}{} Si $\pi:X\rightarrow S$ est un morphisme
g\'eom\'etriquement plat avec $X$ faiblement  normal  alors $S$ est
faiblement normal.\rm\smallskip\noindent \dem Il nous faut montrer
que toute fonction continue et g\'en\'eriquement holomorphe sur $S$
se prolonge globalement en une fonction holomorphe.\par\noindent
Comme le probl\`eme est de nature locale sur $S$,  $\pi$ \'etant
g\'eom\'etriquement plat, donc universellement \'equidimensionnel
\`a fibres de dimension pure $n$  et le probl\`eme de nature locale
sur $S$, on peut se fixer une factorisation locale $X\rightarrow
S\times U\rightarrow S$ dans laquelle $U$ est un polydisque ouvert
relativement compact  de ${\Bbb C}^{n}$. Il va de soi que le
probl\`eme se ram\`ene  au cas o\`u $\pi$ est  g\'eom\'etriquement
plat et  fini.\par\noindent
 Supposons donc $\pi$ fini et consid\'erons  une fonction  $h$ continue et
g\'en\'eriquement holomorphe sur $S$. Comme  $h$ est  m\'eromorphe
localement born\'ee, elle d\'efinit naturellement  une section du
faisceau ${\cal L}^{0}_{S}$. L'image r\'eciproque de $h$ par $\pi$
est continue et d\'efinit  une section du faisceau  ${\cal
L}^{0}_{X}$ (fonctorialit\'e des faisceaux ${\cal
L}^{k}_{\bullet}$). Mais $X$ \'etant faiblement normal et
$\pi^{*}(h)$  m\'eromorphe continue se prolonge, alors,
holomorphiquement  sur $X$.\par\noindent Alors, utilisant les faits
que la platitude g\'eom\'etrique de $\pi$ est \'equivalente \`a la
donn\'ee d'un morphisme trace ${\cal T}^{0}_{f}: f_{*}({\cal
O}_{X})\rightarrow {\cal O}_{S}$ et que la compos\'ee du pull back
et de la {\it{trace}} donne le morphisme  $deg(f).Id$, il nous est
facile de conclure $\blacksquare$.\smallskip\par\noindent {\bf{1.4.
Analogie avec les morphismes plats.}}\smallskip\noindent Les
morphismes analytiquement g\'eom\'etriquement plats poss\`edent des
propri\'et\'es analogues \`a celles des morphismes plats  et de
Tor-dimension finie entre espaces complexes  {\bf r\'eduits}  \`a
savoir :\smallskip\noindent {\bf(i)} ils sont stables par
restriction ouvertes sur $X$ et $S$. La propri\'et\'e d'\^etre
analytiquement g\'eom\'etriquement plat est une condition ouverte
sur $X$.\smallskip\noindent {\bf(ii)} ils sont stables par
changements  de bases. A noter que les morphismes de Tor-dimension
finie ne v\'erifient cette propri\'et\'e que pour les changements de
bases cohomologiquement transversaux.\smallskip\noindent {\bf(iii)}
En tant que morphismes universellement \'equidimensionnels, ils se
factorisent toujours, au moins localement par rapport \`a chacune de
leurs fibres, en morphismes finis suivis  de projections . Mais
\'etant g\'eom\'etriquement plats (essentiellement par construction
de l'espace des cycles), les  morphismes finis  de la factorisation
locale sont toujours g\'eom\'etriquement plats. Signalons que la
factorisation locale d'un morphisme de Tor-dimension finie poss\`ede
d\'ej\`a cette prorpi\'et\'e d'h\'eridit\'e puisque le morphisme
fini de la factorisation est de Tor-dimension finie mais
g\'en\'eralement fausse dans le cas  plat car sinon les fibres
seraient des espaces de Cohen Macaulay!.\smallskip\noindent
{\bf(iv)} ils sont caract\'eris\'es par l'existence de morphismes
trace dans toute factorisation locale; les morphismes plats ou de
Tor-dimension finie induisent des morphismes traces qui ne suffisent
pas \`a les caract\'eris\'es.\smallskip\noindent {\bf(v)} ils ne
sont pas stables par composition alors que les plats le sont! Cette
lacune constitue l'obstruction majeure \`a l'\'elaboration d'une
th\'eorie de l'intersection avec param\`etres dans ce
cadre.\smallskip\noindent {\bf(vi)} Les morphismes de Tor-dimension
finie d'espaces complexes r\'eduits sont analytiquement
g\'eom\'etriquement plats (cf [KII], {\it corollaire 5},
p.35).\bigskip\noindent {\tite{II.  Morphisme g\'eom\'etriquement
plat et int\'egration de classes de cohomologie sur les
fibres.}}\bigskip\noindent {\bf 2.0.} L'int\'egration de classes de
cohomologie de type $(n,n)$ sur une famille analytique de $n$-cycles
d'un espace complexe donn\'e remonte au moins \`a l'article
inaugurale [A.N] de A. Andr\'eotti et F. Norguet. Certes le cadre y
\'etait assez restrictif puisque quasi-projectivit\'e de l'ambiant
et faible normalit\'e de l'espace des param\`etres \'etaient
impos\'ees, mais force est de constater que la technique
fondamentale \'etait clairement mise en lumi\`ere. La
g\'en\'eralisation se fait principalement en trois  \'etapes. En
utilisant l'espace analytique r\'eduit des cycles compacts d'un
espace complexe, Barlet montre dans [B1] que le r\'esultat principal
de [A.N] est encore vrai sous l'hypoth\`ese de quasi-projectivit\'e
de l'ambiant mais sans aucune hypoth\`ese, autre que r\'eduit, sur
l'espace des param\`etres. L'\'etape interm\'ediaire consiste \`a se
d\'efaire de l'hypoth\`ese lourde de quasi-projectivit\'e en
imposant seulement \`a l'ambiant d'\^etre une vari\'et\'e analytique
complexe; ce qui est fait dans [B3]. L'\'etape d\'ecisive entreprise
par Barlet et Varouchas dans [B.V] consite \`a passer outre la
condition de lissit\'e pr\'ec\'edente et la condition de compacit\'e
sur les cycles en consid\'erant des familles analytiques locales.
Les probl\`emes techniques soulev\'es par cette g\'en\'eralisation
ultime sont fondamentalement  de deux ordres \`a savoir : faire face
au d\'efaut de l'isomorphisme de Dolbeault et comment assurer la
finitude des int\'egrales sur des cycles non compacts. Concernant,
l'int\'egration de classes de type quelconque,  on peut signaler
qu'un th\'eor\`eme optimal a \'et\'e  donn\'e par l'auteur  dans
[K2].\smallskip\noindent La  construction repose essentiellement sur
l'utilisation d'un " bon" d\'ecoupage des classes de cohomologie,
technique apparaissant comme  une cons\'equence importante du Lemme
de Reiffen [R]:
 \Prop{} \par\noindent
Soient $ S\subset X\subset Z$ des sous ensembles analytiques
complexes d'un espace analytique complexe $Z$ tels que $dim S <n$ et
$ X-S$ soit lisse de dimension pure $n$. Alors pour tout faisceau
coh\'erent ${\cal F}$ sur $Z$, on a :\par\noindent i)${\rm
H}^{k}_{c}(V, {\cal F}\mid_{V})= 0$ pour tout $k>n$ et tout ouvert
$V$ de $X$.\par\noindent
 ii) Pour tout recouvrement ouvert
$(V_{\alpha})_{\alpha \in A}$ de $X-S$, le morphisme canonique\par
\centerline{$ \bigoplus_{\alpha\in A} {\rm H}^{n}_{c}(V_{\alpha},
{\cal F}\mid_{V_{\alpha}})\rightarrow {\rm H}^{n}_{c}(X, {\cal
F}\mid_{X}) $} \noindent est surjectif.\par\noindent iii) Soient
$\pi: X\rightarrow S$ est un morphisme surjectif d'espaces
analytiques complexes dont les fibres sont de dimension au plus
\'egales \`a $n$, ${\cal F}$ un faisceau coh\'erent sur $X$ et
$(X_{\alpha})_{\alpha\in A}$ un recouvrement ouvert de $X$. Si
$\pi_{\alpha}$ et ${\cal F}_{\alpha}$ d\'esignent les restrictions
de $\pi$ et ${\cal F}$ \`a $X_{\alpha}$, le morphisme canonique \par
 \centerline{$ \bigoplus_{\alpha\in A}{\rm I}\!{\rm R}^{n}{\pi_{\alpha}}_{!}{{\cal
 F}_{\alpha}}
\rightarrow {\rm I}\!{\rm R}^{n}{\pi_{!}}{\cal F}$} \noindent est
surjectif\rm
\smallskip\noindent
Il en r\'esulte les annulations ${\rm I}\!{\rm R}^{j}\pi_{!}{\cal
F}= 0\,\,\forall\,j>n,\,\,\forall\,{\cal F}\in {\rm C}oh(X)$
assurant l'exactitude \`a droite du foncteur  ${\rm I}\!{\rm
R}^{n}\pi_{!}$ permettant de localiser sur $X$.\par\noindent Un cas
particulier de ces r\'esultats appel\'e le lemme du d\'ecoupage dans
[B.V] a \'et\'e introduit dans  [B3] pour pallier l'abscence du
th\'eor\`eme de Dolbeault-Grothendieck dans le cas
singulier.\par\noindent Rappelons que l'int\'egration d'une forme
$\phi$  de type $(n,n)$ \`a coefficients continus  sur un espace
analytique complexe $Z$ de dimension pure $n$ et de partie
r\'eguli\`ere not\'ee ${\rm Reg}(Z)$  est d\'efinie par
$$\int_{Z}\phi:=\int_{Z_{red}}\phi=\int_{{\rm Reg}(Z)}\phi$$
Cette  int\'egrale  a  un sens, en vertu du  th\'eor\`eme de Lelong
[L],  et d\'efinit un courant $d$-ferm\'e appel\'e {\it courant
d'int\'egration} de $Z$ que l'on  note souvent  $[Z]$.
\smallskip\noindent
{\bf 2.1. Classes de cohomologie et repr\'esentants
$\bar\partial$-ferm\'es sur un espace singulier.}\smallskip\noindent
Si $X$ est un espace analytique complexe d\'enombrable \`a l'infini
et $\Phi$ une famille paracompactifiante de supports de $X$,
rappelons que ${\rm H}^{n,n}_{\Phi}(X):={\rm
Ker}({\bar\partial}|_{{\cal A}^{n,n}_{\Phi}(X)})/{\bar\partial} (
{\cal A}^{n,n-1}_{\Phi}(X))$ d\'esigne  le groupe de
$\bar\partial$-cohomologie de type $(n,n)$ de $X$
 \`a support dans $\Phi$ et $\phi_{n,n}:{\rm H}^{n}_{\Phi}(X, \Omega^{n}_{X})
 \rightarrow{\rm H}^{n,n}_{\Phi}(X)$ le morphisme canonique qui est ni injectif ni surjectif en g\'en\'eral
  mais un isomorphisme (appel\'e isomorphisme de Dolbeault) si $X$ est lisse.
\smallskip\noindent {\bf 2.1.1.  Description en {\v C}ech du  morphisme
${\rm H}^{q}_{\Phi}(X, \Omega^{n}_{X})\rightarrow {\rm
H}^{q}(X,(\Gamma_{\Phi}(X, {\cal A}^{n,\bullet}_{X}
))$}\par\noindent Le proc\'ed\'e technique, qui remonte \`a
Andr\'eotti [A] ou  Harvey [H] est clairement  expos\'e dans [B.V].
Concernant la notion de famille paracompactifiante et de familles
duales, recouvrement adat\'e \`a une famille de supports, on renvoie
le lecteur \`a [A.B] ou [A.K]. Mais  [B.V] et [V1] contiennent  les
d\'efinitions essentielles utilis\'ees ici. Soit $X$ un espace
topologique localement compact, paracompact et compl\`etement
paracompact(i.e tous ces ouverts sont paracompacts).\par\noindent
Soit ${\goth U}$ un recouvrement ouvert de $X$ et ${\cal F}$ un
faisceau de groupes ab\'eliens sur $X$. Soit
$\xi:=(\xi_{\alpha_{0}\cdots\alpha_{q}})\in {\cal C}^{q}({\goth U},
{\cal F})$ une $q$- cochaine et $U$ un ouvert de $X$, on d\'efinit
sa restriction \`a $U$ par
$\xi|_{U}:=\xi_{\alpha_{0}\cdots\alpha_{q}}|_{U_{\alpha_{0}\cdots\alpha_{q}}\cap
U}$  et son support{\footnote{$^{(17)}$}{ Tel qu'il est d\'efini, ce
support est g\'en\'eralement plus grand que celui d\'efini par
Godement[G];  mais  pour un recouvrement ouvert adapt\'e \`a la
famille de support $\Phi$, ces d\'efinitions coincident.}} par
$\ds{\Sigma(\xi):=\overline{\bigcup_{\xi_{\alpha_{0}\cdots\alpha_{q}}
\not= 0} U_{\alpha_{0}}\cap\cdots\cap
U_{\alpha_{q}}}}$\smallskip\noindent Soit $({\cal
A}^{n,\bullet}_{X}, \bar\partial)$ le complexe des formes
diff\'erentielles ind\'efiniment diff\'erentiables  (induite par
plongement local) et\par\noindent
 $\ds{\xymatrix{{\cal K}^{\bullet}_{X}:=\Omega^{n}_{X}
\ar[r]^{i}&{\cal
A}^{n,0}_{X}\ar[r]^{\bar\partial}&\cdots\ar[r]^{\bar\partial}&{\cal
A}^{n,q}_{X}\ar[r]^{\bar\partial}&\cdots}}$ son complexe de
Dolbeault. Alors, un \'el\'ement $\psi^{n,q}\in \Gamma_{\Phi}(X,
{\cal A}^{n,q}_{X})$ est appel\'e {\it{repr\'esentant $\bar\partial$
-ferm\'e}} d'une classe de cohomologie  $\xi$ de ${\rm
H}^{q}_{\Phi}(X,\Omega^{n}_{X})$  s'il existe un recouvrement ouvert
${\goth U}$ de $X$ et un  cocycle de \v Cech $(f, \phi^{n,0},
\phi^{n,1},\cdots,\phi^{n,q-1})\in {\cal Z}^{q}({\goth U},{\cal
K}^{\bullet}_{q-1})$ (${\cal K}^{\bullet}_{q-1}$ \'etant le
tronqu\'e en  degr\'e $q$ du complexe de Dolbeault) v\'erifiant (au
signe pr\`es!)\smallskip $\bullet$ $f$ est un repr\'esentant de {\v
C}ech de $\xi$\par $\bullet$ $\delta(f)=0$,
$i(f)=\delta(\phi^{n,0})$, $\bar\partial(\phi^{n,0})=
\phi^{n,1}$,$\cdots$ et (surtout!) $\bar\partial(\phi^{n,q-1})=
\epsilon(\psi^{n,q})$,
$\bar\partial(\psi^{n,q})=0$.\smallskip\noindent o\`u $\ds{\epsilon:
\Gamma_{\Phi}(X,{\cal A}^{n,q}_{X})\rightarrow {\cal
C}^{0}_{\Phi}({\goth U},{\cal A}^{n,q}_{X})}$ est le morphisme
d'augmentation naturel. Cela revient \`a dire que $\psi^{n,q}$ se
``rel\`eve'', via le morphisme canonique ${\rm I}\!\!{\rm
H}^{q}_{\Phi}({\goth U},{\cal K}^{\bullet}_{q-1} )\rightarrow
\Gamma_{\Phi}(X,{\cal A}^{n,q}_{X})$  en un \'el\'ement dont la
``t\^ete'' est le repr\'esentant de \v Cech de
$\xi$.\smallskip\noindent Les $\phi^{n,r}$ se construisent par le
proc\'ed\'e classique consistant \`a prendre  une partition de
l'unit\'e, $(\rho_{\alpha})$, subordonn\'ee \`a ${\goth U}$ et
consid\'erer  les op\'erateurs\par \centerline{$G^{q}_{r}:{\cal
C}^{q}({\goth U}, \Omega^{n}_{X})\rightarrow {\cal C}^{q-r-1}({\goth
U},{\cal A}^{n,r}_{X})$} d\'efinis par
$${G^{q}_{r}}(f)_{\lambda_{0}\cdots\lambda_{q-r-1}}:=
(-1)^{{r(r+1)\over{2}}
+r(n+q)}\sum_{\alpha_{0}\cdots\alpha_{r}}\rho_{\alpha_{0}}{\bar\partial}\rho_{\alpha_{1}}\wedge\cdots\wedge{\bar\partial}\rho_{\alpha_{r}}\wedge
f_{\alpha_{0}\cdots\alpha_{r}\lambda_{0}\cdots\lambda_{q-r-1}}$$
 Pour $r=q$, la formule garde encore un sens  dans  $\Gamma_{\Phi}(X,{\cal A}^{n,q})$ et il
  est facile de v\'erifier la relation $$\delta G^{q}_{r}(f) + G^{q+1}_{r}(\delta f) =
(-1)^{(n+r)}\bar\partial(G^{q}_{r-1}f)$$
 donnant, en particulier,
pour $q=n$,
$$\psi^{n,n}:=
(-1)^{{n(n+1)\over{2}}}\ds{\sum_{\alpha_{0}\cdots\alpha_{n}}}
\rho_{\alpha_{0}}{\bar\partial}\rho_{\alpha_{1}}\wedge\cdots
\wedge{\bar\partial}\rho_{\alpha_{n}}\wedge
f_{\alpha_{0}\cdots\alpha_{n}}$$ Alors, si
$\ds{W=\sum_{j}n_{j}W_{j}}$ est un $n$-cycle de $X$, on d\'efinit
l'int\'egration globale en posant
$$\int_{W}\xi:=\sum_{j}n_{j}\int_{W_{j}}\xi:=\sum_{j}n_{j}\int_{W_{j}}
\psi^{n,n}$$
\par\noindent
{\bf 2.2. Morphismes adapt\'es \`a l'int\'egration.}\par\noindent
Soient  $\pi:X\rightarrow S$ un morphisme $n$- universellement
ouvert, $\Phi$ la  famille de supports constitu\'ee des sous
ensembles $S$-propres\footnote{$^{(18)}$}{Il est connu que $\Phi$
est une famille paracompactifiante (cf [A.B] ou [B.M] \S {\bf 2})} (ou
$\pi$-propre) de $X$, $\Psi$ et $\Theta$ deux  familles  de supports
de $X$. \par\noindent Par analogie avec [A.B] ou
[A.K]\footnote{$^{(19)}$}{Deux familles de supports $\Psi$ et
$\Theta$
 sont dites  {\it{duales}} au sens de Andr\'eotti- Kaas
si\par\noindent $\ds{F\in \Psi\,\Longleftrightarrow F\cap G\,\, {\rm
compact}\,\,\forall\, G\in \Theta}$, et ainsi, $X$ \'etant un espace
analytique complexe donc localement compact, $\ds{\bigcup_{G\in
\Psi}G= X,\,\, \bigcup_{F\in \Theta}F= }X$},  On dira que
:\par\noindent $\bullet$  $\Psi$ et $\Theta$ sont $S$-{\it duales}
si  $\ds{F\in \Psi\,\Longleftrightarrow F\cap G\,\,{\rm est}\,\,
S-{\rm propre }\,\,\forall\, G\in \Theta}$,\par\noindent $\bullet$
$\Psi$ est {\it adapt\'ee} au morphisme $\pi$ si tout point $s$ de
$S$ admet un voisinage ouvert $V$ tel que $\ds{\bigcup_{s\in
V}\pi^{-1}(s)\cap \Psi}\in \Phi$ . Il peut nous arriver de dire que
l'ouvert $V$ est $(\Psi, \Phi)$-adapt\'e.
\par\noindent
$\bullet$ Une installation locale $\Psi$-{\it{adapt\'ee}} est la
donn\'ee, pour tout point $s$ de $S$,
   d'un diagramme commutatif $(\clubsuit)$ du type
$$\xymatrix{W\ar[d]^{\pi} \ar [dr]_{f}\ar [r]^{\sigma}
 & V\times U\times B\ar[d]^{p}\\
 V&\ar[l]^{q} V\times U}$$
 o\`u  $W$ est un ouvert $X$, $V$ un voisinage ouvert $\Phi$-adapt\'e de $s$ dans $S$, $U$ et $B$ sont des polydisques relativement compacts de ${\Bbb
C}^{n}$ et ${\Bbb C}^{p}$ respectivement, $f$ est fini, surjectif et
ouvert, $\sigma$ un plongement local, $p$ et $q$ les projections
canoniques.
\par\noindent
 Soient $\Xi$ une classe de cohomologie de ${\rm H}^{n}_{\Phi}(X, \Omega^{n}_{X})$ et
  $\xi\in \Gamma_{\Phi}(X,{\cal A}^{n,n})$ un de ses repr\'esentants
  ${\bar\partial}$-ferm\'e (i.e  $\phi_{n,n}(\xi)=\Xi$), on pose
$$\int_{\pi^{-1}(s)}\Xi :=\int_{\pi^{-1}(s)}\xi$$
Pour bien comprendre cette op\'eration, on va commencer par la
\smallskip\noindent
{\bf 2.2.1. Forme locale de l' int\'egration sur les
fibres.}\smallskip Soit $\pi:X\rightarrow S$ un morphisme
$n$-universellement ouvert d'espaces complexes r\'eduits. Comme
$\pi$ est $n$-\'equidimensionnel, il se factorise, par rapport \`a
l'une quelconque de ses fibres, sous la forme $(\clubsuit)$
\par\noindent Il parait \'evident que la trace jouera un r\^ole fondamental dans cette op\'eration qui se d\'ecompose en une int\'egration sur les fibres de dimension $0$ suivi d'une int\'egration le long des fibres d'une projection. Avec les notations (et tout ce qui a \'et\'e sur ces morphismes traces) du \S{\bf (1.0.3)},p.17,  on voit que modulo le choix d'une famille de supports adapt\'ee \`a $\pi$ pour que $\pi^{-1}(s)$ rencontre $U\times B$ en des ferm\'es $B$-propres, on peut raisonnablement d\'efinir l'int\'egration sur les fibres. Plus pr\'ecisement, supposons que $(\pi^{-1}(s)_{s\in S}$ d\'efinisse une famille continue de cycles. Alors, relativement \`a la situation d\'ecrite par le diagramme
($\clubsuit$), on se fixe un point $s\in S$ et un polydisque $U'$ de
$U$ ne rencontrant pas l'ensemble de ramification $ {\rm
R}(\pi^{-1}(s))$ d\'ecrit par  les branches locales, $(f_{j})_{1\leq
j\leq k}$, de $\pi^{-1}(s)$ sur $U'$. Si $\xi\in \Gamma(U\times B,
\Omega^{q}_{U\times B})$, on a, pour toute  forme ${\cal
C}^{\infty}$ \`a support compact, $\alpha\in {\cal
A}^{n-q,n}_{c}(U')$,
$$\int_{\pi^{-1}(s)\cap (U'\times B)}\xi\wedge f^{*}\alpha = \int_{\{s\}\times U'}\sum^{j=k}_{j=1}f^{*}_{j}(\xi)\wedge
 \alpha =\int_{\{s\}\times U'}{\rm T}^{q,0}_{f}(\xi)\wedge \alpha$$
\par\noindent Cette quantit\'e varie analytiquement en $s$ si
et seulement si le morphisme associ\'e \`a la famille de
rev\^etement de ramifi\'es $\ds{F: S\times U' \rightarrow {\rm
Sym}^{k}(B)}$ est analytique.
\smallskip\noindent
{\bf{2.2.2. Int\'egration globale.}}\par\noindent Il suffit
proc\'eder localement sur $X$ en utilisant tout ce qui pr\'ec\`ede.
\smallskip\noindent
{\bf 2.2.3.  Construction des fonctions holomorphes sur l`espace des
$n$ - cycles.}\par\noindent Pour r\'ef\'erences, nous aurons  [B.V],
ou [K2]. La construction proc\`ede en deux \'etapes techniques
consistant \`a localiser, traiter le probl\`eme, puis globaliser par
recollement. Dans  la situation locale, c'est- \`a-dire supposer,
dans un premier temps, $Z= U\times B$, avec $U$ et $B$ polydisques
relativement compacts de ${\Bbb C}^{n}$ et ${\Bbb C}^{p}$
respectivement, fixer $s_{0}$ et un de ses voisinage ouvert $S_{0}$
qui soit adapt\'e \`a un recouvrement localement fini du support
$|\pi^{-1}(s_{0})|$ (par d\'efinition, l'int\'egration ne d\'epend
que du support du cycle), on \'etablit le r\'esultat d'int\'egration
pour les classes de ${\rm H}^{n}_{c\times B}(U\times B,
\Omega^{n}_{U\times B})$, et ce par un calcul en termes de cochaines
de { \v C}ech que l'on globalise gr\^ace \`a Reiffen. Ce r\'esultat
est d'une importance capitale pour la suite puisqu'il permet de
localiser sur $X$ l 'int\'egration sur les cycles qui est
\'evidemment  de nature locale sur $S$. Pour nous ramener au cas
local, on choisit deux  familles  de supports $\Phi$
paracompactifiante  et  $\Psi$ telle que pour tout $s$ dans $S$, il
existe un voisinage ouvert $V$ de $s$ tel que $\ds{\bigcup_{s\in
V}|X_{s}|}\subset \Psi$.\par\noindent
 Alors, \'etant donn\'e  un cycle
$X_{s}$, de support not\'e
 $\mid X_{s}\mid$, et  une classe dans
${\rm H}^{n}_{\Phi}(X, \Omega^{n}_{X/S})$, ${\Phi}$ \'etant la
famille des ferm\'es $S$-propres dans $X$.
 Comme $\ds{\Phi\cap\mid X_{s}\mid}$ est contenue dans la famille des compacts de $\mid X_{s}\mid$,
  cette classe a une image naturelle dans ${\rm H}^{n}_{c}(\mid X_{s}\mid, {\Omega}^{n}_ {X/S} \mid _{
\mid X_{s}\mid}) $ \`a laquelle on applique le lemme de Reiffen.
Remarquons que si $X$ est  paracompact et compl\`etement
paracompact, l'\'etendue $\ds{\bigcup_{F\in \Phi}F}$ de la famille
$\Phi$ est un ouvert paracompact; ce
 qui nous permet de prendre $X$ \'egale \`a cette \'etendue, quitte \`a intersecter avec le
  support du cycle. Dans ce cas, ${\Phi\cap\mid X_{s}\mid} =
\Phi\mid_{\mid X_{s}\mid}$
 est contenue dans la famille des compacts de $X$, par hypoth\`ese
sur les supports des cycles.\smallskip\noindent On choisit donc un
recouvrement ouvert localement fini $(X_{\alpha})_{\alpha\in A}$ de
 $X$ muni de diagrammes d'installations locales
 $$\xymatrix{X_{\alpha}\ar[r]^{\sigma_{\alpha}}\ar[rd]_{\pi_{\alpha}}
 &Z_{{\alpha}}\ar[d]^{p_{\alpha}}\\
&S_{\alpha}}$$ L'\'equidimensionnalit\'e de $\pi$ et le th\'eor\`eme
de param\'etrisation locale (th\'eor\`eme de pr\'eparation de
Weierstrass) permettent de factoriser localement $\pi$ en
l'installant
 dans un diagramme du type ($\clubsuit$)
$$\xymatrix{X_{\alpha}\ar[d]^{\pi_{\alpha}} \ar [dr]_{f_{\alpha} }\ar [r]^{\sigma_{\alpha}}
 & S_{0}\times U_{\alpha}\times B_{\alpha} \ar[d]^{p_{\alpha}}\\
 S_{0}&\ar[l]^{q_{\alpha}} S_{0}\times U_{\alpha}}$$
 avec  $f_{\alpha}$ fini et surjectif.\par\noindent
 Une fois le morphisme bien pr\'epar\'e, on " adapte " le recouvrement
 en fonction des supports de classes de cohomologie.
Dans [B.V] ou [K2], le contexte diff\`ere l\'eg\`erement du cadre
actuel  dans le sens o\`u l'on consid\`ere des familles analytiques
de cycles d'un espace complexe donn\'e $Z$ sur lesquelles on
int\`egre des classes de cohomologie. M\^eme si cela n'apparait pas
clairement, le passage par le  sous espace d'incidence donn\'e par
le support du graphe de cette famille est incontournable et
in\'evitable. Pour mener \`a bien les calculs en termes de cochaines
de \v Cech, on se doit de supposer $Z$ paracompact et compl\`etement
paracompact (ou ``{\it{hereditaly paracompact}}'' selon Bredon
[Br]), pour  garantir un minimum de finitude au niveau des groupes
de cohomologie et par voie de cons\'equence la convergence des
int\'egrales, on consid\`ere  deux familles de supports $\Phi$ et
$\Psi$ ad\'equates sur $Z$.
 De toute \'evidence, l'une va contenir les supports des cycles et l'autre sera choisie de
 sorte que leurs images inverses dans le support du graphe de la famille soit contenue dans la famille des ferm\'es $S$-propres. Comme les calculs se font en \v Cech et que l'on localise \`a tour de bras, on prendra deux familles de supports duales $\Phi$ et $\Psi$ comme pr\'ec\'edemment d\'efinies.\par\noindent
 \smallskip \noindent {\bf 2.3.
Propri\'et\'es fonctorielles de l'int\'egration sur les fibres d'un
morphisme
 g\'eom\'etriquement plat.}\smallskip\noindent Si $\pi:X\rightarrow S$ est un morphisme
 analytiquement $n$-g\'eom\'etriquement plat pond\'er\'e par un
 cycle $\goth{X}$ et $\Phi$ la famille de supports de $X$ constitu\'ees de ferm\'es
sur lesquels la restriction de $\pi$ est propre, on dispose encore
d'une int\'egration compatible avec la structure de cycles que l'on
note $\int_{\pi,\goth{X}}$. En effet, $\Phi$ \'etant
paracompactifiante et admettant un recouvrement ouvert adapt\'e
([A.K]), ce qui pr\'ec\`ede s'adapte sans difficult\'e et donne les
principales propri\'et\'es fonctorielles de $\int_{\pi,\goth{X}}$
(cf [B.M], [K2]) \`a savoir :\smallskip\noindent
 {\bf (1) Compatibilit\'e avec les changements de bases sur $S$:}
:\par\noindent
   Soit $\phi: \tilde{S}\rightarrow S$ un morphisme d'espaces analytiques
complexes r\'eduits et
$$\xymatrix{\tilde{X}\ar[r]^{\Theta}\ar[d]_{\tilde{\pi}}& X\ar[d]^{\pi}\\
\tilde{S}\ar[r]_{\phi}& S}$$ le diagramme de changement de base
associ\'e. Notons   $\tilde{\Phi}:=\Theta^{-1}(\Phi)$ la famille
paracompactifiante d\'eduite de $\Phi$ (cf [Br]). Alors, on a  le
diagramme commutatif suivant:
$$\xymatrix{ {\rm H}^{n}_{\Phi}(X, {\Omega}^{n}_{X/S})\ar
[d]^{\int_{\pi,\goth{X}}}
\ar [r]^{\Theta^{*}}& {\rm H}^{n}_{\tilde{\Phi}}
(\tilde{X}, {\Omega}^{n}_{\tilde{X}/\tilde{S}})\ar[d]^{\int_{\tilde\pi, \Theta^{*}\goth{X}}}\\
 {\rm H}^{0}(S, {\cal O}_{S} )\ar [r]_{\phi^{* }}& {\rm H}^{0}(\tilde{S}, {\cal O}_{\tilde{S}} )}$$
{\bf (2)  Compatibilit\'e avec les inclusions ouvertes sur\
$X$:}\par\noindent Si $U$ est un ouvert de $X$ dont l'image par
$\pi$ est l'ouvert  $S_{U}$ et  $j:U\rightarrow X$ l'injection
naturelle. Alors, on a le  diagramme commutatif:
$$\xymatrix{ {\rm H}^{n}_{\Phi|_{U}}(U, {\Omega}^{n}_{U/S_{U}})\ar
[rd]^{\int_{\pi_{U}, \goth{X}|_{S_{U}\times U}}} \ar [rr]^{j_{!}}&&
{\rm H}^{n}_{\Phi}(X, {\Omega}^{n}_{X/S})
\ar[ld]_{\int_{\pi,\goth{X}}}\\
 &{\rm H}^{0}(S, {\cal O}_{S} )&}$$
dans lequel $\pi_{U}$ (resp.$\Phi|_{U}$\footnote{$^{(20)}$}{ qui est
paracompactifiante puisque $X$ est suppos\'e compl\`etement
paracompact}) d\'esigne la restriction de $\pi$ (resp. $\Phi$) \`a
$U$ et  la fl\`eche horizontale d\'esigne le prolongement  par
z\'ero des sections \`a supports. usuelle.\smallskip\noindent
 {\bf (3)  Compatibilit\'e avec les images directes de cycles}\par\noindent
Soit $\xymatrix{ X\ar[rr]^{f}\ar[rd]_{\pi}&&Z\ar[ld]^{\tilde{\pi}}\\
&S&}$  un diagramme commutatif de $S$-morphismes d'espaces complexes
r\'eduits dans lequel $\pi$ est analytiquement
$n$-g\'eom\'etriquement plat pond\'er\'e par $\goth X$,
$\tilde{\pi}$ universellement $n$-\'equidimensionnel, $f$ propre et
surjectif. Soient $\Phi$ (resp. $\tilde\Phi$) une famille
paracompactifiante sur $X$ (resp. $Z$) $S$-adapt\'ee c'est-\`a-dire
telles  que $f^{-1}\tilde\Phi$ soit paracompactifiante et adapt\'ee
\`a $\pi$ ou que $f^{-1}\tilde\Phi$ soit contenue dans $\Phi$. Alors
$\tilde{\pi}$ est analytiquement $n$-g\'eom\'etriquement plat,
pond\'er\'e par le cycle relatif $(f\times Id_{S})_{*}{\goth X}$ et
on a le diagramme commutatif :
  $$\xymatrix{ {\rm H}^{n}_{\tilde\Phi}(Z, {\Omega}^{n}_{Z/S})
\ar[dr]_{\int_{\tilde{\pi},(f\times Id_{S})_{*}{\goth X}}}
\ar[rr]^{f^{*}} &&
 {\rm H}^{n}_{\Phi}(X, {\Omega}^{n}_{X/S})\ar[dl]^{\int_{\pi,\goth{X}}}\\
 &{\rm H}^{0}(S, {\cal O}_{S} )&}$$
\noindent Tout ceci s'interpr\`ete facilement en terme d'image
directe de faisceaux.\par\noindent Signalons, \`a pr\'esent, une cons\'equence du  r\'esultat g\'en\'eral et optimal de [K1]  concernant cette int\'egration.\par\noindent
Soient  $\pi: X\rightarrow S$  un  morphisme analytiquement g\'eom\'etriquement plat et $\Phi$ une famille paracompactifiante de supports adapt\'ee \`a $\pi$ (impliquant, en particulier, que $\forall s\in S$, $\pi^{-1}(s)\cap \Phi$ est contenue dans la famille des compacts de $X$). Alors,  il existe
un (unique)  morphisme d'int\'egration d'ordre sup\'erieur sur les cycles
$$\tilde{\sigma}^{q,0}_{{\Phi}}:{\rm H}^{n}_{\Phi}(X, {\cal
L}^{n+q}_{X})\longrightarrow {\rm H}^{0}(S, {\cal L}^{q}_{S})$$
v\'erifiant les propri\'et\'es fonctorielles cit\'ees plus haut et
 qui, pour $q=0$, rend commutatif le
 diagramme
$$\xymatrix{ {\rm H}^{n}_{\Phi}(X, {\Omega}^{n}_{X/S})\ar
[d]^{\int_{\pi,\goth{X}}}
\ar [r]& {\rm H}^{n}_{\Phi}(X, {\cal L}^{n}_{X})\ar[d]^{\int_{\pi,\goth{X}}}\\
 {\rm H}^{0}(S, {\cal O}_{S} )\ar [r]& {\rm H}^{0}(S,
 {\cal L}^{0}_{S} )}$$\smallskip\noindent
 {\bf 2.4. Remarques.}\par\noindent
{\bf(i)} On pourrait penser que la platitude g\'eom\'etrique analytique de $\pi$ forcerait l'inclusion ${\rm Im}(\tilde{\sigma}^{0,0}_{{\Phi}})\subset{\rm H}^{0}(S, {\cal O}_{S})$. Mais il n'en est rien comme le montre l'exemple trivail
D'ailleurs, consid\'erer simplement $S:=\{(x,y,z)\in{\Bbb C}^{3}/
x^{2}=zy^{2}\}$, $X:=S\times{\Bbb C}_{t}$ et la forme m\'eromorphe
$\ds{x\over{y}}dt$ qui est une section du faisceau ${\cal
L}^{1}_{X}$. En prenant, comme de coutume, une fonction ${\cal
C}^{\infty}$ \`a support compact dans ${\Bbb C}$, $\rho$, et
int\'egrant $\ds{\rho.{x\over{y}}dtd\bar{t}}$, on trouve
trivialement $\ds{x\over{y}}$. Cette fonction est m\'eromorphe
localement born\'ee mais ni holomorphe ni m\^eme continue sur $S$!\par\noindent
On a donc des inclusions strictes
$${\rm H}^{0}(S, {\cal O}_{S})\subset {\rm H}^{0}(S, {\cal
O}^{c}_{S})\subset {\rm Im}\tilde{\sigma}^{0,0}_{{\Phi}}$$
D'ailleurs, si tel \'etait le cas, la compos\'ee de morphismes g\'eom\'etriquement plats serait automatiquement g\'eom\'etriquement plate. Or ceci est g\'en\'eralement faux (cf {\bf(1.1)},{\bf(iii)}, p.29).\par\noindent
{\bf(ii)} En fait,
cette int\'egration g\'en\'eralis\'ee \`a ce type de forme peut
\^etre d\'efinie pour n'importe quel morphisme d'espaces analytiques
complexes r\'eduits ayant des fibres de dimension born\'ees par un
certain entier $n$. Pour le voir, prenons, pour simplifier, un
morphisme propre $\pi:X\rightarrow S$ de dimension relative
g\'en\'erique $n$. Alors, le th\'eor\`eme d'applatissement
g\'eom\'etrique de [B5] montre l'existence d'une modification
$\eta:\tilde{S}\rightarrow S$  et d'un diagramme commutatif
d'espaces complexes
$$\xymatrix{\tilde{X}\ar[r]_{\Theta}\ar[d]_{\tilde{\pi}}&X\ar[d]^{\pi}\\
\tilde{S}\ar[r]_{\eta}&S}$$ avec $\tilde{\pi}$ fortement
analytiquement g\'eom\'etriquement plat (i.e g\'eom\'etriquement
plat au sens de [B5]). En utilisant l'image r\'eciproque
$\Theta^{*}{\cal L}^{n+q}_{X}\rightarrow{\cal L}^{n+q}_{\tilde X}$
et l'isomorphisme $\eta_{*}{\cal L}^{q}_{\tilde{S}}\simeq {\cal
L}^{q}_{S}$, il est facile de produire un morphisme canonique de
faisceaux coh\'erents
$${\rm I}\!{\rm R}^{n}{\pi}_{*}{\cal L}^{n+q}_{X}\rightarrow {\cal
L}^{q}_{S}$$ Ce ph\'enom\`ene s'explique par le fait que les
faisceaux de ce type de formes m\'eromorphes  sont de nature
purement globale et ne peuvent d\'etecter la platitude
g\'eom\'etrique analytique (ou continue) qui est une notion locale!
\bigskip\bigskip\noindent
\tite{III. Faisceaux dualisants, formes r\'eguli\`eres et\par
m\'eromorphes r\'eguli\`eres.} \rm\smallskip\bigskip\noindent {\bf
3.0. Faisceaux dualisants.}\par\noindent Il convient de distinguer
le cas alg\'ebrique ou analytique compact
 du cas analytique g\'en\'eral.\smallskip\noindent
{\bf 3.0.0.}   Pour un  apper\c
cu historique plus complet  (bien que pr\'esent\'e sous un  angle
purement alg\'ebrique), nous  renvoyons le lecteur \`a
[Li].\par\noindent Rappelons qu' un {\it {module}} {\it { dualisant
}}, sur une vari\'et\'e alg\'ebrique projective $X$ de dimension $n$
d\'efinie sur un corps alg\'ebriquement clos de caract\'eristique
nulle $\bf{k}$ (resp. un espace analytique complexe compact  $X$
 de dimension $n$), est la donn\'ee d'un ${\cal O}_{X}$- module
coh\'erent ${\cal K}_{X}$ et d'un  morphisme $\bf{k}$-lin\'eaire
appel\'e {\it{trace}} $\int_{X}: {\rm H}^{n}(X, {\cal
K}_{X})\rightarrow {\bf{k}}\,\, ({\rm resp.} {\Bbb C})$ assurant
l'isomorphisme de foncteurs de la cat\'egorie des ${\cal O}_{X}$-
modules coh\'erents sur la cat\'egorie des $\bf{k}$ (resp. ${\Bbb
C}$) -espaces vectoriels donn\'e par $${\rm H}om({\cal F}, {\cal
K}_{X})\simeq {\rm H}^{n}(X, {\cal F})^{*}={\rm H}om({\rm H}^{n}(X,
{\cal F}), {\bf{k}})$$ Le couple $({\cal K}_{X} , \int_{X})$,
appel\'e {\it paire dualisante} (unique \`a isomorphisme canonique
pr\`es), repr\'esente le foncteur ${\cal F}\rightarrow {\rm
H}^{n}(X, {\cal F})^{*}$ et $\int_{X}$ correspond \`a l'identit\'e
dans ${\rm H}om({\cal K}_{X}, {\cal K}_{X})$.
\smallskip\noindent
Un exemple absolument fondamental a \'et\'e donn\'e par {\bf Grothendieck}. En effet, il montre dans [G1] que  toute
vari\'et\'e projective (pas n\'ecessairement normale) admet une
paire dualisante dans laquelle le faisceau ${\cal K}_{X}$ not\'e
$\omega^{n}_{X}$ est l'unique faisceau coh\'erent de profondeur au
moins deux sur $X$ qui, pour tout plongement $\sigma$, de $X$ dans
un espace projectif ${\Bbb P}_{N}$, coincide avec  $\sigma^{*}{\cal
E}xt^{N-n}_{{\cal O}_{{\Bbb P}_{N}}}({\cal O}_{X}, \Omega^{N}_{{\Bbb
P}_{N}})$. Cette construction s'\'etend sans aucune difficult\'es au cas
de la geom\'etrie analytique complexe comme peut s'en convaincre ais\'ement
le lecteur en consultant le texte clair et  instructif de Monique
Lejeune Jalabert ([LJ]).\smallskip\noindent
{\bf 3.0.1.  Faisceaux dualisants de Andr\'eotti-Kas-Golovin.}\smallskip\noindent
Rappelons que Ramis et Ruget dans  [R.R] ont construit, pour tout espaces analytique complexe $X$, un complexe  ${\cal D}^{\bullet}_{X}$ de ${\cal O}_{X}$-modules  v\'erifiant les propri\'et\'es suivantes:\par\noindent
{\bf(i)} Si $X$ est une vari\'et\'e de dimension $n$, les fibres de ses composantes, en chaque point $x$ de $X$, sont des  ${\cal O}_{X,x}$-modules injectifs et ${\cal D}^{\bullet}_{X}$ est une r\'esolution de $\Omega^{n}_{X}[n]$.\par\noindent
{\bf(ii)} Pour toute immersion ferm\'ee $\sigma:X\rightarrow Z$ d'espaces complexes, on a  ${\cal D}^{\bullet}_{X}\simeq {\cal H}om_{{\cal O}_{Y}}(\sigma_{*}{\cal O}_{X}, {\cal D}^{\bullet}_{Y})|_{X}$.\par\noindent
{\bf(iii)} Il est \`a cohomologie coh\'erente avec
${\cal H}^{k}({\cal D}^{\bullet}_{X})=0$ pour $k<-{\rm Prof}({\cal O}_{X})$.
De plus, si  $X$ est de dimension finie $n$, ${\cal D}^{\bullet}_{X}$ est
 d'amplitude $[-n,0]$ et muni d'une application ${\Bbb C}$-lin\'eaire appel\'ee {\it trace}, ${\rm T}_{X}:{\rm H}^{0}_{c}(X, {\cal D}^{\bullet}_{X})\rightarrow {\Bbb C}$
\smallskip\noindent
Si $X$ est
  de dimension $n$, le $n$-\`eme faisceau d'homologie de ce complexe est
   exactement le faisceau de Grothendieck  (i.e $\omega^{n}_{X}={\cal H}^{-n}({\cal
  D}^{\bullet}_{X})$). Le fait qu'il soit dualisant au
  sens de la g\'eom\'etrie analytique complexe rel\`eve de r\'esultats
  non triviaux de  Andr\'eotti-Kas  et Golovin ([A.K],[G]). On
  profite de cette occasion pour rappeler la notion de faisceau
  dualisant dans ce cadre.\par\noindent Soit $X$ un espace complexe et ${\cal
F}\in{\rm Coh}(X)$. Soit ${\goth U}$ la famille des ouverts de Stein
relativement compacts dans $X$. Alors, \`a toute inclusion
$U'\subset U$ d'ouverts de ${\goth U}$ est associ\'e, en tout
degr\'e $k$,  un morphisme continu d'espaces vectoriels topologiques
de type {\bf D.F.S} $\ds{\rho^{k}_{U',U}:{\rm H}^{k}_{c}(U', {\cal
F})\rightarrow {\rm H}^{k}_{c}(U, {\cal F})}$. Pour tout entier $k$,
on d\'esigne par ${\cal D}^{k}({\cal F})$ le faisceau associ\'e au
pr\'efaisceau\par \centerline{$U\in{\goth U}\rightarrow{\rm {le\,
dual\, topologique\, de}}\,{\rm H}^{k}_{c}(U, {\cal F})$}\par
\centerline{$U'\subset U\rightarrow {\rm {le\, transpos\acute{e}\,
du\, morphisme}}\,\rho^{k}_{U',U}$}\smallskip\noindent Andr\'eotti
et Kas montrent dans [A.K] qu'\`a tout faisceau coh\'erent ${\cal
F}$ sur l'espace complexe $X$, est associ\'e, de fa\c con
fonctorielle, un faisceau analytique ${\cal D}^{k}({\cal F})$
(appel\'ee  {\it faisceau  dualisant} de ${\cal F}$) qui est  ${\cal
O}_{X}$- coh\'erent  et  v\'erifie \par $\bullet$   ${\cal
D}^{k}({\cal F})=0$ si $k\notin[{\rm Prof}({\cal F}), {\rm
Dim}({\cal F})]$,\par $\bullet$   ${\cal D}^{k}({\cal F})\simeq
{\cal E}xt^{-k}_{{\cal O}_{X}}({\cal F}, {\cal D}^{\bullet}_{X})$
et,\par $\bullet$  pour tout $U\in {\goth U}$, $\Gamma(U, {\cal
D}^{k}({\cal F}))$ est isomorphe au dual fort de ${\rm H}^{k}_{c}(U,
{\cal F})$. \smallskip\noindent Golovin ([Go]) g\'en\'eralise ce
r\'esultat au cas o\`u ${\cal F}$ est seulement un faisceau
analytique sur $X$. Il les appellent {\it{faisceaux d'homologie}}
associ\'e \`a ${\cal F}$ et les note  ${\cal H}_{k}({\cal
F})$.\smallskip\par\noindent {\bf 3.1. Formes r\'eguli\`eres et
m\'eromorphes r\'eguli\`eres: Cadre alg\'ebrique }
\smallskip\noindent
{\bf 3.1.0. Cas  absolu.}\par\noindent {\bf 3.1.0.0. Formes
m\'eromorphes r\'eguli\`eres.}\par\noindent {\bf E.  Kunz} montre
dans une s\'erie d'articles ([Ku1], [Ku2], [Ku3])  que toute
vari\'et\'e alg\'ebrique projective $X$ irr\'eductible et de
dimension $n$ porte un faisceau dualisant ${\cal
K}_{X}=\tilde{\omega}^{n}_{X}$ dont la construction se fait par
recollement relativement aux projections finies. Plus
pr\'ecisemment, pour tout morphisme fini et dominant $\pi:
X\rightarrow {\Bbb P}_{n}$, il d\'efinit un faisceau ${\cal O}_{X}$-
coh\'erent $\omega_{\pi}$ donn\'e par l'isomorphisme
$$ \pi_{*}\omega_{\pi} =  {\cal H}om_{{\cal O}_{{\Bbb
P}_{n}}}(\pi_{*}{\cal O}_{X}, \Omega^{n}_{{\Bbb P}_{n}})$$ dont il
montre l'ind\'ependance vis-\`a-vis de la projection finie choisie
en utilisant les traces de diff\'erentielles. Ce faisceau ${\cal
O}_{X}$-coh\'erent $\tilde{\omega}^{n}_{X}$ est de profondeur au
moins deux dans $X$, coincide avec le faisceau des formes usuelles
aux points lisses de $X$ et est muni d'un morphisme canonique
 $ {\cal C}_{X}:\Omega^{n}_{X}\rightarrow \tilde{\omega}^{n}_{X}$
et d'une {\it trace} $\int_{X}: {\rm H}^{n}(X,
\tilde{\omega}^{n}_{X})\rightarrow {\bf k} $.\par\noindent Ce
faisceau est appel\'e le faisceau des {\it formes m\'eromorphes
r\'eguli\`eres}.  Kunz montre que le faisceau
$\tilde{\omega}^{n}_{X}$, s'identifie \`a un sous faisceau du
faisceau des $r$-formes m\'eromorphes enti\`erement  caract\'eris\'e
par la fameuse {\it propri\'et\'e de la trace} disant que  tout
morphisme fini, surjectif et s\'eparable $f:U\rightarrow U'$ d'un
ouvert $U$ de $X$ sur une vari\'et\'e normale $U'$, induit le
diagramme commutatif
$$\xymatrix{f_{*}f^{*}\Omega^{n}_{U'}\ar[r]^{f^{*}}\ar[d]_{{\cal T}^{0}_{f}
\otimes Id}&f_{*}\Omega^{n}_{U}\ar[r]^{f_{*}{\cal C}_{U}}&f_{*}\tilde{\omega}^{n}_{U}\ar[d]^{{\cal T}^{n}_{f}}\\
 \Omega^{n}_{U'}\ar[rr]_{Id}&& \Omega^{n}_{U'}}$$
Par un calcul explicite, Kunz et Waldi (cf [K.W]) \'etendent cette
construction au cas  des $r$- formes en d\'efinissant le faisceau
des $r$- formes  m\'eromorphes r\'eguli\`eres et
 mettant en \'evidence l'isomorphisme canonique $\tilde{\omega}^{r}_{X}\simeq{\cal H}om(\Omega^{n-r}_{X}, \tilde{\omega}^{n}_{X})$.
On peut remarquer que pour tout morphisme $\pi: X\rightarrow Y$ fini
et surjectif de vari\'et\'es compactes de dimension $n$ tel que $X$
soit muni d'un faisceau dualisant ${\cal F}_{X}$, le faisceau ${\cal
H}om_{{\cal O}_{Y}}(\pi_{*}{\cal O}_{X}, {\cal F}_{X})$ est aussi
dualisant. Signalons, d'une part,  que Kunz g\'en\'eralise sa construction au cas
d'une vari\'et\'e alg\'ebrique de Cohen - Macaulay et que, d'autre part,
Lipman ([Li]) montre que $\tilde{\omega}^{n}_{X}$ est encore
dualisant si le corps de base est seulement suppos\'e parfait.
\smallskip\noindent
{\bf 3.1.0.1. Formes r\'eguli\`eres.}\par\noindent
Ces formes \'emanent naturellement d'une th\'eorie de la dualit\'e. A ce titre, le faisceau de Grothendieck est l'exemple type de faisceau de formes r\'eguli\`eres en degr\'e maximal.\par\noindent
On les retrouvent dans le travail de {\bf Elzein} ([E]) bien qu'elles ne soient pas clairement mises  en
\'evidence. En effet, utilisant le complexe r\'esiduel ${\rm
K}^{\bullet}_{X}$  ( c'est-\`a-dire un complexe de ${\cal
O}_{X}$-modules injectifs d'amplitude $[-n, 0]$ et dont l'image dans
la cat\'egorie d\'eriv\'ee est le {\it{ complexe dualisant}}) d'un
sch\'ema $X$ de dimension $n$ et de type fini sur un corps de
caract\'eristique nulle ${\bf k}$, il montre que  le bicomplexe
${\cal K}^{\bullet, \star}_{X}= {\cal H}om(\Omega^{\bullet}_{X},
{\rm K}^{\star}_{X})$ peut \^etre muni d'une structure de complexe
diff\'erentiel de $(\Omega^{\bullet}_{X}, d)$-modules \`a droite
dot\'e de deux diff\'erentielles $\delta$ et $d^{'}_{X}$, la
premi\`ere \'etant naturellement induite par celle de ${\rm
K}^{\bullet}_{X}$, la seconde est, certes, d\'eduite de la
diff\'erentielle ext\'erieure usuelle du complexe de de Rham (qui
n'est pas ${\cal O}_{X}$-lin\'eaire!) mais pas dans le sens naif. Il
est, alors, facile de d\'eduire de cette construction que les
faisceaux coh\'erents $\omega^{\bullet}_{X}= {\cal
H}^{\bullet,0}({\cal K}^{\bullet, \star}_{X}[-n,-n])$ satisfont,
pour tout $r\leq n$, $\omega^{r}_{X} \simeq {\cal H}om_{{\cal O}_{X}}(\Omega^{n-r}_{X},
\omega^{n}_{X})$.
 Ces faisceaux que l'on peut appeler faisceaux des {\it{formes
diff\'erentielles r\'eguli\`eres}} sont tous de profondeur au moins
deux dans $X$ puisque le ${\cal O}_{X}$-module coh\'erent
$\omega^{n}_{X}:={\cal H}^{0}({\rm K}^{\bullet}_{X}[-n])$, qui est
le faisceau de Grothendieck, l'est. Dans le cas r\'eduit de
dimension pure, ils coincident avec les faisceaux de
Kunz.\smallskip\par\noindent {\bf 3.1.1. Cadre alg\'ebrique
relatif.}\smallskip\noindent {\bf 3.1.1.0.}  Rappelons (cf [Ha]) que
si $\pi:X\rightarrow S$ est un morphisme de type fini de sch\'emas
localement noeth\'eriens, le {\it complexe dualisant relatif} ${\cal
D}^{\bullet}_{X/S}$ est un objet de la cat\'egorie d\'eriv\'ee  des
${\cal O}_{X}$-modules \`a cohomologie coh\'erente, d\'efini \`a
quasi-isomorphisme pr\`es et v\'erifiant :
\par {\bf(i)} Pour tout $s$ fix\'e, sa restriction \`a la fibre
$\pi^{-1}(s)$ est le complexe dualisant de cette fibre.\par
{\bf(ii)} Si $\pi$ est lisse, ${\cal D}^{\bullet}_{X/S}$ est
quasi-isomorphe au complexe $\Omega^{n}_{X/S}[n]$.\par {\bf(iii)}Si
$\pi$ se factorise en le diagramme commutatif
$\xymatrix{X\ar[r]^{f}\ar@/_/[rr]_{\pi}&Y\ar[r]^{q}&S}$ o\`u $f$ est
fini et $q$ lisse de dimension relative $n$ sur $S$, $f_{*}{\cal D}^{\bullet}_{X/S}\simeq {\rm I}\!{\rm R}{\rm H}om(f_{*}{\cal O}_{X}, \Omega^{n}_{Y/S}[n])$.
\smallskip\noindent
Si $X$ et $S$ admettent des complexes r\'esiduels ${\rm K}^{\bullet}_{X}$ et ${\rm  K}^{\bullet}_{S}$ respectivement, on
note souvent  ${\cal
 D}^{\bullet}_{X/S}:=\pi^{!}{\cal O}_{S}={\rm I}\!{\rm R}{\rm H}om({\rm I}\!{\rm L}
 \pi^{*}{\rm K}^{\bullet}_{S}, {\rm
 K}^{\bullet}_{X})$ qui, pour $\pi$ plat, admet la repr\'esentation canonique
 compos\'ee de faisceaux flasques ${\cal
 K}^{\bullet}_{X/S}={\rm H}om_{{\cal O}_{X}}(\pi^{*}{\rm
 K}^{\bullet}_{S}, {\rm
 K}^{\bullet}_{X})$. \smallskip\noindent
{\bf 3.2.1.1 . Formes m\'eromorphes r\'eguli\`eres relatives.}\smallskip\noindent Si $X$
et $S$ sont deux sch\'emas noeth\'eriens avec $X$ sans composantes
immerg\'ees (par exemple r\'eduit) et $\pi: X\rightarrow S$ un
morphisme de type fini, $n$-\'equidimensionnel et g\'en\'eriquement
lisse, Kunz et Waldi ont construit dans [K.W], pour tout entier
$k\leq n$, un sous faisceau $\tilde{\omega}^{k}_{X/S}$ du faisceau
des formes m\'eromorphes, dont les sections sont dites  {\it
m\'eromorphes r\'eguli\`eres relatives}.
\par\noindent
On peut remarquer que, si ${\bf k}$ est un corps de
caract\'eristique nulle  et $\pi$ un morphisme  $n$-
\'equidimensionnel et g\'en\'eriquement lisse de ${\bf k}$-sch\'emas
noeth\'eriens sans composantes immerg\'ees et munis de complexes
r\'esiduels, le $n$-\`eme faisceau d'homologie du dualisant relatif
$\omega^{n}_{X/S}:={\cal H}^{-n}(\pi^{!}({\cal O}_{S}))$ est un
faisceau ${\cal O}_{X}$-coh\'erent de Kunz- Waldi. En effet, par
hypoth\`ese sur $X$ et $S$, il est sans torsion. De plus, $\pi$
\'etant g\'en\'eriquement  lisse,  $\omega^{n}_{X/S}$ s'identifie \`a un sous faisceau
du faisceau des formes m\'eromorphes relatives. Il suffit, alors, de montrer que ce faisceau v\'erifie la {\it propri\'et\'e de la trace relative} disant que pour  tout
morphisme fini ou quasi-fini , surjectif et s\'eparable $f:U\rightarrow U'$ d'un
ouvert $U$ de $X$ sur une vari\'et\'e projective lisse sur $S$ et de dimension relative $n$ , on a le diagramme commutatif
$$\xymatrix{f_{*}f^{*}\Omega^{n}_{U'/S}\ar[r]\ar[d]_{{\cal T}^{0}_{f}
\otimes Id}&f_{*}\tilde{\omega}^{n}_{U/S}\ar[d]^{{\cal T}^{n}_{f}}\\
 \Omega^{n}_{U'/S}\ar[r]_{Id}& \Omega^{n}_{U'/S}}$$
Si toutefois l'on ne dispose pas de normalisation de Noether
relative, on peut toujours recouvrir  $X$ par des ouverts de Zariski
munis de morphismes dominants et quasi-finis sur des espaces
projectifs relatifs et recourir au " Main lemma" de Zariski. On en
d\'eduit donc $\omega^{n}_{X/S}= {\tilde{\omega}}^{n}_{X/S}$. On
peut se r\'ef\'erer \`a [Y] ou [L.S] pour de plus amples d\'etails
mais en prenant garde aux diff\'erentes hypoth\`eses puisque l'un
consid\'ere comme base un sch\'ema r\'egulier et l'autre, un peu
plus g\'en\'eral, un sch\'ema excellent sans composantes immerg\'ees
et v\'erifiant la condition de prolongement ${\bf S}_{2}$ de Serre.
\smallskip\noindent
{\bf 3.1.1.2. Formes r\'eguli\`eres relatives.}\smallskip\noindent
Si $\pi$ est un morphisme plat de vari\'et\'es alg\'ebriques sur un
corps de caract\'eristique nulle, Elzein (prop. p.89) montre que le
${\cal O}_{X}$-module bigradu\'e ${\cal
 K}^{\bullet,\star}_{X/S}={\rm H}om_{{\cal O}_{X}}({\Omega}^{\bullet}_{X/S},
 {\rm K}^{\bullet}_{X/S})$ peut \^etre muni d'une structure de
 complexe
 double dot\'e de deux diff\'erentielles $(d_{X/S}, \delta)$. Tout
 comme dans le cas absolu, on en d\'eduit ais\'ement que les
 composantes du complexe diff\'erentiel
 $\omega^{\bullet}_{X/S}= {\cal
H}^{\bullet,0}({\cal K}^{\bullet, \star}_{X/S}[-n,-n])$ satisfont,
pour tout $r\in[0, {dim}X]$,
$$\omega^{r}_{X/S} = {\cal H}om_{{\cal O}_{X}}(\Omega^{n-r}_{X/S},
\omega^{n}_{X/S})$$ et que l'on peut appeler  faisceaux des
{\it{formes diff\'erentielles r\'eguli\`eres relatives }}. L\`a
encore, si $X$ est r\'eduite sur $S$ et $\pi$ \`a fibres de
dimension pure $n$, ces faisceaux coincident avec les faisceaux
relatifs de Kunz- Waldi.\smallskip\noindent Comme nous l'avions
mentionn\'e dans l'introduction et \`a la lumi\`ere de ce qui sera
dit dans la suite, ces constructions nous sugg\`erent la remarque
importante suivante :\par\noindent {\it Si  $\pi:X\rightarrow S$ est
un morphisme de type fini  de sch\'emas noeth\'eriens r\'eduits (ou
sans composantes immerg\'ees avec $S$ excellent)  sur un corps de
caract\'eristique nulle,  g\'en\'eriquement lisse et
$n$-\'equidimensionnel alors $\tilde{\omega}^{n}_{X/S}$ est de
formation compatible \`a tout changement chande base (entre espaces
de m\^eme nature que celle de $S$!) si et seulement si $\pi$
d\'efinit une famille alg\'ebrique de $n$-cycles param\'etr\'ee par
$S$.}\smallskip\noindent On en d\'eduira  un morphisme canonique
${\cal C}_{X/S}:\Omega^{n}_{X/S}\rightarrow
\tilde{\omega}^{n}_{X/S}$  donnant, en particulier, le morphisme
classe fondamentale de Elzein ([E]) si $\pi$ est
plat.\smallskip\par\noindent
 {\bf 3.2. Formes r\'eguli\`eres et m\'eromorphes r\'eguli\`eres: le cadre analytique complexe.} \smallskip\noindent
{\bf 3.2.0. Cas absolu.}\par\noindent {\bf 3.2.0.0. Propri\'et\'e de
la trace absolue.}\par\noindent Elle a d\'ej\`a  \'et\'e
rencontr\'ee dans {\bf(1.0.3)}. Rappelons simplement que, dans le
cadre analytique complexe, le th\'eor\`eme de param\'etrisation
locale (ou de pr\'eparation de Weierstrass) permet toujours de
r\'ealiser localement  un espace analytique de dimension pure $n$
comme un rev\^etement ramifi\'e au dessus d'un certain polydisque
ouvert d'un espace num\'erique ${\Bbb C}^{n}$. Ainsi, tout point $x$
de $X$ admet un voisisnage ouvert $V$ muni d'un morphisme
$f:V\rightarrow U$ fini surjectif et ouvert sur un certain
polydisque ouvert de ${\Bbb C}^{n}$. Notons $j:{\rm
Reg}(X)\rightarrow X$ est l'inclusion naturelle de la partie lisse
de $X$ et $(f_{j})_{1\leq j\leq k}$ les branches locales du
rev\^etement ramifi\'e (que l'on suppose \^etre de degr\'e $k$)
d\'efinit par $f$.\smallskip\noindent
 Alors, un germe $\xi$ de section
en $x$ du faisceau $j_{*}j^{*}\Omega^{n}_{X}$ est dit  v\'erifier la
{\it{ propri\'et\'e de la trace }} absolue  si la forme $\ds{{\cal
T}^{n}_{f}(\xi):=\sum^{j=k}_{j=1}f_{j}^{*}\xi}$, d\'efinie en dehors
de la ramification, se prolonge analytiquement \`a $U$ tout entier
pour tout germe de param\'etrisation locale de $X$ en $x$. (cf
{\bf(1.0.3.1.0)}.\par\noindent Une section du faisceau
$j_{*}j^{*}\Omega^{n}_{X}$, sur un ouvert $V$ de $X$, v\'erifie la
propri\'et\'e de la trace absolue si pour chaque point $x$ de $V$ et
chaque germe de param\'etrisation locale en $x$,  son germe en ce
point v\'erifie la condition pr\'ec\'edente.\par\noindent Plus
g\'en\'eralement (et par  construction),  un germe $\xi$ de section
en $x$ du faisceau $j_{*}j^{*}\Omega^{r}_{X}$ v\'erifie la
propri\'et\'e de la trace si pour tout germe en $x$,  $\alpha$,  de
section du faisceau  $\Omega^{n-r}_{X}$ et tout germe de
param\'erisation locale  $f:V\rightarrow U$ de $X$ en $x$, la forme
m\'eromorphe (holomorphe en dehors de la ramification)
$\ds{\sum_{l=1}^{l= k}{f_{l}}^{*}(\alpha\wedge\sigma)}$
 se prolonge analytiquement sur tout $U$.\smallskip\noindent
{\bf 3.2.0.1.  Formes  m\'eromorphes  r\'eguli\`eres.}\par\noindent
La construction de Kunz s'adapte sans difficult\'es  au cadre analytique complexe donnant pour tout espace analytique complexe r\'eduit de dimension pure $n$ un faisceau ${\cal O}_{X}$-coh\'erent, $\tilde{\omega}^{n}_{X}$,  enti\`erement caract\'eris\'e par la propri\'et\'e de la trace. Il est facile de voir que les r\'esultats pr\'ec\'edemment cit\'es sont encore valables dans ce cadre \`a savoir qu'il coincide avec le faisceau de Grothendieck, qu'en tout degr\'e $k$ on les faisceaux  $\tilde{\omega}^{k}_{X}={\cal H}om(\Omega^{n-k}_{X}, \tilde{\omega}^{n}_{X})$ sont aussi caract\'eris\'es par la propri\'et\'e de la trace et qu'il est muni d'un morphisme canonique ({\it classe fondamentale}) ${\cal C}_{X}:\Omega^{n}_{X}\rightarrow\tilde{\omega}^{n}_{X}$. On peut se r\'ef\'erer aux diff\'erents  articles de Kunz ou \`a ceux  de Kersken ([Ke1], [Ke2]) . \par\noindent
Dans ces \'ecrits apparait d\`ej\`a, de fa\c con plus ou moins explicite, une certaine description de ce faisceau en terme de cette notion de classe fondamentale (sur laquelle nous nous attarderons dans le prochain paragraphe). Barlet exploite cette id\'ee  dans [B2].  En effet, consid\'erant un
plongement local de $X$ (dont le lieu singulier sera not\'e
$\Sigma$) dans une vari\'et\'e lisse $Z$ et notant $j$ l'inclusion
naturelle de la partie lisse de $X$ dans $X$, la classe fondamentale
de $X$ fournit un morphisme {\it cup-produit}
$i_{*}\Omega^{r}_{X}\rightarrow {\cal E}xt^{p}({\cal O}_{X},
\Omega^{p+r}_{Z})$ dont on d\'eduit  un diagramme commutatif
$$\xymatrix{j_{*}j^{*}(\Omega^{r}_{X})\ar[rr]^{\partial}
\ar[rd]_{\tilde{\partial_{k}}}&&{\cal H}^{1}_{\Sigma}
(\Omega^{r}_{X})\ar[ld]^{{\cal H}^{1}_{\Sigma}({\cal C}_{X})}\\
&{\cal H}^{1}_{\Sigma}({\cal E}xt^{p}({\cal O}_{X},
\Omega^{p+r}_{Z}))&}$$ On pose, alors, $\hat{\omega}^{r}_{X}:= {\rm
Ker}\,\tilde{\partial}_{r}$.\par\noindent
 Les propri\'et\'es
intrins\`eques de la classe fondamentale, en particulier son
ind\'ependance vis-\`a-vis du plongement choisi, montre que ce
faisceau est l'incarnation locale d'un faisceau intrins\`eque sur
$X$. Par d\'efinition, ${\rm Ker}\,\tilde{\partial}_{n}$ est le
faisceau de Grothendieck.  Il en r\'esulte  (cf [B2]) que les faisceaux ${\cal O}_{X}$-coh\'erents ${\rm Ker}\,\tilde{\partial}_{r}$ sont ind\'ependants du plongement choisi et sont, de ce fait, intrins\`eques sur tout espace analytique complxe r\'eduit et de dimension pure. Ils sont enti\`erement caract\'eris\'es par la propri\'et\'e de la trace et v\'erifient  $\hat{\omega}^{r}_{X}\simeq {\cal H}om(\Omega^{n-r}_{X},
\omega^{n}_{X})$.
 Il est alors absolument \'evident que ces faisceaux sont exactement les faisceaux  des formes m\'eromorphes r\'eguli\`eres de Kunz. Remarquons, au passage, que,  si ${\cal
C}^{Z}_{X}$ d\'esigne la classe fondamentale de $X$ dans $Z$, [B2]
une $r$- forme g\'en\'eriquement holomorphe  $\sigma$ est une
section de $\omega^{r}_{X}$ si et seulement si $\sigma\wedge {\cal
C}^{Z}_{X}$ se prolonge en  une section globale du faisceau ${\cal
E}xt^{p}_{{\cal O}_{Z}}({\cal O}_{X},
\Omega^{r+p}_{Z})$.\smallskip\noindent {\bf 3.3.0.2. Formes
r\'eguli\`eres.}\par\noindent Le r\^ole principal est tenu par le
faisceau de Grothendieck qui est encore dualisant au sens de la
g\'eom\'etrie analytique complexe comme le pr\'ecise la \Prop{1}{}
Soit $X$ un espace analytique complexe de dimension finie. Alors, il
existe un unique faisceau coh\'erent, $\omega^{n}_{X}$,  sur $X$
v\'erifiant:\par\noindent {\bf(i)} pour tout plongement local
$\sigma:X\rightarrow Z$ dans un ouvert de Stein d'une vari\'et\'e
analytique complexe de dimension $N$, $\omega^{n}_{X}\simeq
\sigma^{*}{\cal E}xt^{N-n}_{{\cal O}_{Z}}({\cal O}_{X},
\Omega^{N}_{Z})$.\par\noindent {\bf(ii)} il est dualisant au sens de
la g\'eom\'etrie complexe, muni d'un morphisme ${\Bbb C}$-lin\'eaire
continu $\ds{\int_{X}: {\rm H}^{n}_{c}(X,\omega^{n}_{X})\rightarrow
{\Bbb C}}$ de sorte que le couple $(\int_{X}, \omega^{n}_{X})$ soit
une paire dualisante c'est-\`a-dire que,  pour tout ouvert $U$ de
Stein  d'un recouvrement ouvert de $X$ et tout faisceau coh\'erent
${\cal F}$ sur $X$, $${\rm I}\!{\rm H}om(U; {\cal F},
\omega^{n}_{X})\simeq {\rm H}^{n}_{c}(U,{\cal F})^{'}$$ (o\`u le
symbole $^{'}$ d\'esigne le dual fort) \dem {\bf(i)} Il est facile
de se convaincre que  le faisceau ${\cal E}xt^{N-n}_{{\cal
O}_{Z}}({\cal O}_{X}, \Omega^{N}_{Z})$ ne d\'epend pas du plongement
choisi. En effet, sans trop entrer dans les d\'etails sur lesquels
nous reviendrons dans [KII], disons que l'on se ram\`ene \`a un
double plongement $X\subset Z_{1}\subset Z_{2}$ muni d'une
r\'etraction de $Z_{2}$ sur $Z_{1}$ et l'on utilise une suite
spectrale classique \`a laquelle on applique les annulations, non
moins classiques des faisceaux ``${\cal E}xt$'' (on peut renvoyer le
lecteur \`a [Go] o\`u \`a [LJ] par exemple).\par\noindent {\bf(ii)}
La construction de la fl\`eche ${\rm H}^{n}_{c}(X,
\omega^{n}_{X})\rightarrow {\Bbb C}$ d\'ecoule du  morphisme
{\it{trace absolu}} d\'efini sur le complexe dualisant de Ramis et
Ruget.  En effet, comme ${\cal H}^{j}({\cal D}^{\bullet}_{X})=0$
pour tout $j<-n$ et que ${\cal H}^{-n}({\cal D}^{\bullet}_{X})=
\omega^{n}_{X}$, on voit que la trace $\ds{{\cal T}: {\rm
H}^{0}_{c}(X,{\cal D}^{\bullet}_{X})\rightarrow {\Bbb C}}$
d\'etermine  enti\`erement le morphisme de faisceaux de groupes
ab\'eliens $\ds{{\cal T}: {\rm
H}^{n}_{c}(X,\omega^{n}_{X})\rightarrow {\Bbb C}}$ puisque la suite
spectrale $\ds{E^{i,j}_{2}= {\rm H}^{i}_{c}(X, {\cal H}^{j}({\cal
D}^{\bullet}_{X}))\Longrightarrow {\rm H}^{i+j}_{c}(X,{\cal
D}^{\bullet}_{X})}$, v\'erifiant $E^{i,j}_{2}=0$ pour $j<-n$, donne
le morphisme lat\'eral $\ds{{\rm H}^{n}_{c}(X, {\cal H}^{-n}({\cal
D}^{\bullet}_{X}))\rightarrow {\rm H}^{0}_{c}(X,{\cal
D}^{\bullet}_{X})}$ que l'on peut d'ailleurs d\'eduire du morphisme
canonique $\ds{{\cal H}^{0}({\cal D}^{\bullet}_{X}[-n])\rightarrow
{\cal D}^{\bullet}_{X}[-n]}$ auquel on applique le foncteur
hypercohomologique exacte \`a droite ${\rm I}\!{\rm
H}^{n}_{c}(X,-)$.\par\noindent Il reste \`a voir qu'il est bien
dualisant au sens de la g\'eom\'etrie analytique complexe. Or cela
r\'esulte imm\'ediatement  de la construction de [A.K] ou [Go]
puisque $\omega^{n}_{X}:={\cal D}^{n}({\cal O}_{X})={\cal
H}_{n}({\cal O}_{X})$. De plus,  pour tout plongement local,
$\sigma$,  de codimension $p$
 de $X$  dans une vari\'et\'e de Stein  lisse $Z$ et tout faisceau
coh\'erent ${\cal F}$ sur $X$, on a   $\ds{{\cal D}^{n}({\cal
F})\simeq {\cal E}xt^{N-n}(i_{*}{\cal F}, \Omega^{N}_{Z})}$ donnant, par recollement,
$${\cal D}^{n}({\cal F})\simeq {\cal
H}om_{{\cal O}_{X}}({\cal F}, \omega^{n}_{X})$$ Mais d'apr\`es [A.K]
ou [Go], pour tout ouvert $U$ d'un recouvrement de $X$,
$\Gamma(U,{\cal D}^{n}({\cal F}))$ est isomorphe au dual fort de
${\rm H}^{n}_{c}(U, {\cal F})$.\par\noindent
Il s'en suit que  $(\int_{X}, \omega^{n}_{X})$ est bien une paire dualisante au sens de la g\'eom\'etrie complexe.$\,\,\,\blacksquare$\smallskip\noindent
Encore une fois, nous conseillons au lecteur de consulter [LJ] pour se faire une id\'ee de l'utilisation des outils de la g\'eom\'etrie analytique complexe dans la construction de la trace par exemple.\smallskip\noindent
Pour les formes r\'eguli\`eres de degr\'es  interm\'ediaires, il est important de noter que la construction de Elzein (cf {\bf (3.1.0.1)} est transposable \`a ce cadre moyennant un
effort consid\'erable effectu\'e par   Kersken ([Ke], [Ke1], [Ke2])
qui r\'eussit le tour de force de d\'efinir,  pour toute alg\`ebre
analytique non n\'ecessairement r\'eduite ${\cal A}$, un complexe
diff\'erentiel de {\it formes r\'eguli\`eres}
$\omega^{\bullet}_{{\cal A}}$ sur l'alg\`ebre diff\'erentielle de de
Rham $(\Omega^{\bullet}_{{\cal A}},d)$ v\'erifiant, en tout degr\'e
$r$, $\omega^{r}_{{\cal A}} = {\cal H}om_{{\cal
A}}(\Omega^{n-r}_{{\cal A}}, \omega^{n}_{{\cal A}})$. Ainsi, par globalisation en vertu du lemme de Grothendieck-Frisch,  pour tout espace analytique complexe de dimension finie $n$, les faisceaux ${\cal O}_{X}$-coh\'erents ${\cal H}om(\Omega^{n-r}_{{X}}, \omega^{n}_{X})$ peuvent \^etre consid\'er\'es comme les faisceaux des formes r\'eguli\`eres de la g\'eom\'etrie alg\'ebrique. Alors,
\Prop{2}{} Si $X$ est un espace analytique complexe r\'eduit de dimension pure $n$, on a,  pour tout entier $r\leq n$,
$$\tilde{\omega}^{r}_{X}=\hat{\omega}^{r}_{X}=\omega^{r}_{X}={\cal
D}^{n}(\Omega^{n-r}_{X})$$
\dem On a d\'ej\`a vu que $\tilde{\omega}^{r}_{X}=\hat{\omega}^{r}_{X}$. L'identification canonique $\tilde{\omega}^{r}_{X}=\omega^{r}_{X}$ d\'ecoule \'evidemment du travail cons\'equent de [Ke]. Enfin, de  [A.K] ou [Go] r\'esulte  l'isomorphisme
$${\cal
D}^{n}(\Omega^{n-r}_{X})\simeq {\cal H}om_{{\cal O}_{X}}
(\Omega^{n-r}_{X}, \omega^{n}_{X})$$ dont on peut se convaincre
ais\'ement en proc\'edant localement puis par recollement sur  $X$.
En effet, supposons ce dernier localement plong\'e en codimension
$p$ dans une vari\'et\'e de Stein  lisse $Z$  de dimension $n+p$,
alors, on a
   $\ds{{\cal D}^{k}({\cal F})\simeq {\cal E}xt^{N-k}(\sigma_{*}{\cal F},
 \Omega^{N}_{Z})}$, pour tout faisceau coh\'erent ${\cal F}$ sur $X$. D'o\`u, en particulier,  pour  $k:=n$ et  ${\cal F}:= \Omega^{n-r}_{X}$,  les faisceaux $\ds{{\cal
E}xt^{N-n}(\sigma_{*}\Omega^{n-r}_{X},
 \Omega^{N}_{Z})}$, qui, en vertu de d\'eg\'en\'erescence de suites spectrales classiques, sont isomorphes aux faisceaux
  ${\cal H}om_{{\cal O}_{Z}}(\sigma_{*}\Omega^{n-r}_{X},{\cal E}xt^{N-n}(\sigma_{*}{\cal O}_{X},
   \Omega^{N}_{Z}))$, qui nous donnent apr\`es recollement sur $X$,
   l'isomorphisme (ou l'identification) canonique d\'esir\'e.
Ces faisceaux coh\'erents ${\cal D}^{n}(\Omega^{n-r}_{X})$ sont des
  faisceaux dualisants au sens de [A.K] et coincident avec les faisceaux  des formes
r\'eguli\`eres.$\,\,\blacksquare$\smallskip\noindent {\bf 3.2.0.3.
Formes m\'eromorphes r\'eguli\`eres  et courants.}\par\noindent On
dispose, gr\^ace \`a [E1], d'un r\'esultat explicite comparant  les
r\'esidus de Grothendieck et de Herrera. En effet, il y est montr\'e
que, pour toute param\'etristation locale d'un espace analytique
complexe de dimension pure $n$
$\ds{\xymatrix{X\ar@/_/[rr]_{f}\ar[r]^{\sigma}&Z:=U\times
B\ar[r]^{q}&U}}$ et toute syst\`eme $(g_{1},\cdots,g_{p})$ de
fonctions holomorphes sur $Z$ et d\'efinissant une suite
r\'eguli\`ere sur $X$ en un point g\'en\'erique $x$. Alors, si l'on
d\'esigne par ${\rm Res}^{G}$ (resp. ${\rm Res}^{H}$) le r\'esidu de
Grothendieck (resp. Herrera), on a, pour toute section $w$ de
$\Gamma(Z,\Omega^{n+p}_{Z})$,
$${\cal T}_{f}\left[{\rm Res}^{H}_{g_{j}}\left({w\over{g_{1}\cdots g_{p}}}
\right) \right] = (2i\pi)^{n}\pi_{*}\left[X\right]\wedge{\rm
Res}^{G}\left[{w\over{g_{1}\cdots g_{p}}}\right]$$ ${\cal T}_{f}$
\'etant l'image directe au sens des courants telle qu'elle a\'et\'e
d\'ecrite dans {\bf(1.0.3.1.0)} (p.17).\smallskip\noindent Par ailleurs,
[B2] met en \'evidence une description du faisceau de Kunz  en terme
de courants ``holomorphes `` sur $X$. Plus pr\'ecisemment, les
consid\'erations locales pr\'ec\'edentes permettent d'\'ecrire toute
$r$-forme m\'eromorphe comme un quotient $\xi= { v\over{g}} $ o\`u
$v$ est une $r$-forme holomorphe sur $Z$ et $g$ une fonction
holomorphe s'annulant sur le lieu singulier $\Sigma$ de $X$.
D'apr\`es  Herrera - Lieberman ([H.L]), $\xi$  d\'efinit un courant
appel\'e  {\bf { valeur principale}} d\'efinit par
 $$\langle{\bf T}_{\xi},\phi\rangle:= lim_{\epsilon\rightarrow 0}\int_{X\cap\{|g|>\epsilon\}}\xi\wedge \phi$$
pour toute forme $\phi$ de type $(n-r,n)$,  $C^{\infty}$ et \`a
support compact dans $X$. Il est, d'ailleurs,  souvent not\'e $\xi\wedge
[X]$.\smallskip\noindent Alors, une $r$-forme m\'eromorphe est une
section du faisceau $\omega^{r}_{X}$ si et seulement si le courant
valeur principale associ\'e est $\bar\partial$- ferm\'e; ce qui nous
permet d'identifier ${\omega}^{r}_{X}$, au faisceau des courants
$\bar{\partial}$ - ferm\'es sur
 $X$ modulo ${\cal O}_{X}$ torsion.\smallskip\par\noindent
{\bf 3.2.1. Cas relatif.}\smallskip\noindent
 {\bf 3.2.1.0. Propri\'et\'e de la trace relative.}\par\noindent
Nous renvoyons le lecteur au pararagraphe \S{\bf(1.0.3.2.1)} (p.21)   dans lequel il trouvera tous les d\'etails qu'il d\'esire sur cette notion. Rappelons simplement que : \par\noindent
\'etant donn\'e un morphisme $\pi:X\rightarrow S$
 universellement $n$-\'equidimensionnel, il admet en chaque
  point $x$ de $X$ une factorisation locale
$\ds{\xymatrix{X\ar@/_/[rr]_{\pi}\ar[r]^{f}&Y:=S\times
U\ar[r]^{q}&S}}$ et $X$ localement plong\'e dans un espace complexe
lisse sur $S$ et du type $S\times U\times B$ avec $U$ (resp. $B$) un
polydisque relativement compact de ${\Bbb C}^{n}$ (resp.${\Bbb
C}^{p}$). Cela nous d\'efinit un rev\^etement ramifi\'e
``g\'en\'erique'' d'un certain degr\'e $k$, de branches locales
$(f_{l}(s,t))_{1\leq l\leq k}$ et de morphisme classifiant $F:
S\times U\rightarrow {\rm Sym}^{k}(B)$. Notons $j$ l'inclusion
naturelle du lieu r\'egulier de $\pi$ dans $X$ (cf
p.37).\par\noindent Alors, \par $\bullet$ si pour tout $s$ de $S$,
$F(\{s\}\times U)$ n'est pas contenue dans le lieu singulier ${\rm
Sym}^{k}(B)$, on dira qu'un germe  de section, en $x$, du faisceau
$j_{*}j^{*}\Omega^{r}_{X/S}$, $\xi$ v\'erifie la {\it propri\'et\'e
de la trace relative} si, pour toute param\'etrisation locale de $X$
en $x$, et tout $\alpha$, germe en $x$ d'une section de
$\Omega^{n-r}_{X/S}$, la forme g\'en\'eriquement holomorphe
$\ds{\sum_{l=1}^{l= k}{f_{l}}^{*}(\alpha\wedge\xi)}$ se prolonge
analytiquement sur $S\times U$ tout entier;\par $\bullet$ sinon, on
demande que ce soit v\'erifi\'ee cette condition sur la composante
g\'en\'erique de $F$ (cf \S{\bf 1.0.3.2}, p.21).\smallskip\noindent
{\bf 3.2.1.1. Formes r\'eguli\`eres relatives.}\par\noindent Si
$\pi:X\rightarrow S$ est un morphisme {\it plat} d'espaces complexes
r\'eduits de dimension finie, la construction de [K.W] s'applique
sans aucune difficult\'e et donne en tout degr\'e $r$, le faisceau
des formes m\'eromorphes r\'eguli\`eres ${\tilde\omega}^{r}_{X/S}$
dont les sections sont caract\'eris\'ees par la {\it propri\'et\'e
de la trace relative}.\par\noindent Toujours sous cet aspect  local
faisant intervenir les alg\`ebres analytiques, Kersken propose dans
[Ke] une m\'ethode compl\`etement diff\'erente reprenant l'id\'ee
qui est \`a la base de la construction des formes r\'eguli\`eres de
Elzein (cf {\bf(3.1.1.2)}). Ce travail utilise les complexes
r\'esiduels et de Cousin relatifs et peut  repr\'esenter un pas non
n\'egligeable vers la dualit\'e analytique relative. Plus
pr\'ecisemment, si $\phi:{\cal P}\rightarrow {\cal A}$ est un
morphisme plat d'alg\`ebres analytiques locales avec ${\cal A}$ non
n\'ecessairement r\'eduite ni de dimension pure
 (ce qui constitue une s\'erieuse  difficult\'e technique!)  muni de l'alg\`ebre diff\'erentielle
$(\Omega^{\bullet}_{{\cal A}/{\cal P}},d)$ des formes
diff\'erentielles relatives dot\'ees  de la diff\'erentielles
relative usuelle, Kersken construit un complexe de Cousin relatif
${\cal C}_{\Omega^{\bullet}}({\cal A}/{\cal P})$ et un complexe
r\'esiduel  diff\'erentiel ${\cal D}_{\Omega}({\cal A}/{\cal P})$ de
$(\Omega^{\bullet}_{{\cal A}/{\cal P}},d)$-modules tel que pour
toute
 ${\cal P}$-param\'etrisation $f:{\cal R}\rightarrow {\cal A}$ de codimension $p$  (i.e ${\cal P}$-morphisme fini
  d'alg\`ebres analytiques locales )
 $${\cal D}_{\Omega}({\cal A}/{\cal P})=
 {\rm H}om_{ \Omega^{\bullet}_{{\cal R}/{\cal P}}}(\Omega^{\bullet}_{{\cal A}/{\cal
 P}}, {\cal C}_{\Omega^{\bullet}}({\cal A}/{\cal P}))[p,p]$$
 et  dont la cohomologie de degr\'e $0$ donn\'ee par ${\rm Ker}: {\cal D}_{\Omega}({\cal A}/{\cal
P})^{0,\star}\rightarrow {\cal D}_{\Omega}({\cal A}/{\cal
P})^{1,\star}$
 est le  $(\Omega^{\bullet}_{{\cal A}/{\cal P}}, d)$- module
des formes (g\'en\'eriquement holomorphes) r\'eguli\`eres
$\omega^{\bullet}_{{\cal A}/{\cal P}}$ (analogue des formes
r\'eguli\`eres du cadre alg\'ebrique).
\par\noindent De ce travail cons\'equent, on peut tirer  le \Cor{1}{} Soit
$\pi:X\rightarrow S$ un morphisme $n$-plat d'espaces
analytiques complexes avec $S$ r\'eduit de dimension pure localement
fini. Alors, il existe un faisceau analytique (unique \`a
isomorphisme canonique pr\`es) $\omega^{n}_{X/S}$ v\'erifiant
:\par\noindent {\bf(i)} il est ${\cal O}_{X}$-coh\'erent et de
profondeur au moins deux fibre par fibre sur $S$,\par \noindent
{\bf (ii)} la famille  de faisceaux coh\'erents   $\omega^{\bullet}_{X/S}={\cal
H}om(\Omega^{n-\bullet}_{X/S}, \omega^{n}_{X/S})$ est munie d'une
diff\'erentielle non triviale ${\rm D}$, faisant de
$(\omega^{\bullet}_{X/S}, {\rm D})$  un complexe diff\'erentiel de
$(\Omega^{\bullet}_{X/S}, d_{X/S})$-modules.\par\noindent
{\bf(iii)}
si $X$ est r\'eduit,  $\omega^{n}_{X/S}$ est le faisceau
$\tilde{\omega}^{n}_{X/S}$ des formes m\'eromorphes r\'eguli\`eres
caract\'eris\'e par la propri\'et\'e de la trace relative au sens de
Kunz-Waldi [K.W].\rm\smallskip\noindent
{\bf 3.2.1.2. Formes m\'eromorphes r\'eguli\`eres.}\par\noindent Comme ces formes \'emanent, par nature,
d'une th\'eorie de la dualit\'e relative, l'occasion nous  est donn\'ee de voir l'apport, dans cette direction, de  la dualit\'e analytique relative. Malheureusement, elle n'est pas encore au point et reste, d'ailleurs' sous une forme incompl\`ete et insatisfaisante dans le cas g\'en\'eral. Toutefois, dans le cas d'un morphisme propre on dispose d'un vrai th\'eor\`eme de dualit\'e relative gr\^ace \`a  Ramis, Ruget et Verdier. En effet, il est montr\'e, dans [R.R.V], que, pour tout morphisme $\pi:X\rightarrow S$  d'espaces analytiques
complexes d\'enombrables
 \`a l'infini de dimension  finie munis de leurs complexes dualisants de
 Ramis-Ruget  ${\cal D}^{\bullet}_{X}$ et  ${\cal D}^{\bullet}_{S}$  , pour tout  complexes  de faisceaux de ${\cal O}_{X}$-modules \`a cohomologie coh\'erente,
   $A^{\bullet}$ et $B^{\bullet}$, ce dernier \'etant \`a cohomologie  born\'ee
    \`a gauche, il existe une fl\`eche canonique
    \footnote{$^{(3)}$}{avec des formules analogues obtenues en \'echangeant ${\rm I}\!{\rm
 R}{\pi}_{!}$ et ${\rm I}\!{\rm
 R}{\pi}_{*}$.} (et fonctorielle en les arguments)
$${\rm I}\!{\rm
 R}{\pi}_{!}{\rm I}\!{\rm R}{\rm H}om(X; A^{\bullet},{\rm I}\!{\rm R}{\rm H}om( {\rm I}\!{\rm
 L}\pi^{*}({\cal D}_{S}(B^{\bullet})), {\cal D}^{\bullet}_{X}))\rightarrow
{\rm I}\!{\rm R}{\rm H}om(S; {\rm I}\!{\rm R}\pi_{*}A^{\bullet},
B^{\bullet}),$$ o\`u ${\cal D}_{S}(B^{\bullet}):={\rm I}\!{\rm
R}{\cal H}om(B^{\bullet}, {\cal D}^{\bullet}_{S})$,  qui est un
isomorphisme (alg\'ebrique) si $\pi$ est propre .\par\noindent
 Alors, en notant  ${\cal D}^{\bullet}_{X/S}:=\pi^{!}({\cal O}_{S}) ={\rm
I}\!{\rm R}{\cal H}om({\rm I}\!{\rm
 L}\pi^{*}{\cal D}^{\bullet}_{S}, {\cal D}^{\bullet}_{X})$, on obtient, gr\^ace \`a la trace relative construite dans [RRV], un morphisme ${\cal O}_{S}$-lin\'eaire continue ${\rm I}\!{\rm
 R}{\pi}_{*}\pi^{!}({\cal O}_{S})\rightarrow {\cal O}_{S}$\smallskip\noindent
 Dans
ce cas, le foncteur $\pi^{!}$ est parfaitement d\'efini,
 v\'erifie les propri\'et\'es d'usage (cf [H], [L1]) et, est l'adjoint \`a droite du
 foncteur ${\rm I}\!{\rm R}\pi_{*}$. Ainsi, [RRV] a pour cons\'equences  imm\'ediates le
 \Cor{2}{} Si $X$ et $S$ sont deux espaces analytiques  complexes d\'enombrables \`a l'infini de dimension finie et $\pi:X\rightarrow S$ un morphisme propre, il existe un foncteur
$$\pi^{!}:{\cal D}^{+}_{coh}(S)\rightarrow {\cal D}^{+}_{coh}(X)$$
muni d'un morphisme de foncteurs ${\rm I}\!{\rm R}\pi_{*}\pi^{!}\rightarrow Id$
 v\'erifiant :\smallskip\noindent {\bf(i)} il est de nature locale
sur $X$ (au sens de Verdier [V]) c'est-\`a-dire pour deux
 morphismes propres  $\pi_{i}:X_{i}\rightarrow S,\,\,i=1,2$ et $U$ un ouvert
de $X_{1}$ muni
  de deux inclusions ouvertes $j_{i}:U\rightarrow X_{i}$ et
  tel que le diagramme
  $$\xymatrix{&U\ar[ld]_{j_{1}}\ar[rd]^{j_{2}}&\\
  X_{1}\ar[rd]_{\pi_{1}}&&X_{2}\ar[ld]^{\pi_{2}}\\
  &S&}$$
  soit commutatif, on a
  $$j_{1}^{*}({\pi_{1}}^{!}F) =
  j_{2}^{*}({\pi_{2}}^{!}F),\,\,\forall\,F\in {\rm C}oh(S)$$
\par\noindent {\bf(ii)} Il est de formation compatible avec la composition des
morphismes propres  d'espaces analytiques complexes d\'enombrables
\`a l'infini c'est-\`a-dire si
$\xymatrix{X\ar[r]^{f}&Z\ar[r]^{g}&S}$, alors $$(gof)^{!}=
f^{!}og^{!}$$\rm\smallskip\noindent qui n'est rien d'autre qu'une
relecture de  [V] dans ce contexte et le \Cor{3}{} Soit
$\pi:X\rightarrow S$ un morphisme universellement
$n$-\'equidimensionnel d'espaces complexes d\'enombrables \`a
l'infini. Alors, le faisceau analytique $\omega^{n}_{X/S}:= {\cal
H}^{-n}(\pi^{!}({\cal O}_{S}))$ v\'erifie les propri\'et\'es
suivantes:\par\noindent {\bf(i)} il est ${\cal O}_{X}$-coh\'erent,
de profondeur au moins deux fibre par fibre sur $S$, de formation
compatible aux inclusions ouvertes sur $X$ et changement de base
plat sur $S$, et coincide avec le faisceau des formes holomorphes
relatives sur la partie r\'eguli\`ere de $\pi$.\par\noindent
{\bf(ii)} il munit $\pi$ d'un morphisme ${\cal O}_{S}$-lin\'eaire
continu $ \int_{\pi}:{\rm I}\!{\rm
 R}^{n}{\pi}_{*}\omega^{n}_{X/S}\rightarrow {\cal O}_{S}$
  de formation compatible \`a tout
  changement de base plat et \`a la composition des morphismes propres
  et universellement \'equidimensionnels dans le sens  o\`u
  tout diagramme commutatif
$$\xymatrix{X_{2}\ar[rr]^{\Psi}\ar[rd]_{\pi_{2}}&&X_{1}\ar[ld]^{\pi_{1}}\\
&S&}$$
 avec $\pi_{1}$, $\Psi$ et $\pi_{2}$ universellement \'equidimensionnels
 et propres de dimension relative respective
 $n_{1}$, $n$ et  $n_{2}:=n_{1}+n$,
 donne un diagramme commutatif de faisceaux coh\'erents
$$\xymatrix{{\rm I}\!{\rm
 R}^{n}{{\pi_{2}}_{*}}\omega^{n_{2}}_{X_{2}/S}\ar[rr]\ar[rd]_{\int_{\pi_{2}}}&&
 {\rm I}\!{\rm
 R}^{n_{1}}{\pi_{1}}_{*}\omega^{n_{1}}_{X_{1}/S}\ar[ld]^{\int_{\pi_{1}}}\\
&{\cal O}_{S}&}$$ \rm
\par\noindent
\dem\par\noindent
{\bf(i)}\par
$\bullet$  La  {\it coh\'erence}  d\'ecoule imm\'ediatement du fait que le complexe
$\pi^{!}{\cal O}_{S}$ est  \`a cohomologie coh\'erente. En effet, le
probl\`eme \'etant de nature locale sur $X$ et  les morphismes
universellement $n$- \'equidimensionnels admettant une r\'ealisation
locale du type
$$\xymatrix{X\ar[r]^{\sigma}\ar[rd]_{\pi}&Z\ar[d]^{q}\\
&S}$$ avec $\sigma$ plongement local, $Z$ lisse sur $S$ de dimension
relative $N$, il est facile de voir que,  pour tout faisceaux
coh\'erents $F$ et $G$,  les isomorphismes (le second \'etant celui de Verdier [V])
$$\sigma_{*}\sigma^{!}F\simeq {\rm I}\!{\rm R}{\cal H}om(\sigma_{*}{\cal
O}_{X}, F),\,\,\,\,{\rm et}\,\,\,q^{!}G\simeq
\Omega^{N}_{Z/S}[N]\otimes q^{*}G$$
 montrent clairement  que
$\pi^{!}{\cal O}_{S}=i^{!}q^{!}{\cal O}_{S}$ est \`a cohomologie
coh\'erente.\par
$\bullet$ Sa nature {\it locale} d\'ecoule directement de celle du
 foncteur  $\pi^{!}$ (cf {\it corollaire 2}) et de la stabilit\'e par changement de base des morphismes universellement \'equidimensionnels.\par
$\bullet$ Pour v\'erifier que ${\rm P}rof(\omega^{n}_{X/S})\geq 2$, on commence par constater que les faisceaux d'homologie du complexe $\pi^{!}({\cal O}_{S})$ sont  tous nuls pour $j<-n$. En effet, si l'on suppose $X$
localement donn\' e par l'annulation de $p$ fonctions holomorphes
sur $Z$, on aura
$${\cal H}^{j}(\pi^{!}{\cal O}_{S})\simeq
{\cal E}xt^{j+N}(\sigma_{*}{\cal O}_{X}, \Omega^{N}_{Z/S})$$ dont
l'annulation est assur\'ee pour $j+N<N-n$ (cf  le {\it lemme 2},p.6 de [KII]).\par\noindent
 On en d'eduit, en
particulier, l'incarnation locale $\omega^{n}_{X/S}\simeq {\cal
E}xt^{N-n}(\sigma_{*}{\cal O}_{X}, \Omega^{N}_{Z/S})$ qui montre bien que  sur la partie ${\rm R}eg(\pi)$, le faisceau $\Omega^{n}_{X/S}$ s'identifie naturellement  avec $\omega^{n}_{X/S}$.\smallskip\noindent
Consid\'erons, \`a pr\'esent,  un sous espace $Y\subset X$  de codimension au moins
$2$ dans $X$ et montrons que
$${\cal H}^{0}_{Y}(\omega^{n}_{X/S})={\cal H}^{1}_{Y}(\omega^{n}_{X/S})=0$$
Pour cela, on peut, soit utiliser la fonctorialit\'e, la stabilit\'e
par changement de base des morphismes universellements
\'equidimensionnels et les annulations de cohomologie du dualisant
relatif, soit proc\'eder localement sur $X$ au voisinage de $Y$ et
se ramener \`a v\'erifier les annulations
 ${\cal H}^{j}_{Y}({\cal E}xt^{N-n}(\sigma_{*}{\cal O}_{X},
\Omega^{N}_{Z/S})=0$ pour $j=0,1$ et  pour lesquels on renvoie le lecteur au  {\it lemme 3}, p.6  de [KII].  \par
$\bullet$ La {\it stabilit\'e par changement de base plat} se justifie en utilisant  les m\^emes arguments que ceux \'evoqu\'es  dans le {\it th\'eor\`eme 1} ([KII], p.13). On rappelle bri\`evement que si $\eta:S_{1}\rightarrow S$ est un morphisme plat d'espaces complexes r\'eduits et \par
\centerline{$\xymatrix{X_{1}\ar[d]_{\pi_{1}}\ar[r]_{\Theta}\ar[r]&X\ar[d]^{\pi}&\\
S_{1}\ar[r]_{\eta}&S&}$} \noindent
le diagramme commutatif (cart\'esien) de changement de base qui en d\'ecoule, la stabilit\'e  par changement de base des  morphismes universellement ouverts et propres  nous ram\`ene \`a \'etablir l'existence d'une fl\`eche (non triviale!)
$\Theta^{*}\omega^{n}_{X/S}\rightarrow
\omega^{n}_{X_{1}/S_{1}}$.\par\noindent Mais $\Theta$ \'etant plat, les \'egalit\'es fonctorielles ${\rm I}\!{\rm L}\Theta^{*}=
 \Theta^{*}$ et ${\rm I}\!{\rm L}\Theta^{*}{\rm I}\!
 {\rm R}{\cal H}om(A,B)= {\rm I}\!
 {\rm R}{\cal H}om({\rm I}\!{\rm L}\Theta^{*}A,{\rm I}\!{\rm
 L}\Theta^{*}B)$ permettent, en utilisant des suites spectrales ad\'equates,  d'aboutir au r\'esultat tout comme il a \'et\'e fait dans le {\it th\'eor\`eme 1}.\smallskip\noindent
{\bf(ii)}\par
$\bullet$ La construction de ce morphisme d' {\it int\'egration} est une cons\'equence de l'existence de la {\it trace relative} construite dans [R.R.V]. En effet,  comme $\pi$ est \`a fibres de
dimension born\'ee par l'entier entier $n$, les annulations des
faisceaux de cohomologie ${\cal H}^{j}(\pi^{!}{\cal O}_{S})$ pour
tout $j<-n$, nous donnent une fl\`eche naturelle (edge) ${\cal
H}^{-n}(\pi^{!}{\cal O}_{S})[n]\rightarrow \pi^{!}({\cal O}_{S})$.
Alors, en appliquant le foncteur ${\rm I}\!{\rm
 R}{\pi}_{*}$ et utilisant la trace d\'eduite de [RRV], on obtient
 le morphisme $${\rm I}\!{\rm
 R}{\pi}_{*}{\cal
H}^{-n}(\pi^{!}{\cal O}_{S})[n]\rightarrow {\cal O}_{S}$$ dont on
prend la cohomologie en degr\'e $0$ pour aboutir au morphisme ${\cal
O}_{S}$-lin\'eaire continu $ {\rm I}\!{\rm
 R}^{n}{\pi}_{*}{\cal H}^{-n}(\pi^{!}({\cal O}_{S}))\rightarrow {\cal O}_{S}$ rendant, naturellement, commutatif le diagramme
$$\xymatrix{{\rm I}\!{\rm
 R}{\pi}_{*}\omega^{n}_{X/S}[n]\ar[rr]\ar[rd]&&{\rm I}\!{\rm
 R}{\pi}_{*} \pi^{!}({\cal O}_{S})\ar[ld]\\
&{\cal O}_{S}&}$$ De plus, en vertu de ce qui pr\'ec\`ede,  cette
fl\`eche est compatible \`a tout changement de base plat sur
$S$.\smallskip
$\bullet$ La {\it compatibilit\'e avec la composition des
morphismes}  est une cons\'equence du {\it corollaire 2}. En effet, il
suffit d'utiliser les annulations des images directes issues du
lemme de Reiffen, des images inverses extraordinaires mentionn\'ees
plus haut et les identifications:
$${\rm I}\!{\rm H}om({\rm I}\!{\rm
 R}{{\pi}_{2}}_{*} \omega^{n_{1}+n}_{X_{2}/S}[n_{1}+n], {\cal O}_{S})\simeq
{\rm I}\!{\rm H}om({\rm I}\!{\rm
 R}{\Psi}_{*} \omega^{n_{1}+n}_{X_{2}/S}[n], {\pi_{1}}^{!}{\cal
 O}_{S}[-n_{1}])$$
 D'o\`u le morphisme canonique de "liaison"
$${\rm I}\!{\rm
 R}^{n}{{\Psi}_{*}} \omega^{n+n_{2}}_{X_{2}/S}\rightarrow
\omega^{n_{1}}_{X_{1}/S}$$ Il s'en suit que, si $\pi_{2}$ et
$\pi_{1}$ ont m\^eme dimension relative (i.e $\Psi$
g\'en\'eriquement fini), on dispose d'une image directe
${\pi_{2}}_{*}:{\Psi}_{*}(\omega^{n}_{X_{2}/S})\rightarrow
\omega^{n}_{X_{1}/S}$
$\,\,\,\,\,\blacksquare$\smallskip\bigskip\bigskip\noindent On peut
facilement  en d\'eduire le \Cor{4}{} Soit $\pi:X\rightarrow S$ un
 morphisme propre et plat dont les fibres sont de dimension $n$.
Alors,
\par\noindent {\bf(i)} pour tout entier $k$, $\omega^{k}_{X/S}={\cal
H}om(\Omega^{n-k}_{X/S}, \omega^{n}_{X/S})$ est un  faisceau ${\cal
O}_{X}$-coh\'erent, de profondeur au moins deux fibre par fibre sur
$S$ et repr\'esente le faisceau des $k$- formes r\'eguli\`eres sur
$X$,\par\noindent {\bf(ii)} si $X$ est r\'eduit sur $S$, ces
faisceaux sont les faisceaux des $k$-formes m\'eromorphes
r\'eguli\`eres au sens de
Kunz-Waldi-Kersken.\rm\smallskip\smallskip\noindent
{\bf 3.2.1.3. Remarques}\par\noindent
{\bf(i)} Dans le cas non propre, bien que l'on ne puisse  invoquer de  dualit\'e analytique relative,  on a tout de
m\^eme un r\'esultat positif (cf {\it th\'eor\`eme 1} de [KII]) mettant en
\'evidence un faisceau $\omega^{n}_{\pi}$ poss\'edant de nombreuses
propri\'et\'es parmi lesquelles celle de  coincider  avec
 $\omega^{n}_{X/S}$ dans le cas propre.\par\noindent
{\bf(ii)} Ces constructions nous sugg\`ere d'introduire la notion
de faisceaux dualisants relatifs analogues des dualisants de
Andr\'eotti-Kas-Golovin, en posant, pour tout faisceau coh\'erent
${\cal F}$ sur $X$,
$${\cal D}_{/S}^{n}({\cal F})= {\cal H}om({\cal F},
\omega^{n}_{X/S})$$  Ce qui a  l'avantage de ne faire intervenir
qu'une partie v\'erifiable de la dualit\'e analytique relative.\smallskip\par\noindent
Pour terminer, rappelons que le
{\it corollaire (3.1)}  donne la conditon  ``optimale'' (g\'en\'eralisant le cas plat) pour que de telles formes m\'eromorphes  r\'eguli\`eres existent. Rappelons le ici pour l'inclure de fa\c con naturelle \`a ce qui pr\'ec\`ede.
 \Cor{3.1}{} Soit $\pi:X\rightarrow S$ un
morphisme $n$-analytiquement g\'eom\'etriquement plat  d'espaces
analytiques complexes r\'eduits de dimension localement fini. Alors,
il existe un faisceau analytique (unique \`a isomorphisme canonique
pr\`es) $\Lambda^{n}_{X/S}$ sur $X$ v\'erifiant :\par\noindent
{\bf(i)} il est ${\cal O}_{X}$-coh\'erent, de profondeur au moins
deux fibre par fibre sur $S$, compatible aux changements de bases et
coincide avec le faisceau des formes holomorphes relatives sur la
partie r\'eguli\`ere de $\pi$,\par \noindent {\bf(ii)} En chaque
point $x$ de $X$ en lequel $\pi$ est plat, le germe
$\Lambda^{n}_{X/S,x}$ coincide avec le faisceau des formes
m\'eromorphes r\'eguli\`eres relatives caract\'eris\'e par la
propri\'et\'e de la trace relative.\par\noindent {\bf(iii)} De plus,
la famille de faisceaux ${\cal O}_{X}$-coh\'erents
$\Lambda^{\bullet}_{X/S}:={\cal H}om(\Omega^{n-\bullet}_{X/S},
\Lambda^{n}_{X/S})$ est munie d'une diff\'erentielle non triviale
${\rm D}$, faisant de $(\Lambda^{\bullet}_{X/S}, {\rm D})$  un
complexe diff\'erentiel de $(\Omega^{\bullet}_{X/S},
d_{X/S})$-modules.\rm\smallskip\noindent ne n\'ecessite aucun
recours \`a la th\'eorie de la dualit\'e analytique relative qui
n'est, d'ailleurs, pas encore au point dans le cas g\'en\'eral ni
\`a la construction de Kersken. \vfill\eject\noindent{\tite IV.
Classe fondamentale.}\rm\bigskip\noindent {\bf 4.0. Cas
absolu.}\smallskip\noindent {\bf 4.0.1. Cadre
alg\'ebrique.}\par\noindent Pour l'aspect alg\'ebrique et enl
caract\'eristique nulle, nous renvoyons le lecteur \`a [E] (Thm 3.1,
p34)  o\`u il est montr\'e qu'\`a  toute vari\'et\'e alg\'ebrique
$X$ admettant un  complexe r\'esiduel  ${\cal K}^{\bullet}_{X}$ de
complexe dualisant ${\cal D}^{\bullet}_{X}$ est canoniquement
associ\'e  un unique  morphisme ${\cal C}_{X}$ de ${\rm I}\!{\rm
H}om(\Omega^{n}_{X}[n],{\cal D}^{\bullet}_{X})$, qui est, en fait,
une  section de l'objet bigradu\'e $ K^{{\bullet},\star}_{X}={\cal
H}om(\Omega^{\star}_{X}, K^{\bullet}_{X})$ annul\'ee  par ses
diff\'erentielles $d_{X}$ et $\delta_{X}$ et v\'erifiant la {\it{
Propri\'et\'e de la trace}} disant que  pour tout morphisme fini et
dominant $f:U\rightarrow U'$, d'un ouvert $U$ de $X$ sur une
vari\'et\'e lisse de dimension $n$, les morphismes traces ${\cal
T}_{f}:f_{*}K^{\bullet}_{U}\rightarrow K^{\bullet}_{U'}$, ${\cal
T}_{f}:f_{*}{\cal O}_{U}\rightarrow {\cal O}_{U'}$ et le morphisme
canonique $f^{*}: f^{*}\Omega^{n}_{U'}\rightarrow \Omega^{n}_{U}$
induisent le diagramme commutatif
$$\xymatrix{f_{*}f^{*}\Omega^{n}_{U'}[n]\ar[r]^{f^{*}}\ar[d]_{{\cal T}_{f}\otimes Id[n]}&f_{*}\Omega^{n}_{U}[n]\ar[r]^{{\cal C}_{U}}&f_{*}K^{\bullet}_{U}\ar[d]^{{\cal T}_{f}}\\
\Omega^{n}_{U'}[n]\ar[rr]_{{\cal
C}_{U'}}&&K^{\bullet}_{U'}}$$\smallskip\noindent {\bf 4.0.2. Cadre
analytique complexe.}\smallskip\noindent {\bf 4.0.2.0. Du global au
local.} Disposant du  complexe dualisant ${\cal D}^{\bullet}_{X}$ de
Ramis-Ruget [RR], il est facile de voir que la construction de
Elzein s'adapte au  cadre analytique complexe (\`a quelques
modifications mineures pr\`es!) et permettent de dire  qu'\`a  tout
$n$-cycle $X$ d'un espace analytique complexe $Z$, est canoniquement
associ\'e ( dans la cat\'egorie d\'eriv\'ee des ${\cal
O}_{X}$-modules)  un unique  morphisme ${\cal C}_{X}:
\Omega^{n}_{X}[n]\rightarrow {\cal D}^{\bullet}_{X}$ compatible aux
inclusions ouvertes sur $X$ (i.e de nature locale sur $X$),
compatible avec la diff\'erentielle ext\'erieure usuelle (d\'efinie
au sens des distributions) et  v\'erifiant {\it{la propri\'et\'e de
la trace}}, c'est-\`a-dire que pour toute param\'etrisation locale
$f:X\rightarrow Y$, $Y$ \'etant un polydisque ouvert de ${\Bbb
C}^{n}$, on a un diagramme {\bf{commutatif}}
$$\xymatrix{f_{*}f^{*}\Omega^{n}_{Y}[n]
\simeq f_{*}({\cal O}_{X}\otimes
f^{*}\Omega^{n}_{Y}[n])\ar[d]_{{\cal T}^{0}_{f}\otimes Id[n]}
\ar[r]&\!\!\!f_{*}\Omega^{n}_{X}[n]
\ar[r]^{\indent{f_{*}{\cal C}_{X}}}&f_{*}{\cal D}^{\bullet}_{X}\ar[d]^{{\cal T}^{\bullet}_{f}}\\
\Omega^{n}_{Y}[n]\ar[rr]_{Id}&&\Omega^{n}_{Y}[n]}$$\par\noindent
 dans
lequel ${\cal T}^{\bullet}_{f}$ est le morphisme trace d\'efini dans
[R.R].
\par\noindent
Par construction, les faisceaux de cohomologie ${\cal H}^{j}({\cal
D}^{\bullet}_{X})$ sont coh\'erents et nuls pour  $j<-n$. Il en
r\'esulte que le morphisme ${\cal C}_{X}$ est en fait enti\`erement
caract\'eris\'e par le morphisme de faisceaux ${\cal
O}_{X}$-coh\'erents $\Omega^{n}_{X}\rightarrow {\cal H}^{-n}({\cal
D}^{\bullet}_{X})$ et que l'on note encore (abusivement) ${\cal
C}_{X}$.  Supposons $X$ irr\'eductible,  localement plong\'e dans une vari\'et\'e de Stein $Z$ de dimension $N$ et posons $p:=N-n$. Alors,  la succesion de morphismes naturels (dans la cat\'egorie
d\'eriv\'ee)
$${\rm I}\!{\rm H}om( \Omega^{n}_{X}[n],{\cal D}^{\bullet}_{X})\simeq {\rm I}\!{\rm H}om( \sigma_{*}\Omega^{n}_{X}[n],{\cal D}^{\bullet}_{Z})\rightarrow
{\rm I}\!{\rm H}om(\sigma_{*}{\cal
O}_{X}\otimes\Omega^{n}_{Z}[n],\Omega^{n+p}_{Z}[n+p])$$
$$\simeq {\rm E}xt^{p}({\cal O}_{X}, \Omega^{p}_{Z})\simeq
\Gamma(Z, {\cal E}xt^{p}({\cal O}_{X}, \Omega^{p}_{Z}))$$
obtenus  gr\^ace \`a la dualit\'e pour un plongement, la lissit\'e
de $Z$ et les  annulations de cohomologie du complexe dualisant,
permet de voir le morphisme classe fondamentale comme une section
globale du faisceau coh\'erent  ${\cal E}xt^{p}({\cal
O}_{X}, \Omega^{p}_{Z})$.  Il est d'usage d'appeler {\it classe
fondamentale } de  $X$ dans $Z$, l'image de ce dernier \'el\'ement
dans ${\rm H}^{p}_{|X|}(Z, \Omega^{p}_{Z})$ via le morphisme
canonique  ${\rm E}xt^{p}({\cal O}_{X},
\Omega^{p}_{Z})\rightarrow{\rm H}^{p}_{|X|}(Z, \Omega^{p}_{Z})$. La
classe obtenue est \'evidemment annul\'ee par tout  id\'eal de
d\'efinition de $X$ dans $Z$, par la diff\'erentielle ext\'erieure
usuelle $d$ et induit un morphisme de ${\cal O}_{X}$-modules gradu\'es $\Omega^{\bullet}_{X}\rightarrow{\cal E}xt^{p}({\cal O}_{X},
\Omega^{p+\bullet}_{Z}) $. \smallskip\noindent
{\bf 4.0.2.1. Du local au global.}  On peut adopter un point de vue ne n\'ecessitant aucun recours au complexe dualisant. Pour cela, la strat\'egie consiste  \`a passer du local au global par recollement de donn\'ees simplicielles. \par\noindent
 On commence donc par associer, \`a tout sous ensemble analytique  $X$ de codimension pure  $p$ dans une vari\'et\'e lisse de Stein  $Z$ de dimension $n+p$, une classe de cohomologie dans ${\rm H}^{p}_{X}(Z, \Omega^{p}_{Z})$ v\'erifiant un certain nombre de propri\'et\'es fonctorielles. On se ram\`ene fondamentalement \`a l'\'etude du cas classique  o\`u $X$ est un point de ${\Bbb C}^{p}$  que l'on trouve ais\'ement dans la litt\'erature (cf [Gri-H]). Les calculs explicite, en terme de r\'esidu de symboles, montrent que l'\'el\'ement d\'esir\'e est en fait une section globale  du faisceau coh\'erent ${\cal E}xt^{p}({\cal O}_{X},
\Omega^{p}_{Z})$ (s'explicitant au moyen des r\'esolutions de Koszul ou de Dolbeault-
Grothendieck) ayant de bonnes propri\'et\'es de stabilit\'e
vis-\`a-vis des morphismes finis. $Z$ \'etant lisse, ces classes
d\'efinissent naturellement
 un morphisme $\tilde{{\cal C}}^{\sigma}_{X}: \Omega^{n}_{X}\rightarrow {\cal
E}xt^{p}({\cal O}_{X}, \Omega^{n+p}_{Z})$.\par\noindent Si
$\ds{X:=\sum_{j}n_{j}X_{j}}$ est un $n$-cycle quelconque d'un
espaces analytique complexe $Z$ arbitraire, on proc\`ede par
localisation en nous ramenant \`a la situation universelle (cf {\bf
(1.0.3.1.1)}) pour laquelle on  montre que la diagonale
g\'en\'eralis\'ee $\sigma_{k}$-invariante et donc le sous espace
d'incidence $\#$ admet une classe fondamentale, au sens o\`u on
l'entend, repr\'esent\'e simplement par un objet de la cohomologie \'`equivariante. Alors, par image r\'eciproque et utilisant les
th\'eor\`emes d'annulations de Siu-Trautmann [S.T], on produit une
classe de ${\rm H}^{p}_{|X|}(U\times B, \Omega^{p}_{U\times B})$ et
m\^eme de ${\rm  E}xt^{p}({\cal O}_{|X|}, \Omega^{n+p}_{U\times
B})$. Ainsi, pour toute installation de $X$ ou \'ecaille adapt\'ee
$(V,\sigma, U, B)$, surgit un morphisme $\tilde{{\cal
C}}^{\sigma}_{X}:\Omega^{n}_{|X|}\rightarrow {\cal E}xt^{p}({\cal
O}_{|X|}, \Omega^{n+p}_{U\times B})$. Mais cette collection de
morphismes locaux se recollent globalement et d\'efinit un morphisme
canonique  $\tilde{{\cal C}}_{X}:\Omega^{n}_{|X|}\rightarrow
\omega^{n}_{|X|}$ puisque la collection de faisceaux locaux ${\cal
E}xt^{p}({\cal O}_{|X|}, \Omega^{n+p}_{U\times B})$ se recollent
globalement en le faisceau dualisant de Grothendieck.
\smallskip\noindent
 La propri\'et\'e
de la trace est v\'erifi\'ee puisque pour toute param\'etrisation locale $f:|X|\rightarrow U$ et tout plongement local  $\sigma$ de $|X|$ dans $Z$ lisse de dimension relative $p$ sur $U$, on dispose de fl\`eches
d'int\'egration ou r\'esidues:\par \centerline{
$f_{*}({\cal E}xt^{p}(\sigma_{*}{\cal O}_{X},
\Omega^{n+p}_{Z}))\rightarrow f_{*}{\cal
H}^{p}_{|X|}(\Omega^{n+p}_{Z})\rightarrow \Omega^{n}_{Y}$}\noindent
D\`es lors, il est facile d'en d\'eduire que $\tilde{{\cal C}}_{X} = {\cal
C}_{X}$. \smallskip\noindent Rappelons que dans le cas d'une
intersection compl\`ete, on a une description simple de cette classe
 en terme de cochaines de \v Cech. Supposons $Z={\Bbb C}^{n+p}$ et
$X=\{z\in{\Bbb C}^{n+p}/f_{1}(z)=f_{2}(z)=\cdots=f_{p}(z)=0\}$ de
dimension pure $n$, les $(f_{i})$ \'etant des fonctions holomorphes
et soit ${\bf{\cal U}}$ le recouvrement ouvert du compl\'ementaire
de X dans ${\Bbb C}^{n+p}$, dont les \'el\'ements sont les ouverts
$U_{i}=\{z\in{\Bbb C}^{n+p}/f_{i}(z)\not= 0\}$. Alors un
repr\'esentant de \v{C}ech-Leray  de $C(X)$ est donn\'e par le
cocycle\par \centerline{$\ds{{df_{1}\wedge df_{2}\wedge\cdots\wedge
df_{p}}\over{f_{1}\cdots\cdots f_{p}}}$}
\smallskip\par\noindent
{\bf 4.0.2.2.  Classe fondamentale et
int\'egration}\smallskip\noindent Rappelons que le  faisceau
$\omega^{n}_{X}$ est  coh\'erent et dualisant pour le faisceau
structural ${\cal O}_{X}$ selon Andr\'eotti- Kas- [A.K]. Alors,
utilisant l'int\'egration usuelle $\ds{{\rm H}^{n}_{c}(X,
\Omega^{n}_{X})\rightarrow {\Bbb C}}$ et le morphisme classe
fondamentale $\ds{\Omega^{n}_{X}\rightarrow \omega^{n}_{X}}$, on
obtient le diagramme commutatif
$$\xymatrix{{\rm H}^{n}_{c}(X,
\Omega^{n}_{X})\ar[rr]\ar[rd]&&{\rm H}^{n}_{c}(X, \omega^{n}_{X})\ar[ld]\\
&{\Bbb C}&}$$ Dans lequel la fl\`eche ${\rm H}^{n}_{c}(X,
\omega^{n}_{X})\rightarrow {\Bbb C}$ n'est rien d'autre que le
morphisme {\it{trace absolu}} d\'efini sur le complexe dualisant de
Ramis et Ruget (cf {\bf (3.1.3.1.1)}).
 On obtient, alors,  un diagramme commutatif
$$\xymatrix{{\rm I}\!{\rm H}om(\omega^{n}_{X},
\omega^{n}_{X})\ar[r]\ar[d]&{\rm I}\!{\rm H}om_{cont}({\rm
H}^{n}_{c}(X,\omega^{n}_{X}), {\Bbb C})\ar[d]\\
{\rm I}\!{\rm H}om(\Omega^{n}_{X}, \omega^{n}_{X})\ar[r]&{\rm
I}\!{\rm H}om_{cont}({\rm H}^{n}_{c}(X,\Omega^{n}_{X}), {\Bbb C})}$$
dans lequel la notation ${\rm I}\!{\rm H}om_{cont}$ d\'esigne les morphismes  ${\Bbb C}$-lin\'eaire et  continues. Ce diagramme
met  en correspondance l'int\'egration et le morphisme classe
fondamentale qui peut \^etre vu comme l'image de la fonction $1$ via
le morphisme injectif
$${\cal H}om(\omega^{n}_{X}, \omega^{n}_{X})\rightarrow
{\cal H}om(\Omega^{n}_{X},\omega^{n}_{X}),$$ On peut pr\'eciser que
${\cal H}om(\omega^{n}_{X},\omega^{n}_{X})$ est un sous faisceau de
${\cal H}om(\Omega^{n}_{X},\omega^{n}_{X})$ qui est isomorphe au
faisceau $\omega^{0}_{X}$ contenant ${\cal O}_{X}$ et coincidant
avec ce dernier si $X$ est normal. Si $X$ est r\'eduit, $\omega^{0}_{X}$ s'identifie au faisceau des fonctions m\'eromorphes localement born\'ees.
\smallskip\noindent
Signalons qu'il est facile de  reconstituer le morphisme
d'int\'egration globale \`a partir des morceaux locaux. En effet, il
suffit pour cela de choisir un recouvrement ouvert localement fini
 $(V_{\alpha})_{\alpha\in A}$de $X$  tel que chaque $V_{\alpha}$ soit muni d'un plongement dans un
polydisque $Z_{\alpha}$ de ${\Bbb C}^{n+p}$ et d'une projection sur
un polydisque $U_{\alpha}$ de ${\Bbb C}^{n}$ faisant de $V_{\alpha}$
un rev\^etement ramifi\'e de $U_{\alpha}$. On supposera $X$
paracompact et compl\`etement paracompact pour pouvoir localiser et
globaliser au moyen des partitions de l'unit\'e les calculs
explicitables en terme de cochaines de \v Cech. Signalons que
l'exactitude \`a droite du foncteur ${\rm H}^{n}_{c}(X,-)$ permet de
localiser sur $X$. \smallskip\par\noindent {\bf 4.1.  Classe
fondamentale relative. }\smallskip\noindent {\bf 4.1.0.  Cadre
alg\'ebrique.}\par\noindent Si ${\bf k}$ est un corps de
caract\'eristique nulle et  $\pi$ un morphisme  $n$-
\'equidimensionnel de ${\bf k}$-sch\'emas noeth\'eriens munis de
complexes r\'esiduels. Alors, Elzein et Ang\'eniol ([E],[A.E]) ont
montr\'e l'existence d' un unique (et donc canonique) \'el\'ement
dans ${\rm I}\!{\rm H}om(\Omega^{n}_{X/S}[n], {\cal
 D}^{\bullet}_{X/S})$ v\'erifiant {\it{la propri\'et\'e de la
 trace relative}} dans les cas suivants:\par
$\bullet$ $X$ intersection compl\`ete relative sur $S$ quelconque ([E]),\par
$\bullet$ $X$ est plat sur $S$ r\'eduit ([E]),\par
$\bullet$ $X$ \'equidimensionnel sur $S$ normal ([E]),\par
$\bullet$ $X$ de {\it Tor-dimension finie} sur $S$ r\'eduit ([E.A]]).\par\noindent
Cette  construction est  compatible aux changements de bases quelconques, r\'eduits, normaux et cohomologiquement transversaux \`a $\pi$ respectivement. \smallskip\noindent
Par analogie avec le cas absolu et eu \'egard \`a l'\'equidimensionnalit\'e,
la {\it propri\'et\'e de la trace relative} exprime le fait que pour
toute installation
  $\xymatrix{X\ar[r]^{f}\ar@/_/[rr]_{\pi}&Y\ar[r]^{q}&S}$, comme ci-dessus, on a un  diagramme commutatif
$$\xymatrix{f^{*}f_{*}\Omega^{n}_{Y/S}[n]\simeq f_{*}({\cal O}_{X})
\otimes_{{\cal O}_{Y}}\Omega^{n}_{Y/S}[n]\ar[d]_{{\cal T}^{0}_{f}
\otimes Id[n]}\ar[r]&f_{*}({\cal O}_{X}\otimes \Omega^{n}_{Z/S}[n])
\ar[r]^{\,\,\,\,\,f_{*}{\cal C}_{X}}&f_{*}{\cal D}^{\bullet}_{X/S}\ar[d]^{{\cal T}_{f}}\\
\Omega^{n}_{Y/S}[n]:={\cal D}^{\bullet}_{Y/S}\ar[rr]_{Id}&&\Omega^{n}_{Y/S}[n]}$$ dans lequel
le morphisme ${\cal T}_{f}$ peut \^etre  d\'efini pour tout
morphisme propre et surjectif alors que ${\cal T}^{0}_{f}:
f_{*}{\cal O}_{X}\rightarrow {\cal O}_{Y}$ n'existe pas en
g\'en\'eral. \smallskip\noindent
{\bf 4.1.1.  Cadre analytique complexe.}\par\noindent
En g\'eom\'etrie complexe, la th\'eorie de la  dualit\'e analytique relative est encore tr\`es loin d'\^etre au point hormis le cas d'un morphisme propre ou  lisse en plus de  quelques cas sporadiques tr\`es particuliers. L'obstacle majeur est qu'actuellement on est dans l'incapacit\'e d'\'etendre  le `` foncteur''  de la g\'eom\'etrie alg\'ebrique $\pi^{!}$, avec toutes les propri\'et\'es qui le caract\'erisent, \`a la cat\'egorie des espaces complexes (hormis les cas cit\'es). Nous nous
attarderons plus longuement sur ces questions  dans [KIII].\par\noindent
N\'eanmoins,  si $\pi$ est donn\'e avec un plongement de $X$ dans $Z$ lisse et  de dimension relative $n+p$ sur $S$ r\'eduit,
 D.Barlet  montre dans [B4] que la classe fondamentale
 relative de $X$ sur $S$ existe si et seulement si $X$ est le
 support du graphe d'une famille analytique de $n$ -cycles
$(X_{s})_{s\in S}$ de $Z$. De plus, cette construction est
compatible aux changement de bases quelconques entre espaces
complexes r\'eduits quelconques. On peut proposer l' \'enonc\'e suivant, en apparence,   sensiblement plus g\'en\'eral: \Th{0}{}([B4]) Soient $S$, $Z$
et $X$ des espaces analytiques complexes r\'eduits install\'es dans
le diagramme commutatif
$\xymatrix{X\ar[r]^{\sigma}\ar@/_/[rr]_{\pi}&Z\ar[r]^{q}&S}$ dans
lequel $q$ est lisse de dimension relative $n+p$, $\sigma$ est un
plongement local et $\pi$ un morphisme contin\^ument
$n$-g\'eom\'etriquement plat (donc  universellement $n$-\'equidimensionnel)  muni d'une certaine pond\'eration
$\goth{X}$.
\par\noindent Alors $\pi$ est  analytiquement
$n$-g\'eom\'etriquement plat si et seulement si il existe  une classe de cohomologie  ${\rm C}_{\pi,\goth{X}}$ de ${\rm H}^{p}_{X}(Z, \Omega^{p}_{Z/S})$ v\'erifiant les propri\'et\'es suivantes:\par\noindent
elle est de nature locale sur $X$ et $S$, de  formation compatible \`a l'additivit\'e des pond\'erations et \`a tout changement de base entre espaces complexes r\'eduits donnant, en particulier, pour chaque $s$ fix\'e la classe fondamentale absolue du cycle $[\pi^{-1}(s)]$. De plus, si $\goth X$ est la pond\'eration standard, ${\rm C}_{\pi,\goth{X}}={\cal C}_{X/S}$ construite dans [B4].
\rm\smallskip\noindent
\dem Le lecteur d\'esireux d'approfondir peut consulter [B4]. Nous pr\'esentons ici  les arguments fondamentaux. Nous commencerons par quelques remarques d'ordre g\'en\'erales.\par
La version relative des th\'eor\`emes
d'annulations de cohomologie de Siu- Trautmann ([S.T]) donn\'ee dans [K1], {\it proposition}  p.294 et assurant l'isomorphisme\par
\centerline{ $\ds{{\rm H}^{p}_{X}( Z,\Omega^{p}_{Z/S})\simeq {\rm
H}^{0}(Z,{\cal H}^{p}_{X}(\Omega^{p}_{Z/S})}$}\par\noindent
permet de nous ramener \`a un probl\`eme local sur $Z$ au voisinage de $X$ et de supposer $S$ irr\'eductible. En effet, en raisonnant par r\'ecurrence sur la dimension de $S$ et utilisant la variante relative de Siu-Trautmann, il est facile de voir que cette classe, si elle existe, est unique et enti\`erement d\'etermin\'ee par les points g\'en\'eriques des composantes irr\'eductibles de $S$. Par ailleurs,  pour tout recouvrement ouvert
localement fini $(X_{\alpha})_{\alpha\in a}$ de $X$ que l'on peut choisir adapt\'e \`a  $\pi$ (qui est
universellement $n$-\'equidimensionnel) dans le sens o\`u  chacun de ses
\'el\'ements est  muni d'installation locale du type $(\spadesuit)_{\alpha\in A}$
$$\xymatrix{X_{\alpha}\ar[rdd]_{\pi_{\alpha}}\ar[rr]^{\sigma_{\alpha}}\ar[rd]^{f_{\alpha}}&&Z_{\alpha}\ar[ld]_{q'_{\alpha}}\ar[ldd]^{q_{\alpha}}\\
&Y_{\alpha}\ar[d]^{q''_{\alpha}}&\\&S&}$$ o\`u $\sigma_{\alpha}$ est
un plongement local, $\pi_{\alpha}$ la restriction de $\pi$ \`a
$X_{\alpha}$, $Y_{\alpha}$ lisse sur $S$ et de dimension  relative $n$,  $Z_{\alpha}$ des espaces complexes
lisses sur $S$ de dimension relative $n+p$ ou ouvert de carte que l'on peut  toujours supposer $S$-adapt\'e dans le sens de {\bf(1.0.2.2)}, $f_{\alpha}$ fini, ouvert  et surjectif, $q_{\alpha}$,$q'_{\alpha}$
et $q''_{\alpha}$ lisses, le morphisme canonique
  \par\noindent \centerline{$\ds{{\rm
H}^{p}_{X}(Z,\Omega^{p}_{Z/S})\rightarrow \prod_{\alpha\in A}{\rm
H}^{p}_{X_{\alpha}}(Z_{\alpha},\Omega^{p}_{
Z_{\alpha}/S})}$}\par\noindent
est injectif. Ainsi, il nous suffit de construire un tel objet seulement  localement en v\'erifiant ses principales propri\'et\'es dans cette situation.\smallskip
$\blacklozenge$ $\Rightarrow$:\par\noindent
Comme  $\pi$ est analytiquement g\'eom\'etriquement plat, il en est de m\^eme pour chaque  $\pi_{\alpha}$ et
$f_{\alpha}$ de  la ``d\'ecomposition'' locale pr\'ec\'edente. On peut, d\`es lors, supposer  $Z:=S\times U\times B$, $Y:=S\times U$,  $\pi:=\pi_{\alpha}$ et $S$ irr\'eductible. Ayant localiser ainsi $\pi$  autour de l'une quelconque de ses fibres tel qu'en {\bf(1.0.3.2.1)}, on obtient  diagramme commutatif du type
$$\xymatrix{&S\times U\times B\ar[r]^{\tilde{F}:=F\times Id_{B}}\ar[dd]&{\rm Sym}^{k}(B)
 \times B\ar[dd]&\\
X\ar[ur]\ar[rd]_{f}\ar[rrr]&&& {\rm Sym}^{k}(B) \# B\ar[ld]\ar[ul]
\\&S\times U\ar[r]_{F}&{\rm Sym}^{k}(B)&}$$
Cela a pour effet imm\'ediat de nous ramener  \`a la situation universelle et \`a de la
cohomologie \'equivariante pour laquelle on construit  la classe
 fondamentale locale du sous espace d'incidence $\#$  dans
 $\rm{Sym}^{k}(B)\times B$  qui est  un \'el\'ement
  $C_{\#}$ de ${\rm H}^{p}_{|\#|}(\rm{Sym}^{k}( B)\times B,
  \omega^{p}_{\rm{Sym}^{k}(B)\times{\Bbb C}^{p} })$
  coincidant g\'en\'eriquement (sur les points r\'eguliers de
  $\#$) avec la classe fondamentale de la diagonale
  g\'en\'eralis\'ee de $(B)^{k}\times B$. Rappelons que
 le faisceau $\omega^{p}_{{\rm{Sym}}^{k}({\Bbb C}^{p})}$ que Barlet a appel\'e le faisceau des formes ``truqu\'ees'' n'est rien d'autre que l'image directe par l'application
   quotient $q: ({\Bbb C}^{p})^{k}\rightarrow{\rm{Sym}}^{k}({\Bbb C}^{p})$ du faisceau des formes holomorphes $\Omega^{p}_{{\Bbb C}^{kp}}$ dont on prend la partie $\sigma_{k}$-invariante, en d'autre termes:
$$\omega^{p}_{{\rm{Sym}}^{k}({\Bbb C}^{p})}=(q_{*}(\Omega^{p}_{{\Bbb C}^{kp}}))^{\sigma_{k}}$$
\par\noindent
Or, comme on est dans le cas d'une V-vari\'et\'e, ce faisceau
coincide avec le faisceau ${\cal L}^{p}$ des $p$ formes
m\'eromorphes se prolongeant
 analytiquement sur toute d\'esingularis\'ee de ladite V-vari\'et\'e. Mais ce
  faisceau a la propri\'et\'e d'\^etre stable par image r\'eciproque quelconque
   entre espaces complexes, d'o\`u morphisme
   $\tilde{F}^{*}\omega^{p}_{\rm{Sym}^{k}(B)\times B}\rightarrow {\cal L}^{p}_{S\times U\times B}$. Alors, utilisant la projection naturelle  ${\cal L}^{p}_{S\times U\times B}\rightarrow {\cal L}^{0}_{S}\otimes\Omega^{p}_{U\times B}$ et l'analyticit\'e g\'eom\'etrique de la pond\'eration assurant que la variation en le param\`etre est holomorphe, on en d\'eduit une image r\'eciproque plus pr\'ecise \`a savoir  $\tilde{F}^{*}\omega^{p}_{\rm{Sym}^{k}(B)\times B}\rightarrow  \Omega^{p}_{U\times B/S}$.\par\noindent
Ainsi, $C_{\Delta}$ fournira
 par image r\'eciproque par $\tilde{F}$ un \'el\'ement de $H^{p}_{X}(S\times U\times B, \Omega^{p}_{/S})$
  qui coincide aux points g\'en\'eriques de $X$ (et donc aux point r\'eguliers de $\pi$) avec la classe fondamentale usuelle
   d'une sous vari\'et\'e lisse.\smallskip\noindent
On aura
alors, par recollement, un objet global ${\rm C}_{\pi,\goth{X}}$ de
${\rm H}^{p}_{X}(Z, \Omega^{p}_{Z/S})$ dont on v\'erifie
ais\'ement les propri\'et\'es fondamentales.\smallskip
  $\blacklozenge$ $\Leftarrow$: Supposons donn\'ee
   une telle classe ${\rm C}_{\pi,\goth{X}}$ pour $\pi$ et un recouvrement
   ouvert $(X_{\alpha})_{\alpha\in a}$ de $X$ comme pr\'ec\'edemment d\'efini.
   Il nous faut montrer  que cette classe
induit, pour toute projection finie
$f_{\alpha}:X_{\alpha}\rightarrow S\times U_{\alpha}$, un morphisme
trace $\ds{{\cal T}^{0}_{f_{\alpha}}:{f_{\alpha}}_{*}{\cal
O}_{X_{\alpha}}\rightarrow {\cal O}_{S\times U_{\alpha}}}$ se
prolongeant en un morphisme d'alg\`ebres gradu\'ees
\par\noindent
$\ds{{\cal
T}^{\bullet}_{f_{\alpha}}:{f_{\alpha}}_{*}\Omega^{\bullet}_{X_{\alpha}/S}\rightarrow
\Omega^{\bullet}_{S\times U_{\alpha}/S}}$. \par\noindent En
ommettant les indices et en ins\'erant ces morphismes finis et
surjectifs dans le diagramme
$\xymatrix{X\ar[r]^{\sigma}\ar@/_/[rr]_{f}&Z\ar[r]^{q'}&Y}$, la
construction de ces traces repose essentiellement  sur une formule
int\'egrale de type Cauchy. En efet, on utilise une forme
particuli\`ere de l'int\'egration sur les cycles consistant \`a
int\'egrer  des classes de cohomologie de type $(p,p)$ \`a support
dans $X$ (elle peut aussi \^etre vue  comme un morphisme {\bf
r\'esidu}). Plus pr\'ecisemment, on a le diagramme commutatif
$$\xymatrix{f_{*}{\cal O}_{X}\ar[rd]^{{\cal T}_{f}}\ar[r]^{\!\!\!\!\!\!\!
{\rm C}_{\pi,\goth{X}}}&f_{*}{\cal H}^{p}_{X}
(\Omega^{p}_{Z/S})\ar[d]^{\phi}\\
&{\cal O}_{Y}}$$ dans lequel ${\cal T}^{0}_{f}$ est le morphisme
trace recherch\'e, la fl\`eche horizontale est le cup-produit par la
classe fondamentale ${\rm C}_{\pi,\goth{X}}$ et $\phi$ le morphisme
d'int\'egration.\par\noindent Il faut simplement int\'egrer sur une
famille dont les supports rencontre $X$ en des ferm\'es $S$-propres
de $Z$. Pour cela, on se ram\`ene \`a int\'egrer sur  la famille
analytique triviale $({s}\times {t}\times B_{i})_{(s,t)\in
S_{i}\times U_{i}}$; ce qui nous  donne une fonction holomorphe sur
$S_{i}\times U_{i}$. On pose, alors, pour  toute fonction holomorphe
$f$ sur $X$,
$${\cal T}_{E,S}(f)=\int_{Y_{s}}f.{\cal C}(X_{E}/S)\,\,\,\,\,\,\,\,\,\,\,\blacksquare $$\smallskip\par\noindent
{\tite V. Quelques petits r\'esultats  utiles.}\bigskip\noindent
 \Lemme{5.0}{}
([E],p.75)\par\noindent Soient ${\cal A}$ un anneau et ${\cal B}$
une alg\`ebre finie et de Tor-dimension finie sur ${\cal A}$. Il
existe un morphisme trace de ${\cal A}$-modules $Tr: {\cal
B}\rightarrow {\cal A}$, d\'efini comme suit: si ${\cal
P}^{\bullet}$ est une r\'esolution projective finie de ${\cal B}$
sur ${\cal A}$, on fait correspondre \`a tout \'el\'ement  $b$ de
${\cal B}$ un morphisme $b_{\bullet}:{\cal P}^{\bullet}\rightarrow
{\cal P}^{\bullet}$ induisant la multiplication par $b$ dans ${\cal
B}$ et on pose  $Tr(b) = \sum_{i}(-1)^{i}Tr(b_{i})$. Cette trace
commute aux changements de base cohomologiquement transversaux au
morphisme ${\cal A}\rightarrow {\cal B}$, elle est stable par
composition de morphismes de Tor.dimension finie c'est-\`a-dire que
pour tout morphisme ${\cal B}\rightarrow {\cal C}$ de Tor.dimension
finie, on a $Tr_{{\cal C}/ {\cal A}} = Tr_{{\cal B}/ {\cal A}}o
Tr_{{\cal C}/ {\cal B}}$, et elle coincide avec la trace usuelle
dans le cas plat.\rm\smallskip\noindent \Lemme{5.1}{} ([E], \S{\bf 4.4},
lemme 1, p 60) Soient $Z$, $X$, $Y$ et $S$ des vari\'et\'es munies
d'un diagramme commutatif $$\xymatrix{X\ar[rdd]_{\pi}\ar[rr]^{\sigma}\ar[rd]^{f}&&Z\ar[ld]_{q'}\ar[ldd]^{q}\\
&Y\ar[d]^{q''}&\\
&S&}$$ avec $q$, $q'$, $q"$ lisses de dimension relatives
respectives $n+p$, $p$ et $n$, $\sigma$ plongement ferm\'e et $f$
quasi-fini et s\'eparable en un point g\'en\'erique $x$ d'une fibre
$X_{s}$ (autour de laquelle a lieu cette factorisation). Soit
$(t_{1},\cdots,t_{n})$ des \'el\'ements de ${\cal O}_{Z,x}$
fournissant une suite r\'eguli\`eres de param\`etres dans ${\cal
O}_{X_{s},x}$. Alors, \par\noindent {\bf(i)}
$\ds{\Omega^{n}_{X/S,x}=\sum_{J}{\cal O}_{X,x}dt^{J}\wedge
\Omega^{n}_{Z/S,x}}$ pour tout ensemble d'indice $J$ de longeure
$j\leq n$
\par\noindent {\bf(ii)} Il existe une famille finie
$f_{\alpha}:X\rightarrow Y$ de morphismes quasi-finis et
s\'eparables en $x$ telle que l'on ait
$\ds{\Omega^{n}_{X/S,x}=\sum_{\alpha}(f^{*}_{\alpha}\Omega^{n}_{Z/S})_{x}}$
\Lemme{5.2}{} ([E], lemme III.1.2, p83){\it{ Passage des projections
lin\'eaires aux projections quelconques.}}\par\noindent Soient
$(l_{i})_{i\in I}$ une famille finie de $S$-morphismes lin\'eaires
dont les membres sont des composantes de $S$-morphismes lin\'eaires
$Z\rightarrow Y$ finis sur $X$ au dessus de $Y$, $k_{i}$ des entiers
et $\alpha^{j}_{i}\in\Gamma(S,{\cal O}_{S})$ pour $j\in
\{1,\cdots,n+p\}$ presque partout nulles. Alors le morphisme
$\Psi:Z\rightarrow Y$ de composantes
$\psi_{j}:=\sum_{i}{\alpha^{j}_{i}}l_{i}^{k_{i}}$ est fini sur $Y$
et induit  un morphisme trace sur $Y$. \rm\smallskip\noindent
\Lemme{5.3}{} Soit $\pi_{0}: {\Bbb C}^{n+p}\rightarrow {\Bbb C}^{n}$
la projection canonique. Pour $u\in L({\Bbb C}^{n+p},{\Bbb C}^{n})$,
on pose $\pi_{u} := \pi_{0} + u$. Alors, pour tout couple d'entiers
proposition $(a,b)$ v\'erifiant $a\leq n$ et $b\leq n$ et tout
voisinage $\Omega_{0}$ de $0$ dans $L({\Bbb C}^{n+p},{\Bbb C}^{n})$,
on a:
$$\Lambda^{a,b}({\Bbb C}^{n+p}) = \sum_{u\in \Omega_{0}}\pi_{u}^{*}\lbrack\Lambda^{a,b}
({\Bbb C}^{n})^{*}\rbrack$$ \noindent o\`u $\Lambda^{a,b}(E)^{*}$
d\'esigne, pour un espace vectoriel $E$
 sur ${\Bbb C}$, l'espace des formes $a$- lin\'eaires et $b$- antilin\'eaires
  altern\'ees sur $E$.\rm\smallskip\noindent
  Ce lemme montre que si $X$ est un sous ensemble analytique de dimension pure $n$  dans ${\Bbb
  C}^{n+p}$, sur lequel $\pi_{0}$ induise un rev\^etement ramifi\'e
  d'un certain degr\'e $k$, alors il existe un ouvert dense $V_{0}$  de la
  grassmanienne ${\cal G}(n+p,n)$ tel que, pour tout $u\in V_{0}$, $\pi_{u}$ soit encore un
  rev\^etement ramifi\'e de degr\'e $k_{u}$.\bigskip\noindent
\smallskip\bigskip\noindent
\centerline{\bf R\'ef\'erences bibliographiques}

\medskip
{\baselineskip=9pt\eightrm
\item{[Ang]}Ang\'eniol~B., Familles de cycles alg\'ebriques---Sch\'emas de Chow. Lecture Notes in Mathematics.~896, Springer Verlag, (1981).
\smallskip

\item{[A.E]}Ang\'eniol~B, Elzein~F., La classe fondamentale relative d'un cycle,Bull.Soc.Math.France Mem~58 (1978) 63---93.
\smallskip

\item {[A.N]}Andreotti~A., Norguet~F., La
convexit\'e holomorphe dans l'espace analytique des cycles d'une
vari\'et\'e  alg\'ebrique, Ann. Sc. Norm. Pisa, tome  21, (1965),
811---842.
\smallskip

\item {[A.K]}Andreotti~A., Kaas~A. Duality on complex spaces, Ann. Sc. Norm. Pisa, tome  27, (1973) 187---263.
\smallskip

\item{[A.S]} Andreotti~A.,Stoll~W., Analytic and algebraic dependence of meromorphic functions. Lecture Notes in Mathematics, ~234. Springer-Verlag, Berlin-New York, (1971).
\smallskip

\item{[B1]}Barlet~D.,
Espace analytique r\'eduit des cycles analytiques complexes compacts
d'un espace analytique
 r\'eduit,  S\'em. F. Norguet-Lecture Notes in Mathematics.~482,
 Springer Verlag, (1975)  1---158.
\smallskip

\item{[B2]}Barlet~D.Faisceau $\omega^{.}_{X}$ sur un espace analytique de dimension pure,  S\'em. F. Norguet-Lecture Notes in Mathematics ~670, Springer Verlag, (1978) 187---204.
\smallskip

\item{[B3]}Barlet~D. Convexit\'e au voisinage d'un cycle. Sem.Fran\c cois Norguet-Lecture Notes in Mathematics ~807, Springer-Verlag, (1977-1979),p.102---121.
\smallskip

\item{[B4]}Barlet~D., Famille analytique de cycles et classe
fondamentale relative, S\'em. F. Norguet-Lecture Notes in
Mathematics ~807, Springer Verlag,  (1980)  1---24.
\smallskip

\item{[B5]}Barlet~D., Majoration du volume des fibres g\'en\'eriques et forme g\'eom\'etrique du th\'eor\`eme d'applatissement. Seminaire P. Lelong-H.Skoda. Lecture Notes in Mathematics~822, Springer verlag, (1980), p.1---17
\smallskip

\item{ [B.M]}Barlet~D.,  Magnusson~J.,
Integration de classes de cohomologie m\'eromorphes et diviseur
d'incidence, Ann. sc. de l'E.N.S: S\'erie 4, tome 31 fasc. 6, (1998)
811---842.
\smallskip

\item{[B.V]}Barlet~D., Varouchas~J.  Fonctions holomorphes sur l'espace des cycles, Bulletin de la Soci\'et\'e Math\'ematique de France ~117, (1989) 329---341\smallskip

\item{[E]}Elzein~F., Complexe dualisant et applications \`a la classe fondamentale d'un cycle, Bull.Soc.Math.France Mem~58 (1978).
\smallskip

\item{[E1]}Elzein~F., Comparaison des r\'esidus de Grothendieck et
de Herrera. C.R.Acad.Sc. Paris,~278, (1974), p.863---866
\smallskip

\item{[C]}Cassa~A.,The Topology of the Space of Positive Analytic
Cycles. Annali di Math.Pura.Appl, tome 112,~1 (1977) p1---12.
\smallskip

\item{[E.G.A.4]}Grothendieck~A., Dieudonn\'e~J.,El\'ements de g\'eom\'etrie alg\'ebrique. Etude locale des sch\'emas et des morphismes de sch\'emas III. Inst.Hautes Etudes Sci.Publ.Math~28, (1966).
\smallskip

\item{[Fi]}Fischer~G.,Complex Analytic Geometry. Lect.Notes.Math~538, Berlin Heidelber New York, Spriner (1976).
\smallskip

\item{[Fu]}Fujiki~A., Closedeness of the Douady Spaces of compact
K\"ahler Spaces. Publ.Rims, Kyoto-Univ,tome 14, (1978) p1---52.
\smallskip

\item{[G1]}GRothendieck~A., Th\'eor\`emes de dualit\'e pour les les faisceaux alg\'ebriques coh\'erents. Seminaire Bourbaki ~149 (1957).
\smallskip

\item{[Go]}Golovin~V.D.,On the homology theory of analytic sheaves, Math.USSR.Izvestija~16,(2)  (1981),p239---260.
\smallskip

\item{[G-R]}Grauert~H, Riemenschneider~O., Verschwindungssätze für analytische kohomologiegruppen auf komplexen räumen. Invent. Math. ~11 (1970), p.263---292.\smallskip

\item{[G.R1]} Grauert~H., Remmert~R., Komplexe Ra\"ume. Math.Ann~136, (1958),p245---318.
\smallskip

\item{[G.R2]} Grauert~H., Remmert~R.,Coherent analytic sheaves.Grundl.Math.Wiss,Bd~265, Berlin heidelberg New York, Springer (1984).
\smallskip

\item{[Gri]}Griffiths~P., Variations on a theorem of Abel. Invent.Math~35(1976),p321---390.
\smallskip

\item{[Gr.H]}Griffiths~P, Harris~J., Principles of algebraic geometry. Wiley Classics Library. John Wiley and  Sons, Inc, New York, 1994.
\smallskip

\item{[G.R.P]}Remmert~R., Local Theory of Complex Spaces. Several Complex Variables VII, Encyclopaedia of Mathematical Sciences,~74, Springer-Verlag,  p.7---96
\smallskip

\item{[H]}Harvey~F.Reese.,Integral formulae connected by Dolbeault's isomorphism. Rice.Univ.Studies~56,(n2), (1970) p77---97.
\smallskip

\item{[Ha]}Hartshorne~R,.Residues and Duality. Lect. Notes in Mathematics~20, (1966), Springer-Verlag, Heidelberg.
\smallskip

\item{[H.L]}Herrera~M, Libermann~D,.Residues and principal values on complex spaces. Math. Ann~194, (1971),p. 259---294.
\smallskip

\item{[H.S]} H\"ubl~R, Sastry~P., Regular differential forms and relative duality. American Journal of Mathematics~115, (1993), p.749---787.
\smallskip

\item{[K1]}Kaddar~M.Classe fondamentale relative en cohomologie de Deligne et application. Math.Ann,~306,(1986),p285---322.
\smallskip

\item{[K2]}Kaddar~M., Int\'egration d'ordre sup\'erieure sur les cycles en g\'eom\'etrie analytique complexe, Ann. Sc. Norm. Pisa Cl.Sci(4), tome  29, (2000)
187---263.
\smallskip

\item{[Ke]}Kersken~M.,Der Residuencomplex in der lokalen algebraischen und analytischen Geometrie, Math.Ann~265,(1983) 423---455.
\smallskip

\item{[Ke1]}Kersken~ Regul\"are Diferentialformen. Manus. Math~46, (1984), p 1---25.
\smallskip

\item{[Ke2]}Kersken~ Some Applications of the trace mapping for differentials. Topis in Algebra. Banach Center Publications ~26 (2), (1990), p141---148.
\smallskip

\item{[K.W]}Kunz~E, Waldi~R.,  Regular differential forms, Contemporary Mathematics~79, (1988), Amer.Math.Soc., Providence.
\smallskip

\item{[Ki]}Kiehl~R., \"Aquivalenzerelationen in analytischen Ra\"umen. Math.Zeit. (1968), p.1---20
\smallskip

\item{[Kl]}Kleiman, Steven~L., Relative duality for quasicoherent sheaves, Compositio.Math.~41 (1980),(1) 39---60.
\smallskip

\item{[Ku1]}Kunz~E., Holomorphe Differentialformen auf Algebraischen Variet\"aten mit Singularit\"aten I.  Manuscripta.Math~15,(19750, p.91---108.
\smallskip

\item{[Ku2]}Kunz~E., Holomorphe Differentialformen auf Algebraischen Variet\"aten mit Singularit\"aten I.  Abh.Math.Sem.Univ.Hamburg~47,(1978), p.42---70.
\smallskip

\item{[Ku3]}Kunz~E., Residuen von  Differentialformen auf Cohen Macaulay  Variet\"aten mit Singularit\"aten I.  Math. Zeitschrift~152,(1977), p.165---189.
\smallskip

\item{[L]}Lipman~J.,  Dualizing  sheaves, differentials and residues on algebraic varieties,  Asterisque~117, (1984).
\smallskip

\item{[L.J]}Lejeune-Jalabert M., Le th\'eor\`eme ``$AF+BG$'' de Max Noether. Seminar on Singularities (Paris, 1976/1977),  pp. 97--138, Publ. Math. Univ. Paris VII, 7, Univ. Paris VII, Paris, 1980.
\smallskip

\item{[L.S]}Lipman~J, Sastry~P., Residues and duality for Cousin complexes, Preprint.
\smallskip

\item{[Lo]}Lojasiewicz~S.,Introduction to Complex Analytic Geometry, Birkh\"auser, Basel,(1991).
\smallskip

\item{[Ma]}Mathieu~D., Universal reparametrizationof family of
analytic cycles: a new approach to meromorphicc equivalence
relations. Analles de l'Inst.Fourier~50 (44),(2000) p1115---1189.
\smallskip

\item{[Mo]}Mochizuki~N., Quasi-normal analytic spaces, Proc.Japan Acad.,~48, (1972), p.181---185.
\smallskip

\item{[Par]}Parusinski~A., Constructibility of the set of points where a complex analytic morphism is open. Proc.Amer.Math.Society.~117(1), (1993),p205---211.
\smallskip

\item{[R]}Reiffen~H.J., Riemmansche Hebbartkeitss\"atze f\"r Kohomologieklassen mit kompactem Tr\"ager, Math.Ann~164, (1966) 272---279.
\smallskip

\item{[R.R]}Ramis~J.P, Ruget~G., Complexe dualisant et th\'eor\`emes de dualit\'e en g\'eom\'etrie analytique complexe. Publ. Math de l'IHES.,~38, (1970), p77---91.
\smallskip

\item{[RRV]}Ramis~J.P, Ruget~G, Verdier~J.L., Dualit\'e relative en g\'eom\'etrie analytique complexe, Invent.Math ~13 (1971) 261---283.
\smallskip

\item{[Re]}Remmert~R., Holomorphe und meromorphe Abbildung komplexer Ra\"ume. Math.Ann.~133, (1957),p 286---319.
\smallskip

\item{[Si]}Siebert~B.,Fiber cycles of holomorphic maps I-II, Math.Ann~296 et 300,(1993,1994),p269----283, p243---271.
\smallskip

\item{[V]}Verdier~J.L. Base change for twisted inverse image of coherent sheaves, Internat.Colloq, Tata Inst.Fund.Res, Bombay,(1968) 393---408.
\smallskip

\item{[Y]}Yekutieli~A., Smooth formal embeddings and the residue complex. Canad.J.Math.~50, (1998),(4), p.863---896.\bigskip\bigskip
M. Kaddar,\par\noindent
Institut Elie Cartan, UMR 7502\par\noindent
Universit\'e Nancy 1, BP 239\par\noindent
Vandoeuvre-l\`es-Nancy Cedex. France.\smallskip\noindent
e-mail: kaddar@iecn.u-nancy.fr

\end
 pour, au
moins, trois  raisons fondamentales, \`a savoir :\par $\star$ on ne
dispose pas  d'un analogue du th\'eor\`eme de compactification de
Nagata,  \par $\star$ la cat\'egorie des faisceaux coh\'erents sur
un espace analytique complexe arbitraire n'a pas assez d'objets
projectifs et donc pas de r\'esolution libre globale pour un
faisceau coh\'erent donn\'e,\par $\star$ la th\'eorie de la
dualit\'e analytique relative, qui n\'ecessite de la topologie,
n'est `` praticable'' que dans le cas d'un morphisme propre gr\^ace
\`a [R.R.V]. Le cas g\'en\'eral, tel qu'il est \'enonc\'e dans
[R.R1],  reste  encore incomplet et par endroit insastisfaisant. Il
nous parait pr\'ef\'erable d'installer cette th\'eorie dans le cadre
des espaces bornologiques...Mais ceci est une autre histoire
!\smallskip\noindent A ceci, s'ajoutent, entre autres, les
probl\`emes s\'erieux caus\'es par les questions de prolongement de
la coh\'erence, la pr\'esence de singularit\'es essentielles, la non
coh\'erence des faisceaux  images directes sup\'erieures pour un
morphisme arbitraire...
On a les propri\'et\'es suivantes:\par\noindent
\Th{1}{}
1) pour tout $j\geq 0$, $\omega^{j}_{Z}$ est de profondeur au moins $2$
dans $Z$.\smallskip\noindent
2) Si $Z$ est une $V$- vari\'et\'e, ces faisceaux coincident avec ceux
introduits par
 Steenbrink [St].\smallskip\noindent
3) $\omega^{k}_{Z}\simeq{\cal H}om_{{\cal O}_{Z}}(\Omega^{m-k}_{Z},
\omega^{m}_{Z})$ \smallskip\noindent
4)  $\omega_{Z}^{m}$ est le faisceau de Grothendieck, dualisant au sens
de Kunz  .\smallskip\noindent
5) Si $Z$ est un espace de Cohen Macaulay, le complexe $
\omega^{m}_{Z}[m]$ est quasi
 isomorphe au complexe dualisant de la g\'eom\'etrie analytique construit
par
 Ramis et Ruget .\smallskip\rm
Nous renvoyons \`a [B3]  pour la preuve des  propri\'et\'es 1 et 3, \`a
[St] pour 2 , \`a [L] pour 4 et  \`a [R.R] pour 5.\bigskip\noindentNous renvoyons le
lecteur au texte de Lipman [Li], pour un historique int\'eressant et
instructif, m\^eme si seul l'aspect alg\'ebrique est consid\'er\'e.
Notons que la plupart des  notions \'evoqu\'ees trouvent leurs
analogues ou s'\'etendent  au cadre de la g\'eom\'etrie analytique
complexe. Toutefois, on se gardera de croire que ce n'est qu'un
exercice de style ! Dans ce type de question, le cadre analytique
complexe diff\`ere consid\'erablement du cadre alg\'ebrique mais nous aurons tout le loisir de revenir sur le sujet dans [KIII]
Comme, par d\'efinition et par
analogie avec la construction de Kunz, on a $ \ds{f_{*}\omega_{f} =
{\cal H}om_{{\cal O}_{U}}(f_{*}{\cal O}_{X}, \Omega^{n}_{U})}$ pour
toute param\'etrisation locale $f:X\rightarrow U$, on en d\'eduit
facilement que  $\tilde{\omega}^{n}_{X}$ coincide avec le faisceau
dualisant de Grothendieck $\omega^{n}_{X}$Dans le cadre alg\'ebrique, leur \'etude
a \'et\'e entreprise, principalement,  par Kunz, Waldi et Platte,
par  Elzein dans [E] sous une forme "d\'eguis\'ee", et dans le cadre
de la g\'eom\'etrie analytique complexe dans ceux de Kersken [Ke] et
Barlet [B2].\smallskip\noindent{\bf 3.1.3.
Propri\'et\'e de la trace et formes m\'eromorphes
r\'eguli\`eres.}\smallskip\noindent
L'etude d'un tel faisceau a, d'abord, \'et\'e men\'e, dans le cadre alg\'ebrique, par Kunz dans une s\'erie d'articles ou il montre qu'il est dualisant et coincide, de ce fait, avec le faisceau de Grothendieck. Suit une g\'en\'eralisation \`a la situation relative donn\'ee par un morphisme {\bf plat} (cf [K.W]). Le cadre analytique traitant la situation absolue est presque une transcription de la construction de Kunz dans le langage des alg\`ebres analytiques locales. On peut citer Platte [Pl] et surtout Kersken [K1], [K2] pour l'\'elaboration de cette approche fid\`ele \`a celle de Kunz mais qui pr\'esente tout de m\^eme une complication technique non n\'egligeable, par rapport au cadre alg\'ebrique,  si $X$ n'est pas r\'eduit. Cependant, l'aspect ``dualisant `` d'un tel faisceau n'a jamais  \'et\'e abord\'e auparavant (du moins \`a notre connaissance). Mais il est facile de combler cette lacune en utilisant l'important travail de  Andr\'eotti-Kaas-Golovin. Le lecteur fera ais\'ement le lien avec  leurs  faisceaux dualisant.\par\noindent
Signalons, par ailleurs, que les calculs de Kunz-Waldi montrent l'existence, en tout degr\'e $r$, d'un faisceau de $k$-formes m\'eromorphes $\tilde{\omega}^{r}_{X}$ coincidant, bien s\^ur, avec $\tilde{\omega}^{n}_{X}$ en degr\'e maximal et caract\'eris\'es par la propri\'et\'e de la trace absolue (cf {\bf(1.0.3.1.0)}. Ils sont d\'efinis comme \'etant les sections $k$-m\'eromorphes dont le cup-produit par n'importe quelle $(n-r)$-forme holomorphe sur $X$ est une section du faisceau  $\tilde{\omega}^{n}_{X}$. Dans ce cas particulier o\`u $X$ est r\'eduit et de dimension pure, ils apparaissent naturellement, en g\'eom\'etrie analytique complexe,  comme les dualis\'es de Andr\'eotti-Kaas-Golovin des faisceaux $\Omega^{n-r}_{X}$ et ainsi $\tilde{\omega}^{r}_{X}=\omega^{k}_{X}={\cal H}om(\Omega^{n-r}_{X}, \omega^{n}_{X})={\cal
D}^{n}(\Omega^{n-r}_{X})$.\par\noindent
La donn\'ee d'une param\'etrisation locale $f:X\rightarrow U$ en un point $x$ de $X$, donne  le morphisme trace
trace ${\cal T}^{n}_{f}:f_{*}\tilde{\omega}^{n}_{X}\rightarrow
\Omega^{n}_{U}$, duquel r\'esulte l'isomorphisme
$\ds{f_{*}{\tilde\omega}^{n}_{X}\simeq {\cal H}om_{{\cal
O}_{X}}(f_{*}{\cal O}_{X}, \Omega^{n}_{U})}$ et par suite  le
morphisme $\ds{{\cal T}^{r}_{f}: f_{*}{\cal
D}^{n}(\Omega^{n-r}_{X})\rightarrow \Omega^{r}_{U}}$ qui est
\'evidemment un v\'eritable morphisme trace.\smallskip\noindent
Bien qu'apparente  en filigrane, il existe une description de ces faisceaux au moyen de la classe fondamentale.  En anticipant un
peu sur cette notion qui sera vu dans le
prochain paragraphe, l'approche de Barlet dans [B2]  consiste
simplement \`a utiliser les propri\'et\'es de cette classe. En effet,
consid\'erant un plongement local de $X$ dans une vari\'et\'e lisse
$Z$ et notant $j$ l'inclusion naturelle de la partie lisse de $X$
dans $X$, la classe fondamentale de $X$ fournit un morphisme {\it
cup-produit} $i_{*}\Omega^{k}_{X}\rightarrow {\cal E}xt^{p}({\cal
O}_{X}, \Omega^{p+k}_{Z}))$ dont on d\'eduit  un diagramme commutatif
$$\xymatrix{j_{*}j^{*}(\Omega^{k}_{X})\ar[rr]^{\partial}\ar[rd]_{\tilde{\partial_{n}}}&&{\cal H}^{1}_{\Sigma}(\Omega^{k}_{X})\ar[ld]^{{\cal H}^{1}_{\Sigma}({\cal C}_{X})}\\
&{\cal H}^{1}_{\Sigma}({\cal E}xt^{p}({\cal O}_{X},
\Omega^{p+k}_{Z}))&}$$ On pose, alors, $\omega^{k}_{X}:= {\rm
Ker}\,\tilde{\partial}_{k}$.\par\noindent
 Les propri\'et\'es
intrins\`eques de la classe fondamentale, en particulier son
ind\'ependance vis-\`a-vis du plongement choisi, montre que ce
faisceau est l'incarnation locale d'un faisceau intrins\`eque sur
$X$. Par ailleurs,on a, essentiellement par d\'efinition,
$$\omega^{k}_{X}\simeq i^{*}{\cal E}xt^{p}(i_{*}\Omega^{n-k}_{X},
 \Omega^{n+p}_{Z})$$
et par cons\'equent $\ds{\omega^{k}_{X} = {\cal
D}^{n}(\Omega^{n-k}_{X})}$.\smallskip\noindent
 Le faisceau coh\'erent $\omega^{n}_{X}:={\cal H}^{-n}({\cal
D}^{\bullet}_{X})$ jouit de propri\'et\'es remarquables. En effet,
il est dualisant au sens de la g\'eom\'etrie analytique complexe
puisque dual du faisceau structural ${\cal O}_{X}$ selon [A.K] ou
[Go] et, si  $X$ est de dimension pure $n$ et  muni d'un plongement
$\sigma$ de codimension $p$ dans une vari\'et\'e analytique complexe
$Z$ de dimension $n+p$, il s'identifie (en vertu des propri\'et\'es
des complexes dualisants) au faisceau  ${\cal E}xt^{p}({\cal O}_{X},
\Omega^{n+p}_{Z})$ s'identifiant  au  faisceau dualisant de Grothendieck
$$\vbox{\boxit{5pt}{${\rm Dim}{\cal B} = {\rm Dim}{\cal A} + {\rm Dim}{{\cal B}/{\goth m}{\cal B}}$}}$$
  Ainsi, \`a tout morphisme
plat, \`a fibres de dimension pure $n$, d'espaces analytiques
complexes ( r\'eduits) $ g:Y\rightarrow T$, est canoniquement
associ\'e un unique faisceau (\`a isomorphisme canonique pr\`es)
${\cal O}_{X}$-coh\'erent  $\tilde{\omega}^{n}_{Y/T}$ enti\`erement
caract\'eris\'e par la propri\'et\'e de la trace
relative.\par\noindent Une fa\c con directe pour nous a \'et\'e de
v\'erifier point par point les d\'etails de la construction de
Kersken pour nous appercevoir qu'elle n'est pas le privil\`ege de la
platitude mais de la platitude g\'eom\'etrique analytique. Mais il
est moins douloureux de construire un tel faisceau pour un morphisme
analytiquement  g\'eom\'etriquement plat gr\^ace aux diff\'erentes
caract\'erisations que l'on donne de ces morphismes.
 En d'autre termes, on peut \'enoncer le \smallskip\noindent
 Cor{3.1}{} Soient  $n$ un entier naturel, $X$ et $S$ des espaces
complexes r\'eduits avec $S$ de dimension localement pure. Soit
$\pi:X\rightarrow S$ un morphisme $n$-analytiquement
g\'eom\'etriquement plat.  Soit $X_{0}$ (resp. $S_{0}$) l' ouvert
dense de $X$ (resp. $S$) sur lequel la restriction $\pi_{0}$ de
$\pi$ est plate sur $S_{0}$. Alors, pour tout point $x$ de $X_{0}$,
on a un isomorphisme canonique  $$\omega^{n}_{\pi,x}=
\tilde{\omega}^{n}_{X_{0}/S_{0},x}$$
{\bf 0.2. Remarques.}\par\noindent {\bf(i)} Il est facile de
d\'eduire de ce qui pr\'ec\`ede que la platitude g\'eom\'etrique est
aussi \'equivalente au fait que  le faisceau ${\cal
O}_{X}$-coh\'erent $\omega^{0}_{\pi}$ soit  stable par changement de
base quelconque entre espaces
 complexes r\'eduits,  muni d'un morphisme canonique ${\cal O}_{X}\rightarrow \omega^{0}_{\pi}$ et v\'erifiant {\it la propri\'et\'e de la trace relative}.\par\noindent
$\bullet$ Si  $X_{1}$ est  l'ensemble des points de $X$ en lesquels
$\pi$ est g\'eom\'etriquement plat, $S_{1}$ son image dans $S$,
$\pi_{1}$ la restriction de $\pi$ \`a $X_{1}$ et  ${\cal
C}_{\pi_{1}}:\Omega^{n}_{X_{1}/S_{1}}\rightarrow
\omega^{n}_{\pi_{1}}$ le morphisme classe fondamentale donn\'e par
le th\'eor\`eme 3, alors,  $\pi$ est g\'eom\'etriquement plat si et
seulement si  le morphisme ${\cal C}_{\pi_{1}}$ se prolonge sur $X$
tout entier avec les m\^emes propri\'et\'es.\par\noindent
\smallskip\smallskip\noindent {\bf 0.3.  Conventions et
notations.}\par\noindent {\bf(i)} Sauf mention expresse du
contraire, les espaces analytiques complexes consid\'er\'es  sont
r\'eduits et de dimension localement  finie.\par\noindent
{\bf(ii)}Si $n$ est un entier naturel donn\'e, on notera ${\cal
E}(S,n)$ (resp. ${\cal G}(S,n)$) l'ensemble des morphismes
$\pi:X\rightarrow S$ $n$-\'equidimensionnels et ouverts  avec $X$
d\'enombrable \`a l'infini (l'ensemble des  morphismes de ${\cal
E}(S,n)$ qui sont g\'eom\'etriquement plats). On dira que $X$ est
$S$-g\'eom\'etriquement plat ou que $\pi$ est g\'eom\'etriquement
plat.\par\noindent {\bf(iii)} Il nous arrivera de dire ``changement
de base  r\'eduit'' pour signifier ``changement de base entre
espaces complexes r\'eduits''.\par\noindent {\bf(iv)} Une
pram\'etrisation locale d'un espace analytique complexe de dimension
pure $n$ en un point $x$ de $X$ est la donn\'ee d'un voisinage
ouvert $U$ de $x$ dans  $X$, d'un polydisque ouvert $V$ de ${\Bbb
C}^{n}$ et d'un morphisme fini et surjectif  $f:U\rightarrow V$
(dont l'interpr\'etation en termes d'alg\`ebres analytiques est
ais\'ee). La variante relative s'\'enonce de fa\c con similaire et
on dira, par abus de langage,  qu'un morphisme fini et surjectif
$f:X\rightarrow Y$ est une pram\'etrisation locale sous-entendant
par l\`a une donn\'e ``germifi\'ee'' comme pr\'ec\'edemment
d\'efinie. Sans entrer  dans les d\'etails de la terminologie disons que   cette classe de
morphismes r\'epond \`a la question suivante:\par\noindent
\smallskip\noindent
 \vrule height 66pt depth 0pt width 3pt$\,\,\,${\vbox{Soit
$\pi:X\rightarrow S$ un morphisme d'espaces analytiques complexes
avec $S$ r\'eduit et ${\rm dim}(X)\geq {\rm dim}(S)$.
Est-il-possible de munir les fibres ensemblistes de $\pi$ (plus
pr\'ecisemment leurs composantes irr\'eductibles) de multiplicit\'es
convenables de sorte que la famille $(\pi^{-1}(s))_{s\in S}$
devienne (ou soit induite par ) une famille analytique de cycles
  au sens de [B1]?}}
Un peu d'intuition g\'eom\'etrique montre que sans certaines
conditions sur le morphisme $\pi$, il est illusoire d'esperer une
issue positive \`a ce type de question. Apr\`es avoir introduit les
notions essentielles \`a la compr\'ehension du probl\`eme, nous
montrerons sur des exemples simples comment apparaissent les
conditions n\'ecessaires que doit v\'erifier un morphisme
susceptible de r\'epondre au probl\`eme.\smallskip\noindent
 La construction de cette derni\`ere se fait, ici, au niveau des faisceaux
  bien qu'elle \'emane  de la th\'eorie de la dualit\'e analytique relative qui  sera  abord\'e ult\'erieurement dans un papier  consacr\'e \`a la notion
de paire dualisante analytique (notion introduite par Kleiman ([Kl])
dans le cadre alg\'ebrique sous certaines conditions). Remarquons
que si $\pi:X\rightarrow S$ (\`a fibres de dimension pure $n$)
appartient \`a une classe de morphismes d\'efinis  dans le cadre
alg\'ebrique ou analytique complexe pour laquelle le foncteur {\it
image inverse extraordinaire} $\pi^{!}$ est bien d\'efini et
v\'erifie les propri\'et\'es d'usage alors $\omega^{n}_{\pi} :={\cal
H}^{-n}(\pi^{!}({\cal O}_{S}))$. Comme nous le constaterons dans le
\S 3, ce sera le cas pour un morphisme propre d'espaces analytiques
complexes d\'enombrables \`a l'infini de dimension finie avec des
fibres de dimension born\'ee par un entier $n$. Supposons  que la
restriction de $\pi$ \`a  $Y$ que l'on note $\pi_{Y}$ est encore
universellement \'equidimensionnelle  sur $S$ de fibres de dimension
au plus $n-2$.\smallskip\noindent Le probl\`eme \'etant de nature
locale sur $X$ au voisinage de $Y$, on peut supposer donner une
installation locale
$$\xymatrix{Y\ar@^{{(}->}[r]\ar[rrd]_{\pi_{Y}}&X\ar@^{{(}->}[r]\ar[rd]^{\pi}&S\times
U\times B\ar[d]\\
&&S}$$ Alors, on a
$${\cal H}^{1}_{Y}(\omega^{n}_{\pi})= {\cal H}^{1}_{Y}({\cal
H}^{-n}(\pi^{!}{\cal O}_{S})={\cal H}^{-n+1}(\pi_{Y}^{!}{\cal
O}_{S})$$ en donnant un sens \`a ces formules de dualit\'e gr\^ace
au diagramme pr\'ec\'edent.\par\noindent Mais les  fibres de
$\pi_{Y}$  \'etant  de dimension au plus $n-2$, on a
$${\cal H}^{j}({\pi_{Y}}^{!}{\cal O}_{S})=0\,\,\forall\,j\leq -n+2$$
Let $X$, $S$, $n$ be, respectively, complex
spaces with
 $X$ countable at infinity, $S$ locally finite pure dimensional reduced and
 $n$ an integer. Let $\pi:X\rightarrow S$ be an open  morphism with constant
 pure dimensional fibers (we call such morphism
  {\it universally-$n$-equidimensional}). If there is a cycle $\goth{X}$
  of $X\times S$ such that, fiberwise, coincide set-theorically
  with the fibers  of $\pi$ and endowed this with a good
  multiplicities in  such a  way that $(\pi^{-1}(s))_{s\in S}$ becomes
  a local analytic (resp. continuous) family of cycles in the sense of [B.M],
  we say that $\pi$ is {\it analytically geometrically flat} (resp. {\it continuously geometrically flat})
   according to the {\it weight} $\goth{X}$.}}}$$\rm
  The purpose of this article is to  characterize the set of analytically  geometrically flat maps as the biggest class, contained in the set of universally equidimensional morphism,  which members  admit a
weighted relative fundamental  class morphism satisfies many nic
e functorial properties.
 In particular, for a {\it flat} or {\it  finite Tor-dimension} morphism, we obtain the analog (in this setting) of the  relative fundamental class of [E] or [E.A].  At first, we show that any universally-$n$-equidimensional  $\pi$  is endowed with a coherent sheaf $\omega^{n}_{\pi}$ which coincide, on the regular part of $\pi$,  with the top degre relative holomorphic forms and for $\pi$ proper with a dualizing sheaf. Many others properties of this sheaf are giving in  theorem 1. Then we deduce this important corollary :\par\noindent
 {\it an  universally-$n$-equidimensional morphism $\pi$ between reduced
 complex spaces of locally finite dimension is analytically geometrically flat if
 and only if $\omega^{n}_{\pi}$ is the  Kunz-Waldi sheaf of regular meromorphic
  relative forms  compatible with arbitrary  base change between reduced
   complex spaces}, for which we have the analog  algebraic statement :\par\noindent
 {\it  Let $X$ and  $S$ be reduced (or more generally with no embedded points)
  locally noetherian schemes of finite Krull dimension  with $S$ excellent
  on field  ${\rm k}$ of  characteristic 0 . Let $\pi:X\rightarrow S$ be
  a  finite type, universally open with $n$- pure dimensional fibers, generically
   smooth morphism. Then the Kunz- Waldi relative sheaf of meromorphic regular
   forms, $\tilde{\omega}^{n}_{X/S}$,  is  compatible with any base change
   between ${\rm k}$-schemes like  $S$ and is dualizing if and only if $\pi$
   define an  algebraic family  of  $n$-cycles parametrized by
   $S$.}}\rm
\smallskip\noindent
In the present paper [KI], we give all definition, conceptual materials with
many examples or counter examples and some results  we need for understand better the proof of our results that we give in the next paper [KII]. In [KIII] we examine, in this setting of complex analytic geometry,  the concept of {\it dualizing pairs} which are introduced by Kleinman ([Kl]) in the algebraic setting.In the present paper [KI], we give all definition, conceptual materials with
many examples or counter examples and some results  we need for understand better the proof of our results that we give in the next paper [KII]. In [KIII] we examine, in this setting of complex analytic geometry,  the concept of {\it dualizing pairs} which are introduced by Kleinman ([Kl]) in the algebraic setting.
 le \Cor{1}{} Soit
$\pi:X\rightarrow S$ un morphisme $n$-plat d'espaces
analytiques complexes avec $S$ r\'eduit de dimension pure localement
fini. Alors, il existe un faisceau analytique (unique \`a
isomorphisme canonique pr\`es) $\omega^{n}_{X/S}$ v\'erifiant
:\par\noindent {\bf(i)} il est ${\cal O}_{X}$-coh\'erent et de
profondeur au moins deux fibre par fibre sur $S$,\par \noindent
{\bf(ii)} si $\pi$ est propre, $\omega^{n}_{X/S}= {\cal
H}^{-n}(\pi^{!}({\cal O}_{S}))$ \par\noindent {\bf(iii)} la famille
de faisceaux coh\'erents   $\omega^{\bullet}_{X/S}={\cal
H}om(\Omega^{n-\bullet}_{X/S}, \omega^{n}_{X/S})$ est munie d'une
diff\'erentielle non triviale ${\rm D}$, faisant de
$(\omega^{\bullet}_{X/S}, {\rm D})$  un complexe diff\'erentiel de
$(\Omega^{\bullet}_{X/S}, d_{X/S})$-modules.\par\noindent {\bf(iv)}
si $X$ est r\'eduit,  $\omega^{n}_{X/S}$ est le faisceau
$\tilde{\omega}^{n}_{X/S}$ des formes m\'eromorphes r\'eguli\`eres
caract\'eris\'e par la propri\'et\'e de la trace relative au sens de
Kunz-Waldi [K.W].\rm\smallskip\noindent